%% file: QFrob2.tex
\documentclass{amsart}
\usepackage[dvipsnames]{xcolor}

\usepackage{amssymb}
\usepackage{tikz}
\usetikzlibrary{arrows}
\usetikzlibrary{fit}

\usepackage{float}

\RequirePackage{color}
\definecolor{myred}{rgb}{0.75,0,0}
\definecolor{mygreen}{rgb}{0,0.5,0}
\definecolor{myblue}{rgb}{0,0,0.65}

\RequirePackage{ifpdf}
\ifpdf
  \IfFileExists{pdfsync.sty}{\RequirePackage{pdfsync}}{}
  \RequirePackage[pdftex,
   colorlinks = true,
   urlcolor = myblue, % \href{...}{...} external (URL)
   citecolor = mygreen, % \cite{...}
   linkcolor = myred, % \ref{...} and \pageref{...}
   pagebackref,
   bookmarksopen=true]{hyperref}
\else
  \RequirePackage[hypertex]{hyperref}
\fi

\RequirePackage{ae, aecompl, aeguill} % to have pretty pdf

\def\un{\underline}

\newtheorem{thm}{Theorem}[section]
\newtheorem{lem}[thm]{Lemma}
\newtheorem{prop}[thm]{Proposition}
\newtheorem{cor}[thm]{Corollary}

\theoremstyle{definition}
\newtheorem{defn}[thm]{Definition}
\newtheorem{notation}[thm]{Notation}
\newtheorem{ex}[thm]{Example}

\theoremstyle{remark}
\newtheorem{rem}[thm]{Remark}

\def\ZZ{\mathbb{Z}}
\def\QQ{\mathbb{Q}}
\def\CC{\mathbb{C}}

\DeclareMathOperator{\rt}{root}

\DeclareMathOperator{\Sym}{Sym}
\DeclareMathOperator{\End}{End}

\DeclareMathOperator{\aff}{aff}
\DeclareMathOperator{\fin}{fin}

\def\ze{\zeta}
\def\La{\Lambda}

\def\si{\sigma}
\def\Om{\Omega}

\newcommand{\pa}{\partial}
\newcommand{\NC}{NC}

\newcommand{\stair}{P}

\DeclareMathOperator{\term}{magicterm}
\DeclareMathOperator{\magic}{magic}
\DeclareMathOperator{\foobar}{exp}
\DeclareMathOperator{\bottom}{bot}
\DeclareMathOperator{\blahblah}{foobar}
\newcommand{\whyme}{\lambda_{49}}

\newcommand{\thesign}{\mu}

\title[]{A closed formula in the deformed affine nilHecke algebra}

\author[]{Ben Elias}
\address{University of Oregon.}
\email{belias@uoregon.edu}%\thanks{Blah blah}

\author[]{Daniel Juteau}
\address{LAMFA, Universit\'{e} de Picardie Jules Verne.}
\email{daniel.juteau@u-picardie.fr}%\thanks{Blah blah}

\author[]{Benjamin Young}
\address{University of Oregon.}
\email{bjy@uoregon.edu}%\thanks{Blah blah}

\begin{document}

\begin{abstract} There is a $q$-deformation of the reflection representation of the affine symmetric group (of type $\tilde{A}_{n-1}$), which arises in the quantum geometric Satake
equivalence, and in the study of the complex reflection groups $G(m,m,n)$. Demazure operators (often called divided difference operators) act on the polynomial ring of this deformed representation. When $n=3$ we prove an explicit closed formula for the scalar one obtains when applying a degree $-\ell$ Demazure operator to a monomial of degree $\ell$. We also prove a simpler formula for the scalar obtained after specializing $q$ to a root of unity. \end{abstract}

\maketitle

\tableofcontents

\input{IntroductionTwo.tex}

\input{Magic.tex}

\input{Preliminaries.tex}
\input{Symmetry.tex}
\input{NastyFour.tex}

\input{FormulaROUpapertwo.tex}

\bibliographystyle{plain}
\bibliography{mastercopy}

\end{document}

%% file: IntroductionTwo.tex
%!TEX root = QFrob2.tex

%%%%%%%%%%%%%%%%%%%%%%%%%%%%%%%%%%%%%%%%%%%%%%%%%%%%%%%%%
%========================================================
\section{Introduction} \label{sec:intro}
%========================================================
%%%%%%%%%%%%%%%%%%%%%%%%%%%%%%%%%%%%%%%%%%%%%%%%%%%%%%%%%

The affine Weyl group $W_{\aff}$ can be viewed as the semidirect product of the finite Weyl group $W_{\fin}$ with its root lattice $\La_{\rt}$. In \cite{EQuantumI}, the first author
introduced a $q$-deformation of the affine Cartan matrix in type $\tilde{A}_{n-1}$, leading to a $q$-deformation of the reflection representation of $W_{\aff}$. When $q$ is specialized
to a primitive $2m$-th root of unity, this representation factors through the quotient $W_m := W_{\aff}/(m \cdot \La_{\rt})$, a complex reflection group known as $G(m,m,n)$.

Let $R$ denote the polynomial ring of this $q$-deformed reflection representation. For each subset $I$ of the simple reflections of $W_{\aff}$, one has a parabolic subgroup $W_I$,
and a subring $R^I := R^{W_I}$ of invariant polynomials. 
When $W_I$ is finite, the inclusion $R^I \subset R$ is a Frobenius extension, implying that induction and restriction are biadjoint functors (up
to a grading shift). One can build a $2$-category of bimodules by iteratively tensoring induction and restriction bimodules, and taking direct summands; these are known as singular
Soergel bimodules. The primary reason to study this $q$-deformation is that it leads to a $q$-deformation of the geometric Satake equivalence, matching representations of the
quantum group of $\mathfrak{sl}_n$ with $q$-deformed singular Soergel bimodules. For more details, see \cite{EQuantumI}. The recent work \cite{EWKSBim} explains this $q$-deformation
as an example of a conjectural $K$-theoretic geometric Satake equivalence.

In \cite{EJY1}, we study the $q$-deformed reflection representation after $q$ is specialized to a root of unity, laying the groundwork for the study of singular Soergel bimodules for
$G(m,m,n)$. It is proven there that $R^{W_{\aff}} \subset R$ is a Frobenius extension when $q$ is a root of unity. Morever, when $n=3$ it is proven (using results from this paper) that
the Frobenius trace $J$ can be constructed within the algebra generated by divided difference operators. Henceforth we focus on the case $n=3$. Given a word $\un{w}$ of length $\ell$ in
the simple reflections of $W_{\aff}$, one obtains an operator $\pa_{\un{w}}$ of degree $-\ell$ on $R$, by composing divided difference operators. When $\ell = 3m$, $\pa_{\un{w}}$ agrees
with the Frobenius trace $J$ up to scalar. This scalar agrees with the evaluation of $\pa_{\un{w}}$ on a certain monomial $\stair = x_1^{2m} x_2^m$, as $J(\stair) = 1$. For many words,
$\pa_{\un{w}}(\stair) = 0$. Knowing for which words the scalar is invertible, and knowing the nonzero scalars precisely, will be crucial in future work which describes the $2$-category
of singular Soergel bimodules for $G(m,m,n)$ by generators and relations. For more detailed motivational material, including an explanation of why knowing the precise scalars is
important, see the introduction to \cite{EJY1}.

It is in the current paper, a technical companion to \cite{EJY1}, where we compute these scalars $\pa_{\un{w}}(\stair)$ and supply the final ingredient in the proofs of \cite{EJY1}.
Note that these scalars are Laurent polynomials in $q$ (or more precisely, in $q^{\frac{1}{3}}$). We discovered a lovely formula for these scalars $\pa_{\un{w}}(\stair)$ by extensive
computer experiments. Up to unit, the result is a quantum binomial coefficient (evaluated at the root of unity) times $m^2$.

However, proving our formula directly is not easy! We strongly recommend that the reader try out Example \ref{ex:dangitscomplex}, which illustrates the remarkable complexity of the
problem, hidden by the remarkable simplicity of the final answer. The operator $\pa_{\un{w}}$ is a composition of divided difference operators, so one might imagine that the formula
could be proven by applying divided difference operators one at a time, and keeping track of the result at each step along the way. While the final result is a scalar, the intermediate
steps yield polynomials of higher degree, and these polynomials are insanely complicated! We were unable to make any headway using this method.

Our formula is also not directly amenable to induction on $m$, as there is no obvious relationship between specializations to different roots of unity.

Instead, we were able to progress as follows. Given any word $\un{w}$ of length $\ell$ and any monomial $f$ of the same degree, we found a much more complicated but explicit formula
for the scalar $\pa_{\un{w}}(f)$ \emph{before} $q$ is specialized to a root of unity. Then, one can evaluate the formula at a root of unity and simplify the result (which is not a
straightforward task). The generic case of our formula for $\pa_{\un{w}}(f)$ is given in Theorem \ref{thm:mainnasty}. However, there are also many special cases depending on the
properties of $\un{w}$ and $f$. Remarks on the discovery of the formula can be found in the next section of the introduction.

Our formula for generic $q$ is amenable to inductive proof. Using properties of divided difference operators (e.g. the twisted Leibniz rule) we can show that the scalars we examine
satisfy certain recursive formulas. By proving that our formulas also satisfy these recursive formulas, we prove their correctness. The disadvantage of this approach is that one
must treat all possible words $\un{w}$ and polynomials $f$ in order to crank the recursion. This leads to a lot of case by case analysis. The advantage of this approach is that it
works.

% The advantage of this approach is that it works, but the disadvantage is that we are required to study all
%
% % Ultimately, we are able to prove this formula, but only by constructing a very large and obnoxiously technical hammer, in order to hit this very small nail.
%
%  The generic case of this formula is given in Theorem \ref{thm:mainnasty}. However, there are also many special cases depending on the
% properties of $\un{w}$ and $f$. The upshot is that this complicated formula is amenable to an inductive proof. (The downshot is that the proof has a lot of case by case
% analysis.)

We present and prove our formula in \S\ref{sec:closedformula}. Then in \S\ref{sec:closedformularou} we prove a much simpler, fairly nice formula for the specialization at a $2m$-th root of unity, when applied to the particular monomial $\stair$. Only the generic case of the formula yields a nonzero answer when applied to $\stair$ at a root of unity.
%, and up to unit the result is a quantum binomial coefficient (evaluated at the root of unity) times $m^2$.

\begin{rem} One might expect a simple formula to arise when $q$ is specialized to $1$ as well. Indeed, there is a simple formula: whenever the length of $\un{w}$ is at
least $4$, then the operator $\pa_{\un{w}}$ is zero, c.f. Remark \ref{rem:whenqis1}. This might seem surprising, since the nilHecke algebra in affine type is infinite-dimensional.
However, setting $q = 1$ does not recover the reflection representation of $W_{\aff}$; instead, it recovers the reflection representation of $W_{\fin}$, inflated to $W_{\aff}$ via
the quotient map $W_{\aff} \to W_{\fin}$ which kills translations. When $n=3$, $W_{\fin} = S_3$, and the longest element has length $3$. \end{rem}

%========================================================
\subsection{On the proof}
%========================================================

Having discussed the result and outlined the proof, we now wish to discuss which parts of the proof are difficult and interesting. To that end, let us describe roughly the formula
for $\pa_{\un{w}}(f)$ given in Theorem \ref{thm:mainnasty}. To find this formula, we first wrote computer code to compute $\pa_{\un{w}}(f)$, and then stared at countless examples. Our code can be found at~\cite{EJYcode}.

In \S\ref{sec:prelims} we explain why we need only consider the words $\un{w} = \un{w}(a,b,i)$ of
a certain form: a clockwise word of length $a+1$ followed by a counterclockwise word of length $b+1$; the two words overlap in {\color{red} one index}, so the overall length is $\ell =
a+b+1$. The final index is $i$. For example, $\un{w}(7,5,1) = (3,1,2,3,1,2,{\color{red} 3},2,1,3,2,1)$. Let $\pa_{(a,b,i)}(f)$ denote $\pa_{\un{w}(a,b,i)}(f)$. We also explain in \S\ref{sec:prelims} why we need only compute $\pa_{(a,b,i)}(f)$ when $f = x_1^k x_2^{\ell-k}$, where $\ell$ is the length of $\un{w}(a,b,i)$ and $0 \le k \le \ell$. For degree reasons, $\pa_{(a,b,i)}(f)$ is a scalar.

Up to a unit, $\pa_{(a,b,i)}(f)$ is equal to a product of three factors we denote $\gamma_r$ for $r = 1, 2, 3$. Our descriptions of $\gamma_r$ in this introduction are accurate only up to a sign, a power of $q$, and a power of $(q-q^{-1})$.

We observed that the scalar $\pa_{(a,b,i)}(f)$ is often divisible by quantum factorial numbers. We found the largest factorials which always divide $\pa_{(a,b,i)}(f)$ for a fixed
triple $(a,b,i)$ as $f$ varies. This ends up being $\gamma_1(a,b) = [\alpha]_q! [\beta]_q!$, where $\alpha$ and $\beta$ are roughly half of $a$ and $b$ respectively. More precisely, $a = 2\alpha+1$ or $a = 2\alpha + 2$ depending on parity, and similarly for $\beta$.

We also observed that, as $k$ varies between $1$ and $\ell-1$, the zeroes of $\pa_{(a,b,i)}(x_1^k x_2^{\ell-k})$ are predictable. For example, suppose that $a$ is odd and $b$ is even, and both are positive. Then the result is zero if and only if
\begin{itemize}
	\item $i=1$ and $k = \frac{\ell}{2}$, or
	\item $i=2$ and $k = \ell-1$, or
	\item $i=3$ and $k = 1$.
\end{itemize}
In our formulas, these zeroes are explained by the fact that
\begin{itemize}
	\item if $i=1$ then $[\frac{\ell}{2}-k]_q$ divides $\pa_{(a,b,1)}(x_1^k x_2^{\ell-k})$ for all $k$,
	\item if $i=2$ then $[\ell-1-k]_q$ divides $\pa_{(a,b,2)}(x_1^k x_2^{\ell-k})$ for all $k$,
	\item if $i=3$ then $[k-1]_q$ divides $\pa_{(a,b,3)}(x_1^k x_2^{\ell-k})$ for all $k$.
\end{itemize}
For other parities of $a$ and $b$ the zeroes are in different locations, but regardless, all zeroes are explained by a multiplicative factor of this sort. These factors are what we denote $\gamma_2(a,b,i,k)$.

What remains when dividing $\pa_{(a,b,i)}(x_1^k x_2^{\ell-k})$ by $\gamma_1$ and $\gamma_2$ is what we call $\gamma_3(a,b,i,k)$. Let us introduce a function\footnote{Outside of the introduction, we use different notation for the inputs to $\magic$.} called $\magic$, a sum of double quantum binomial coefficients.
\begin{equation} \magic(\alpha,\beta,k,\epsilon) = \sum_{j} {k-1 \brack \beta-j}_q {\alpha+\beta+1-k \brack j}_q q^{j(2k-3\alpha - 3 \beta - 6 - 2 \epsilon)}. \end{equation}
Finitely many terms in the sum are nonzero. For example, when $a$ is even and $b$ is odd, we have
\begin{equation}\label{gamma3intro} \gamma_3(a,b,i,k) = \magic(\alpha, \beta, k, 0). \end{equation}
In all other cases we can define $\gamma_3$ using $\magic$, see \eqref{gamma3} for details. The $\epsilon$ variable is a small adjustment to the power of $q^j$, allowing us to tweak the function to account for the parity of $a$ and $b$ and the value of $i$.

While this formula for $\gamma_3$ is relatively nice, discovering it was an act of stubbornness and sorcery and sheer serendipity. See Example \ref{magicexample} for an example of what this complicated Laurent polynomial looks like in practice.

Some readers may be familiar with the $q$-Chu-Vandermonde identity, which states the following for any $M, N \ge 0$ and $\beta \in \ZZ$ (both sides are zero for $\beta < 0$):
\begin{equation} \label{qchuvan} \sum_j {M \brack \beta-j}_q {N \brack j}_q q^{j(M+N)} = q^{N\beta} {M+N \brack \beta}_q. \end{equation}
This is a closed formula for a sum of double quantum binomial coefficients. Sadly, $\magic$ has the wrong power of $q^j$, and (except in special cases) one can not apply the $q$-Chu-Vandermonde identity, nor do we know of a closed form for $\magic$. One should think of $\magic$ as an alternate $q$-deformation of $\binom{\alpha+\beta}{\beta}$.

\begin{rem} The $q$-Chu-Vandermonde identity has a bijective proof, where both sides are a weighted count of subsets of size $\beta$ inside a larger set of size $M+N$. Our function
$\magic$ is also counting subsets of size $\beta$, but with different weights. The scalars we compute relate conjecturally, via quantum geometric Satake, to certain computations in
the $K$-theory of the affine Grassmannian. It would be interesting to find a matching combinatorial interpretation of both scalars. \end{rem}

One way to prove the $q$-Chu-Vandermonde identity uses generating functions. Recall the $q$-binomial theorem, which states (for $N \ge 0$) that
\begin{equation} \label{qbinomialthmintro}
\sum_{j=0}^N
{N \brack j}_q q^{-j(N-1)}t^j
=
\prod_{c=0}^{N-1}\left(1+q^{-2c}t\right).
\end{equation}
By setting $t = q^{K} x$ for an appropriate value of $K$, the $q$-binomial theorem implies that
\begin{equation} \sum_{\beta} q^{N \beta} {M+N \brack \beta}_q x^{\beta} = \prod_{\lambda \in X} (1 + q^{\lambda} x) \end{equation}
where $X$ is some particular set of integers having size $M+N$. More precisely, $X$ is a \emph{parity interval}, containing all the even (resp. odd) integers between two even (resp. odd) integers.  Meanwhile, by letting $u = \beta - j$, we can rewrite the double sum as follows:
\begin{equation} \sum_{\beta} \sum_j {M \brack \beta-j}_q {N \brack j}_q q^{j(M+N)} x^{\beta} = \left( \sum_{u} {M \brack u}_q x^u \right) \left( \sum_j {N \brack j}_q q^{j(M+N)} x^j \right). \end{equation}
Each sum on the right side can be rewritten using the $q$-binomial theorem as a product of terms $(1+q^{\lambda} x)$ as $\lambda$ ranges over a disjoint union of parity intervals $Y$ and $Z$ of sizes $M$ and $N$ respectively. To prove the $q$-Chu-Vandermonde identity, one need only observe that $Y \cup Z = X$.

Indeed, $\magic$ has a similar generating function, a product of terms $(1+q^{\lambda} x)$ as $\lambda$ ranges over a disjoint union of two parity intervals $Y \cup Z$. Unlike the $q$-Chu-Vandermonde setting, $Y \cup Z$ is not itself a parity interval, so one cannot ``undo'' the $q$-binomial process to get an easier description of any given coefficient.

\begin{rem} The generating function we actually consider is not the sum over $\beta$ of $\magic(\alpha,\beta,k,\epsilon) x^\beta$ as hinted above. Instead, we fix the sum $\alpha
+ \beta$ while letting $\beta$ vary. For details see \S\ref{ssec:magicgenerating}. \end{rem}

Nonetheless, we are able to use this generating function to great effect. In order to prove that our formula for $\pa_{(a,b,i)}(x_1^k x_2^{\ell-k})$ is correct, we must prove that
our formula satisfies certain symmetries, and that it obeys certain recursive formulas. These are best proven using the generating function. For example, in the case
governed by \eqref{gamma3intro}, one symmetry implies that $\gamma_3(a,b,i,k)$ and $\gamma_3(a,b,i,2\alpha+2\beta+4-k)$ agree up to a power of $q$. This is not at all obvious from
the formula \eqref{gamma3intro}, but it follows readily from the generating function.

The hardest part of the paper is the proof of a certain recursive formula in \S\ref{ssec:recursion1}. In the case where $b$ is even and $a$ is odd it reduces to the following statement (note that $\ell = a+b+1$ is even): when $\frac{\ell}{2} \le k < \ell$ we have
\begin{align} \nonumber (-1)^{k} (q^{-2k} - q^{-\ell})&  \magic(\alpha,\beta, k,0) q^{k(k - \beta - \ell+1)} \\ & = q^{-2 \ell + 2} \sum_{c = \ell-k}^{k-1} (-1)^c \magic(\alpha,\beta,c,-1) q^{c(c - \beta - \ell +3)}.
\end{align}
We are unaware of any similar formulas related to sums of double quantum binomial coefficients in the literature. Our proof goes via the generating functions but is still not easy, using a technique we like to call the \emph{double telescope}. We find an explicit formula for the partial sums on the right-hand side, explicitly find the ratio between successive partial sums, and thereby prove that the partial sums are equal to partial products in a telescoping product.  See \S\ref{ssec:magictelescope} for more details.

In \S\ref{sec:closedformularou} we have yet more fun manipulating these formulas when evaluated at a root of unity.

These sorts of manipulations form the most interesting part of our proofs. Sadly, in addition to the factors $\gamma_r$, there is also a power of $q$ and a sign to keep track of. The
power of $q$ appearing in $\pa_{(a,b,i)}(x_1^k x_2^{\ell-k})$ is a non-homogeneous degree $2$ polynomial in the variables $\{\alpha, \beta, k\}$, which also depends annoyingly on the
parity of $a$ and $b$ and the value of $i$. Keeping track of these powers of $q$ is a tedious bookkeeping exercise. For sanity, we have written the power of $q$ as a product of various
factors $\kappa_r$ and $\lambda_r$ which are easier to individually analyze. To separate the interesting manipulations from the tedious bookkeeping, we have placed most discussion of
the function $\magic$ in its own chapter \S\ref{sec:magic}.

For edge cases the formulas are simpler. For example, when $k = \ell$, the factors $\gamma_r$ are collectively replaced with a factor of $[\alpha+\beta+1]_q!$. Nonetheless, there are
many recursive formulas to check, and each must be confirmed for many cases; each of these confirmations requires its own tedious bookkeeping exercise. We have indeed performed each
exercise. For reasons of length we have omitted many of these confirmations, focusing on the justification that our formulas hold up to unit. We have
also made our computer code publically available, and ensured it is readable and well-commented \cite[\texttt{VerificationMagic.m}]{EJYcode}. Our MAGMA code is simple enough that one can run it using MAGMA's free online interface\footnote{Found here: \texttt{http://magma.maths.usyd.edu.au/calc/}}. When it comes to bookkeeping and arithmetic, perhaps computer verification of thousands of cases is more convincing than anything we could write.

%========================================================
\subsection{Organization of the paper and other comments}
%========================================================

In \S\ref{sec:prelims} we set our notation and recall a few results from \cite{EJY1}, in order to make this paper mathematically self-contained (though again, we rely on \cite{EJY1}
for further motivation). We also give some technical introductory remarks. In \S\ref{sec:symmetry} we discuss the many symmetries which these computations possess. We also discuss
various recursive formulae which uniquely pin down the scalars we seek, and outline the proof that our formulas are correct. While it may seem strange to outline the proof of the
theorem before stating the theorem itself, having a thorough discussion of symmetry does streamline the statement of the theorem significantly. In \S\ref{sec:closedformula} we state
our complicated formula, and in \S\ref{sec:closedformularou} we evaluate our formula at a root of unity. In \S\ref{sec:magic} we prove various properties of the function $\magic$
which are used in the proof. Most of the fun is in \S\ref{sec:magic} and \S\ref{sec:closedformularou}.

\begin{rem} To permit additional symmetry, we change variables from $q$ to $z$, where $z^3 = q^{-2}$. We use a more symmetric variant of the $q$-deformed reflection representation of \cite{EQuantumI}, a variant introduced in \cite{EJY1}. See \cite[\S 2]{EJY1} for a thorough discussion of this topic. \end{rem}

\textbf{Acknowledgments} The first author was supported by NSF grant DMS-2201387. The first and third authors appreciate the support given to their research group by NSF grant DMS-2039316.

%% file: Magic.tex
%!TEX root = QFrob2.tex

%%%%%%%%%%%%%%%%%%%%%%%%%%%%%%%%%%%%%%%%%%%%%%%%%%%%%%%%%
%========================================================
\section{A little bit of magic} \label{sec:magic}
%========================================================
%%%%%%%%%%%%%%%%%%%%%%%%%%%%%%%%%%%%%%%%%%%%%%%%%%%%%%%%%

Let us recall the definition of $\magic$ from the introduction and tweak it slightly, now using a variable $\nu$ instead of $\alpha$.

\begin{defn} Suppose $\nu, k \ge 1$. Let $\term(\nu,k,\beta,\epsilon,j)$ be the scalar
\begin{equation} \label{eq:termdef} \term(\nu,k,\beta,\epsilon,j) = {k-1 \brack \beta-j} {\nu-k-1 \brack j} q^{j(-3\nu + 2k - 2 \epsilon)}. \end{equation}
Let $\magic(\nu, k, \beta, \epsilon)$ be the scalar
\begin{equation} \label{eq:magicdef} \magic(\nu, k,\beta,\epsilon) = \sum_{j = 0}^{\beta} \term(\nu,k,\beta,\epsilon,j). \end{equation}
Note that $\term(\nu,k,\beta,\epsilon,j)$ vanishes for $j < 0$ or $j > \beta$, so we can sum over all $j \in \ZZ$ if desired. This simplifies certain arguments.
\end{defn}

To shed some light on how $\magic$ is used later, let us briefly discuss the factor $\gamma_3(a,b,i,k)$ in our formulas for $\pa_{(a,b,i)}(x_1^k x_2^{\ell-k})$. Letting $\alpha$ (resp. $\beta$) be roughly half of $a$ (resp. $b$) as in the introduction, and letting $\nu = \alpha+\beta+2$, we have
\begin{equation} \label{gamma3}
\gamma_3(a,b,i,k)  = \begin{cases}
	\magic(\nu,k,\beta,-1) & \text{ if $a$ and $b$ are odd} \\
	\magic(\nu,k,\beta,0) & \text{ if $a$ is even and $b$ is odd } \\
	\magic(\nu,k,\beta,0) & \text{ if $a$ is odd and $b$ is even and } i=1 \\ 
	\magic(\nu,k,\beta,-1) & \text{ if $a$ is odd and $b$ is even and } i=2 \\
	q^{\beta} \magic(\nu,k-1,\beta,-1) & \text{ if $a$ is odd and $b$ is even and } i=3 \\
	\magic(\nu,k,\beta,+1) & \text{ if $a$ and $b$ are even and } i=1 \\
	\magic(\nu,k,\beta,0) & \text{ if $a$ and $b$ are even and } i=2 \\
	q^{\beta} \magic(\nu,k-1,\beta,0) & \text{ if $a$ and $b$ are even and } i=3 	
	 \end{cases}
\end{equation}
In applications $\epsilon \in \{-1,0,1\}$, and we think of it as a small \emph{offset}. In most cases, $\epsilon$ is such that the length of the word in question, $\ell = a+b+1$, is equal to $2 \nu + \epsilon$. However, when $b$ is even and $i \in \{2,3\}$, the length is equal to $2 \nu + 1 + \epsilon$. The values of $k$ we plug in will range from $1$ to $\ell-1$.

{\bf We assume throughout this paper that} $\epsilon \in \{-1,0,1\}$.

The goal of this chapter is to develop methods to manipulate $\magic$, which we use in later chapters to prove properties of our formula for $\pa_{(a,b,i)}(x_1^k x_2^{\ell-k})$,
such as symmetries or recursive formulae. These symmetries and recursions will change the triple $(a,b,i)$, and hence will change the inputs to $\magic$ (especially the offset
$\epsilon$). The reason we demonstrated $\gamma_3$ now was so that one can expect formulas which mix the different versions of $\gamma_3$ above. We will not mention $\gamma_3$ any
further in this chapter.

Additional quantum number fun can be found in \S\ref{ssec:magicrou}, when $\magic$ is specialized at a root of unity.

\begin{ex} \label{magicexample} This example is not particularly illuminating, though we welcome the curious reader to explore $\magic$ using our code.
\begin{align}\nonumber \magic(8,4,3,0) & = q^{-48} + q^{-36} + 2q^{-34} + 3q^{-32} + 2q^{-30} + q^{-28} \\ & + q^{-20} + 2q^{-18} + 3q^{-16} + 2q^{-14} + q^{-12} + 1. \end{align}
\begin{align}\nonumber \magic(8,3,3,0) & = q^{-57} + q^{-55} + q^{-53} + q^{-51} + q^{-41} + 2q^{-39} + 3q^{-37} \\ & + 3q^{-35} + 2q^{-33} + q^{-31} + q^{-21} + q^{-19} + q^{-17} + q^{-15}. \end{align}
\end{ex}

%========================================================
\subsection{Special cases} \label{ssec:magicspecial}
%========================================================

\begin{lem} \label{lem:magicwhenbetazero} When $\beta = 0$ we have $\magic(\nu, k, 0, \epsilon) = 1$. When $\beta < 0$ we have $\magic(\nu,k,\beta,\epsilon) = 0$. \end{lem}

\begin{proof} Note that $\term(\nu,k,\beta,\epsilon,j)$ is zero unless both $j$ and $\beta-j$ are positive. When $\beta = 0$ the only nonzero term in the sum is $\term(\nu,k,0,\epsilon,0) = 1$. \end{proof}
	
\begin{lem} After specializing $q = 1$, $\magic(\nu, k, \beta, \epsilon) = \binom{\nu-2}{\beta}$. \end{lem}

\begin{proof} This is an immediate consequence of the Chu-Vandermonde identity. \end{proof}

In \eqref{qchuvan} we recalled the $q$-Chu-Vandermonde identity. By swapping $q$ and $q^{-1}$, there is a variant on this formula with factors $q^{-j(M+N)}$ and $q^{-N\beta}$
instead. If we set $M = k-1$ and $N = \nu-k-1$, then $M+N = \nu-2$ independent of $k$. The function $\magic$ looks like one side of the $q$-Chu-Vandermonde identity, except that $q^{j(-3 \nu + 2k - 2 \epsilon)}$ is not typically equal to $q^{\pm j (\nu - 2)}$. Sometimes, however, these powers of $q$ are equal. The following lemma treats these special cases. 

\begin{lem} \label{lem:magicchuvan} Fix $\beta, \nu > 0$. Then
\begin{equation} \label{magicchuvan} \magic(\nu,k,\beta,\epsilon) = q^{\pm (\nu-k-1) \beta} {\nu-2 \brack \beta} \quad \text{whenever } 2(k - \epsilon) = 3 \nu \pm (\nu - 2) . \end{equation}
\end{lem}

\begin{proof} This is straightforward using \eqref{qchuvan} or its variant discussed above. \end{proof}

Special cases covered by this lemma include when $k = 2 \nu$ and $\epsilon = +1$, or $k= \nu$ and $\epsilon = -1$, etcetera.

%========================================================
\subsection{Factorial manipulations} \label{ssec:magicmanipulation}
%========================================================

The following lemma can be viewed roughly as a recursive formula for $\magic$. It is used to prove \eqref{recursiveformula2}, part of the recursive formula for $\pa_{(a,b,i)}(x_1^k x_2^{\ell-k})$.

\begin{lem}\label{lem:magicmanip1} When $\epsilon \in \{-1,0\}$ We have
\begin{align} \label{magicmanip1} \nonumber [\beta] \magic(\nu,k,\beta,\epsilon) = & [k-1] \magic(\nu-1,k-1,\beta-1,\epsilon) \quad + \\ & q^{2k - 3\nu - \beta - 2 \epsilon + 1} [\nu-k-1] \magic(\nu-1,k,\beta-1,\epsilon+1). \end{align}
\end{lem}

\begin{proof}
Let us write all three versions of the $\magic$ function in terms of a similar-looking sum. Recall that $\magic = \sum_{j \in \ZZ} \term$. We have 
\begin{align} \label{foobar3} [k-1] \term(\nu-1,k-1,\beta-1,\epsilon,j) = \\ 
	\nonumber [k-1] {k-2 \brack \beta-1-j} {\nu-k-1 \brack j} q^{j(-3(\nu-1)+2 (k-1) - 2 \epsilon)} = \\
	\nonumber \frac{[\beta-j]}{[\beta-j]} \frac{[k-2]! [k-1]}{[\beta-1-j]![k+j-\beta-1]!} {\nu-k-1 \brack j} q^{j} q^{j(-3\nu+2 k - 2 \epsilon)} = \\
	\nonumber q^{j} [\beta-j] {k-1 \brack b-j}  {\nu-k-1 \brack j} q^{j(-3\nu+2 k - 2 \epsilon)} = \\	
	\nonumber q^{j} [\beta-j]  \term(\nu,k,\beta,\epsilon,j). \end{align}

For the next equation, set $j' = j+1$.
\begin{align} \label{foobar1} [\nu-k-1] \term(\nu-1,k,\beta-1,\epsilon+1,j) =  \\ 
	\nonumber[\nu-k-1] {k-1 \brack \beta-1-j} {\nu-k-2 \brack j} q^{j(-3(\nu-1)+2k - 2 - 2\epsilon)}  = \\
	\nonumber \frac{[k-1]! [\nu-k-2]! [\nu-k-1]}{[\beta-1-j]![k+j-\beta]![j]![\nu-k-j-2]!} q^{j} q^{3 \nu - 2k + 2 \epsilon} q^{(j+1)(-3\nu+2k - 2 \epsilon)} =\\
	\nonumber q^{j'} q^{3 \nu - 2k + 2 \epsilon - 1} \frac{[k-1]! [\nu-k-1]!}{[\beta-j']![k-1+j'-\beta]![j'-1]![\nu-k-1-j']!} q^{j'(-3\nu+2k)} =\\
	\nonumber q^{j'} q^{3 \nu - 2k + 2 \epsilon - 1} [j'] \binom{k-1}{\beta-j'}\binom{\nu-k-1}{j'} q^{j'(-3\nu+2k)} =\\
	\nonumber q^{3 \nu - 2k + 2 \epsilon - 1} q^{j'} [j'] \term(\nu,k,\beta,\epsilon,j'). \end{align}
When we sum such terms over all $j \in \ZZ$, we will replace $j'$ with $j$ harmlessly.

So altogether, \eqref{magicmanip1} is equivalent to
\begin{equation} \sum_j \term(\nu,j,\beta,0,j) \left([\beta] -  q^{j} [\beta-j] - q^{j - \beta} [j] \right) = 0. \end{equation}
In fact, each term in the sum is already zero, as one can easily verify that
\begin{equation} [\beta] =  q^{j} [\beta-j] + q^{j - \beta} [j]. \end{equation}
\end{proof}

%========================================================
\subsection{Generating functions} \label{ssec:magicgenerating}
%========================================================

\begin{notation} Let $C$ and $D$ be integers with the same parity. Write $[[C,D]]$ for the set of integers $r$ satisfying
	\begin{itemize}
		\item $C \le r \le D$, and
		\item $r$ has the same parity as $C$ and $D$.
	\end{itemize}
Thus when $C = D$ one has $[[C,D]] = \{C\}$, and when $D < C$ then $[[C,D]]$ is empty. We call $[[C,D]]$ a \emph{parity interval}. \end{notation}

\begin{defn} For the purpose of Theorem \ref{thm:generatingmagic}, let
\begin{equation} \label{range1} X = [[2-k,k-2]] \cup [[3k-4\nu-2 \epsilon +2,k-2\nu-2\epsilon-2]] \qquad \text{ if } 1 \le k \le \nu-1, \end{equation}
\begin{equation} \label{range2} X' = [[2-k, k-2\nu-2\epsilon-2]] \cup [[3k-4\nu-2\epsilon+2,k-2]] \qquad \text{ if } \nu+1+\epsilon \le k \le 2\nu-1+\epsilon. \end{equation}
\end{defn}

\begin{lem} The unions in both $X$ and $X'$ are disjoint. \end{lem}
	
\begin{proof} We treat the case of $X$, leaving the other case to the reader. Because $k \le \nu-1$, we have $k - 2 \nu - 2 \epsilon - 2 \le -k - 2 \epsilon - 4 \le -k$, where the last equality holds since $\epsilon \in \{-1,0,1\}$.  Thus the top of one parity interval in $X$ is always strictly less than the bottom of the other. \end{proof}
	
\begin{lem} \label{lem:range2alt} When $k$ is in the appropriate region, the set $X'$ satisfies
\begin{equation} \label{range2alt} X' = [[2-k,k-2]] \setminus [[k-2\nu-2\epsilon,3k-4\nu-2\epsilon]], \end{equation}
where the second parity interval is contained within the first. \end{lem}

\begin{proof} Left to the reader. \end{proof}

\begin{thm} \label{thm:generatingmagic} Fix positive integers $\nu$ and $k$. Assume that $\epsilon \in \{-1,0,1\}$. Let $x$ be a formal variable. If $1 \le k \le \nu-1$, then
\begin{equation} \label{generating} \sum_{\beta} \magic(\nu,k,\beta,\epsilon) x^{\beta} = \prod_{\lambda \in X} (1 + q^\lambda x). \end{equation}
If $\nu+1+\epsilon \le k \le 2\nu-1+\epsilon$, then
\begin{equation} \label{generating2} \sum_{\beta} \magic(\nu,k,\beta,\epsilon) x^{\beta} = \prod_{\lambda \in X'} (1 + q^\lambda x). \end{equation}
If neither condition holds (e.g. if $k = \nu$ and $\epsilon = 0$), this theorem says nothing.
\end{thm}

\begin{rem} In the proof below, the integers $M$ and $N$ which appear relate to the size of the parity intervals in $X$ and $X'$ as follows. In $X$, $M$ is the size (i.e. number of elements) of the first parity interval, and $N$ is the size of the second. In $X'$, one removes a parity interval of size $N$ from inside a parity interval of size $M$. \end{rem}

\begin{proof} Let us first treat the case when $1 \le k \le \nu-1$. Set $M = k-1$ and $N = \nu-k-1$; both are positive. In fact, $M$ is the size of the first parity interval in \eqref{range1}, and $N$ is the size of the second.
	
Set $u = \beta - j$. Then 
\begin{equation} \magic(\nu,k,\beta,\epsilon) = \sum_{j} {M \brack u} {N \brack j} q^{j(-3\nu + 2k - 2 \epsilon)}. \end{equation}
Taking the generating function over $\beta$, we have
\begin{align} \label{twosums} \sum_{\beta} & \magic(\nu,k,\beta,\epsilon) x^{\beta} = \\ \nonumber & \left(\sum_u {M \brack u} q^{-u(M-1)} q^{u(M-1)} x^u \right)  \left( \sum_j {N \brack j}  
q^{-j(N-1)} q^{j(-3\nu+2k-2\epsilon + N-1)} x^j \right). \end{align}

Let us recall the $q$-binomial theorem. Let $t$ be a formal variable. If $N \ge 0$ then
\begin{equation} \label{qbinomialthm}
\sum_j
{N \brack j} q^{-j(N-1)}t^j
=
\prod_{c=0}^{N-1}\left(1+q^{-2c}t\right).
\end{equation}
We can apply \eqref{qbinomialthm} to both sums in \eqref{twosums}. Letting $t = q^{M-1} x$ and noting that $M-1 = k-2$, we get
\begin{equation} \sum_u {M \brack u} q^{-u(M-1)} q^{u(M-1)} x^u = \prod_{c = 0}^{k-2} (1+q^{-2c+k-2} x). \end{equation}
Letting $t = q^{-3\nu + 2k - 2 \epsilon + N - 1} x$, and noting that 
\begin{equation} N-1 = \nu - k - 2, \qquad -3\nu + 2k - 2 \epsilon + N - 1 = k - 2 \nu - 2 \epsilon - 2, \end{equation}
we get
\begin{equation} \sum_j {N \brack j}  q^{-j(N-1)} q^{j(-3\nu + 2k - 2 \epsilon + N - 1)} x^j = \prod_{c = 0}^{\nu-k-2} (1+q^{-2c+k-2\nu-2\epsilon - 2} x). \end{equation}
Altogether, we obtain a product of $1 + q^\lambda x$ for various values of $\lambda$. It is easy to verify that 
\begin{equation} \{-2c+k-2\}_{c = 0}^{k-2} \cup \{-2c+k-2\nu-2\epsilon - 2\}_{c = 0}^{\nu-k-2} = X, \end{equation}
as desired.

Now suppose that $k \ge \nu$. Set $M = k-1$ and $N = k+1-\nu$ (the opposite of what it was before). Now one obtains $X'$ by removing a parity interval of size $N$ from one of size $M$.

Set $u = \beta - j$. Then
\begin{equation} \magic(\nu,k,\beta,\epsilon) = \sum_{j} {M \brack u} {-N \brack j}  q^{j(-3\nu + 2k - 2 \epsilon)}. \end{equation}
Therefore
\begin{align} \label{twosums2} \sum_{\beta} & \magic(\nu,k,\beta,\epsilon) x^{\beta} = \\ \nonumber & \left(\sum_u {M \brack u} q^{-u(M-1)} q^{u(M-1)} x^u \right)  \left( \sum_j {-N \brack j}  q^{j(N-1)} q^{j(-3\nu + 2k - 2 \epsilon - N + 1)} x^j \right). \end{align}

Let us remind the reader of the negative $q$-binomial theorem. If $N \ge 0$ then
\begin{equation} \label{qbinomialthmneg}
\sum_j
{-N \brack j} q^{j(N-1)}t^j
=
\prod_{c=0}^{N-1}\frac{1}{1+q^{2c}t}.
\end{equation}
We can apply \eqref{qbinomialthm} and \eqref{qbinomialthmneg} to the two sums in \eqref{twosums2} respectively. Exactly as above, letting $t = q^{M-2} x$ we get
\begin{equation} \label{firstsum} \sum_u {M \brack u} q^{-u(M-1)} q^{u(M-1)} x^u = \prod_{c = 0}^{k-2} (1+q^{-2c+k-2} x). \end{equation}
Letting $t = q^{-3\nu + 2k - 2 \epsilon - N + 1} x$, and noting that
\begin{equation} N - 1 = k - \nu, \qquad -3\nu + 2k - 2 \epsilon - N + 1 = k -2 \nu - 2 \epsilon, \end{equation}
we get
\begin{equation} \label{secondsum} \sum_j {-N \brack j}  q^{j(N-1)} q^{j(-3\nu + 2k - 2 \epsilon - N + 1)} x^j = \prod_{c = 0}^{k-\nu} \frac{1}{1+q^{2c+k-2\nu-2\epsilon} x}. \end{equation}
So the result is a product of factors $1 + q^\lambda x$ for $\lambda \in [[2-k,k-2]]$, divided by factors $1 + q^\lambda x$ for $\lambda \in [[k-2\nu-2\epsilon,3k-4\nu-2\epsilon]]$. This is precisely the product over $\lambda \in X'$, see Lemma \ref{lem:range2alt}. \end{proof}

\begin{rem} The product of $1 + q^\lambda x$, over all $d$ in a single parity interval with $M$ elements, is equal to a sum over $j$ of $x^j {M \brack j}$ times a power of $q^j$. This is the quantum binomial theorem \eqref{qbinomialthm} when $t$ is a power of $q$ times $x$. The product over a union of two parity intervals does not seem to have such a nice closed formula for the coefficient of $x^j$. In the special cases of Lemma \ref{lem:magicchuvan}, the union of the two intervals in $X$ or $X'$ is actually one interval. \end{rem}

\begin{rem} Theorem \ref{thm:generatingmagic} does not apply in all cases, but it does apply in all cases which are relevant! For example, if $a$ is even and $b$ is odd then $\gamma_3(a,b,i,k) = \magic(\nu,k,\beta,0)$, and Theorem \ref{thm:generatingmagic} does not apply when $k = \nu$. However, in this case $[k-\nu]$ divides $\gamma_2$, so we have $\gamma_2 \cdot \gamma_3 = 0$ when $k = \nu$. More generally, the factor $\gamma_2$ is zero whenever the factor $\gamma_3$ is not determined by Theorem \ref{thm:generatingmagic}. \end{rem}

Looking at \eqref{gamma3}, there are certain cases where $\gamma_3 = q^{\beta} \magic(\nu,k-1,\beta,\epsilon)$. Let us record the generating function of this variant on $\magic$. The ultimate effect is to move one factor $1+q^\lambda x$ from one parity interval to the other.

\begin{thm} \label{thm:generatingmagicfor3} Fix positive integers $\nu$ and $k$. Assume that $\epsilon \in \{-1,0,1\}$. Let $x$ be a formal variable. If $1 \le k-1 \le \nu-1$, then
\begin{equation} \sum_{\beta} q^{\beta} \magic(\nu,k-1,\beta,\epsilon) x^{\beta} = \prod_{\lambda \in X} (1 + q^\lambda x), \end{equation}
where
\begin{equation} \label{range1for3} X = [[4-k,k-2]] \cup [[3k-4\nu-2 \epsilon,k-2\nu-2\epsilon-2]] \qquad \text{ if } 1 \le k \le \nu-1. \end{equation}
If $\nu+1+\epsilon \le k-1 \le 2\nu-1+\epsilon$, then
\begin{equation} \sum_{\beta} q^\beta \magic(\nu,k-1,\beta,\epsilon) x^{\beta} = \prod_{\lambda \in X'} (1 + q^\lambda x), \end{equation}
where
\begin{equation} \label{range2for3} X' = [[4-k, k-2\nu-2\epsilon-2]] \cup [[3k-4\nu-2\epsilon,k-2]] \qquad \text{ if } \nu+1+\epsilon \le k \le 2\nu-1+\epsilon. \end{equation}
Moreover, the parity intervals in both $X$ and $X'$ are disjoint.

If neither condition holds (e.g. if $k-1 = \nu$ and $\epsilon = 0$), then this theorem does not apply.
\end{thm}

\begin{proof} Let us take \eqref{generating}, and plug in $k-1$ for $k$, and $q x$ for $x$. The result is as stated. For example, to obtain \eqref{range1for3} from \eqref{range1} one should replace $k$ with $k-1$, and then add $1$ everywhere. \end{proof}

%========================================================
\subsection{Symmetry} \label{ssec:magicsymmetry}
%========================================================

The following result will be used to prove that $\Xi$ satisfies certain symmetries (namely \eqref{Xi1symmetry} and \eqref{Xi23symmetry}).

\begin{thm} \label{thm:magicsymmetry} For any $\beta, \nu \ge 0$ and $\epsilon \in \{-1,0,1\}$, let $L = 2 \nu + \epsilon$. For any $1 \le k \le L-1$, unless $\nu \le k \le \nu+\epsilon$, we have
\begin{equation} \magic(\nu,k,\beta,\epsilon) = q^{\beta(2k-L)} \magic(\nu,L-k,\beta,\epsilon). \end{equation}
\end{thm}

This symmetry is not obvious (to our eyes) from the definition of $\magic$, and our proof uses the generating function.

\begin{proof} Note that the result for $k$ is equivalent to the result for $L-k$. Note also that both this theorem and Theorem \ref{thm:generatingmagic} have the same exceptions, e.g. when $\epsilon = 1$ then $\{k,L-k\} = \{\nu,\nu+1\}$ is forbidden. Ignoring these exceptions, either $k$ or $L-k$ is less than $\nu$. So we need only prove the result for $k \le \nu-1$. 

Since we need to prove this result for all $\beta$, we can instead prove that both sides produce the same generating function. Using Theorem \ref{thm:generatingmagic}, the generating function of the left side is
\begin{equation} \sum_\beta \magic(\nu,k,\beta,\epsilon) x^\beta = \prod_{\lambda \in X} (1 + q^\lambda x) \end{equation}
with $X$ defined by \eqref{range1}. Meanwhile, the generating function of the right side is
\begin{equation} \sum_\beta \magic(\nu,L-k,\beta,\epsilon) q^{\beta(2k-L)} x^\beta = \sum_{\beta} \magic(\nu,L-k,\beta,\epsilon) (x')^{\beta}, \end{equation}
where $x' = q^{2k-L} x$. Theorem \ref{thm:generatingmagic} would rewrite this a product over various terms $1 + q^\lambda x'$, or in other words, over $1 + q^{\lambda + 2k-L} x$. Consequently,
\begin{equation} \sum_\beta \magic(\nu,L-k,\beta,\epsilon) q^{\beta(2k-L)} x^\beta = \prod_{\lambda \in Y} (1 + q^\lambda x), \end{equation}
where $Y$ is obtained from \eqref{range2} by plugging in $L-k$ for $k$. We leave it as an exercise to verify that $X = Y$, whence the desired symmetry.
% \begin{align}\nonumber  Y & = [[2-(L-k) + 2k-L, (L-k)-2\nu-2\epsilon-2 + 2k-L]] \cup [[3(L-k)-4\nu-2\epsilon+2+ 2k-L,(L-k)-2 + 2k-L]]\\ & = [[3k-4\nu-2\epsilon+2,k-2\nu-2\epsilon-2]] \cup [[-k+2,k-2]]. \end{align}
% Since $X = Y$, we have the desired symmetry.
\end{proof}

For easier citation, let us reformulate the above as a statement about $\gamma_3$.

\begin{cor} \label{cor:gammasym1} Whenever $\ell = a+b+1 = 2 \nu + \epsilon$ we have either $\nu \le k \le \nu+\epsilon$ or
	\begin{equation} \gamma_3(a,b,1,k) = q^{\beta(2k-\ell)} \gamma_3(a,b,1,\ell-k). \end{equation}
Moreover, whenever $b$ is odd we have either $\nu \le k \le \nu+\epsilon$ or
	\begin{equation} \gamma_3(a,b,2,k) = q^{\beta(2k-\ell)} \gamma_3(a,b,3,\ell-k). \end{equation}
\end{cor}

\begin{proof} One can verify (c.f. \eqref{gamma3}) that $\gamma_3(a,b,1,k) = \magic(\nu,k,\beta,\epsilon)$ where $\ell = 2 \nu + \epsilon$. Now we simply apply the previous theorem. \end{proof}

\begin{cor} \label{cor:gammasym2} Whenever $\ell = a+b+1 = 2 \nu + \epsilon$ and $b$ is even we have either $k = \ell-1$ or $\nu \le k \le \nu+\epsilon-1$ or
	\begin{equation} \gamma_3(a,b,2,k) = q^{\beta(2k-\ell)} \gamma_3(a,b,3,\ell-k). \end{equation}
\end{cor}

\begin{proof} In this case we have 
	\begin{equation} \gamma_3(a,b,2,k) = \magic(\nu,k,\beta,\epsilon-1), \qquad \gamma_3(a,b,3,\ell-k) = q^{\beta} \magic(\nu,k-1,\beta,\epsilon-1)\end{equation}
for some $\epsilon \in \{0,1\}$ depending on $a$. Now let us apply Theorem \ref{thm:magicsymmetry} for $L = \ell-1$, and with $\epsilon - 1$ instead of $\epsilon$. Note that this rules out $k = \ell-1$ since we require $k \le L-1$, and it also rules out $\nu \le k \le \nu+\epsilon-1$. We find that
\begin{equation} \magic(\nu,k,\beta,\epsilon-1) = q^{\beta(2k-\ell+1)} \magic(\nu,\ell-1-k,\beta,\epsilon-1), \end{equation}
which is the desired equality. \end{proof}

%========================================================
\subsection{Telescoping sums and telescoping products} \label{ssec:magictelescope}
%========================================================

In our recursive formulas one must take sums of $\magic$ itself, and the generating function versions of these sums telescope in a tidy way.  Here is one example (used to prove \eqref{recursiveformula1} in the case where $a$ is odd and $b$ is even).  

\begin{thm}\label{thm:telescoping sum}
Suppose that $\ell = 2 \nu$ and $\nu \le k < \ell$. Then
\begin{align} \nonumber (-1)^{k} & (q^{-2k} - q^{-2\nu}) \magic(\nu,k,\beta,0) q^{k(k - \beta - \ell+1)} = \\ \label{telescope1} & q^{-2 \ell + 2} \sum_{c = \ell-k}^{k-1} (-1)^c \magic(\nu,c,\beta,-1) q^{c(c - \beta - \ell +3)}. 
\end{align}
\end{thm}

There are many combinatorial formulas where a sum of binomial coefficients equals another binomial coefficient. In the above theorem, a sum of double binomial coefficients equals
another sum of double binomial coefficients. We found this relation to be surprising and intriguing.

\begin{proof}
First we observe what happens for the special case $k = \nu$. The left side is zero thanks to the factor $q^{-2k} - q^{-2\nu}$. The right side is zero because the sum is empty.

Now assume $\nu+1 \le k \le \ell-1$. We prove the result using the generating functions from Theorem~\ref{thm:generatingmagic}. We write $LHS(\beta)$ for the left-hand side of \eqref{telescope1}, and $LHS(x)$ for its generating function $LHS(x) = \sum_{\beta} LHS(\beta) x^\beta$. We have
\begin{equation} \begin{split}
LHS(x) = 
(-1)^k 
(q^{-2k}-q^{-2\nu})
q^{-k(2\nu-2k-1)}
\sum_{\beta} (xq^{-k})^{\beta} \magic(\nu,k,\beta,0) \\ = 
(-1)^k 
(q^{-2k}-q^{-2\nu})
q^{-k(2\nu-2k-1)}
\prod_{\lambda \in \mathcal{L}_1}(1+q^\lambda x),
\end{split}
\end{equation}
where
\begin{equation}
\mathcal{L}_1 = [[2 - 2k, -2\nu-2]] \cup [[2k-4\nu+2, -2]]. \end{equation}
Note that we obtained $\mathcal{L}_1$ by taking $X'$ from \eqref{range2} and subtracting $k$ from all entries.

On the right-hand side of \eqref{telescope1}, $c$ takes on both small and large values, so that we should use different generating functions \eqref{generating} and \eqref{generating2}
for different parts of the sum. Instead, we simplify using the symmetry of Theorem \ref{thm:magicsymmetry}. For this purpose, given any value of $c$, let $d = 2\nu-1 -c$. Half of the
sum has $\nu \le c \le k-1$, while the other half can be viewed as a sum over $d$ with $\nu \le d \le k-1$. On this half we change our index of summation from $d$ to $c$, which also has
the effect of sending $c = L-d$ to $d = L-c$. Thus the two halves combine to form one sum.

Note that $(-1)^c = (-1)^{d-1}$. Theorem \ref{thm:magicsymmetry} with $L = 2 \nu -1$ states for any $c$ that
\begin{equation} \magic(\nu,c,\beta,-1) = q^{\beta(c-d)} \magic(\nu,d,\beta,-1). \end{equation}
So we have
\begin{align}
q^{2 \ell-2} RHS(\beta) = \sum_{c = \ell-k}^{k-1} (-1)^c \magic(\nu,c,\beta,-1) q^{c(c - \beta - \ell +3)} \\
\nonumber = \sum_{c = \nu}^{k-1} (-1)^c \magic(\nu,c,\beta,-1) q^{c(-d - \beta +2)} + \sum_{c = \ell-k}^{\nu-1} (-1)^c \magic(\nu,c,\beta,-1) q^{c(-d - \beta +2)} \\
\nonumber = \sum_{c = \nu}^{k-1} (-1)^c q^{-c\beta} \magic(\nu,c,\beta,-1) q^{-cd} q^{2c}  \qquad + \qquad \qquad \qquad \qquad \qquad \qquad \\ \nonumber \sum_{c = \ell-k}^{\nu-1} (-1)^{d-1} q^{-c\beta} q^{\beta(c-d)} \magic(\nu,d,\beta,-1) q^{-cd} q^{2c} \\
\nonumber = \sum_{c = \nu}^{k-1} (-1)^c q^{-c\beta} \magic(\nu,c,\beta,-1) q^{-cd} q^{2c} - \sum_{d = \nu}^{k-1} (-1)^d q^{-d\beta} \magic(\nu,d,\beta,-1) q^{-cd} q^{2c} \\
\nonumber = \sum_{c = \nu}^{k-1} (-1)^c q^{-c\beta} \magic(\nu,c,\beta,-1) q^{-cd} (q^{2c} - q^{2d}).
\end{align}
At the last step, one swaps $c$ and $d$ in the $d$-indexed sum, and then combines the sums. Now we can take the generating function and we obtain
\begin{equation} \begin{split}
q^{2 \ell-2} RHS(x) = \sum_{\beta} \sum_{c = \ell-k}^{k-1} (-1)^c \magic(\nu,c,\beta,-1) q^{c(c - \beta - \ell +3)} x^{\beta} \\
= \sum_{c = \nu}^{k-1} (-1)^c  q^{-cd} (q^{2c} - q^{2d}) \sum_{\beta} (xq^{-c})^{\beta} \magic(\nu,c,\beta,-1) \\
= \sum_{c = \nu}^{k-1} (-1)^c  q^{-cd} (q^{2c} - q^{2d}) \prod_{\lambda \in \mathcal{R}_1(c)}(1+q^\lambda x), \end{split}
\end{equation}
where
\begin{equation}
\mathcal{R}_1(c) = [[2 - 2c, -2\nu]] \cup [[2c - 4 \nu + 4, -2]].
\end{equation}
We can rewrite this union of parity intervals as
\begin{equation} \mathcal{R}_{1}(c) = [[2 - 2c, -L-1]] \cup  [[2-2d, -2]]. \end{equation}

Note that the interval $[[2k-4\nu+2,-2]]$ is contained in $\mathcal{L}_1$, and in $\mathcal{R}_1(c)$ for all $c$ with $\nu \le c \le k-1$. Dividing both sides by this common factor (i.e. the product of $(1+q^{\lambda} x)$ over $\lambda$ in this interval), it remains to prove that
\begin{equation}
\label{eqn:first tricky recurrence case}
\begin{split}
(-1)^k 
(q^{-2k}-q^{-2\nu})
q^{-k(2\nu-2k-1)}
\prod_{\lambda \in \mathcal{L}}(1+q^\lambda x)
\\=
\sum_{\nu \leq c \leq k-1, d = 2\nu-1-c}
q^{-4\nu+2}
(-1)^cq^{-cd}(q^{2c}-q^{2d}) 
\prod_{\lambda \in \mathcal{R}(c)}(1+q^\lambda x)
\end{split}
\end{equation}
where
\begin{align*}
\mathcal{L} &= [[-2k+2, -2\nu-2]], \\
\mathcal{R}(c) &= [[2-2c, -2\nu]] \cup [[2c - 4\nu+4, 2k-4\nu]].
\end{align*}
We can rewrite the latter as
\begin{equation} \mathcal{R}(c) = [[2-2c,-L-1]] \cup [[2-2d, 2k-2L-2]]. \end{equation}
The sizes of $\mathcal{L}$ and $\mathcal{R}(c)$ are both $k-\nu-1$, one fewer than the number of terms in the sum.

Let us make some observations about the sets $\mathcal{R}(c)$ as $c$ varies. It has two (parity) intervals, the bottom interval $[[2-2c,-1-L]]$ and the top interval $[[2-2d,2k-2L-2]]$. When $c$ takes its minimum value $\nu$, the bottom interval is empty, and the top interval begins at $3-L$ and has size $B = k-\nu-1$. As $c$ grows by $1$, the top interval loses its bottom-most term and the bottom interval gains a new bottom-most term. When $c$ takes its maximum value $k-1$, then the top interval is empty, and the bottom interval has size $B$. This matches with the fact that there are $B+1$ terms in the sum.

Here is a pleasant way to repackage the right-hand side of \eqref{eqn:first tricky recurrence case}. Let $Z = 2-L$, let $B = k-\nu-1$, and let $a = c-\nu$. Then the bottom interval has $a$ terms, and they are of the form $Z-1-2i$ for $1 \le i \le a$. The top interval has $B-a$ terms, and they are of the form $Z+1+2i$ for $a \le i \le B-1$. Now
\begin{equation} q^{2c}-q^{2d} = q^{L+1+2a} - q^{L-2a-1} = q^L(q^{2a+1}-q^{-(2a+1)}), \end{equation}
\begin{equation} q^{-cd} = q^{-(\nu+a)(\nu-1-a)} = q^{a(a+1)-\nu(\nu-1)} = q^{2 \binom{a+1}{2} - 2 \binom{\nu}{2}}. \end{equation}
Finally, $q^{-4\nu+2} = q^{-2L}$. For reasons of space, we may sometimes write $q^\lambda - q^{-\lambda}$ as $(q-q^{-1})[\lambda]$, although it is the expansion $q^{\lambda} - q^{-\lambda}$ which makes it easier to see certain cancellations in the proof.

Combining these observations, the right-hand side of \eqref{eqn:first tricky recurrence case} is equal to
\begin{equation} \label{rewriteRHSoftelescope} (-1)^{\nu} q^{-L} q^{-2 \binom{\nu}{2}} (q-q^{-1}) \sum_{a=0}^{B} (-1)^a  [2a+1] q^{2 \binom{a+1}{2}} \prod_{i=1}^{a} (1+q^{Z-1-2i} x) \prod_{i=a}^{B-1} (1+q^{Z+1+2i} x). \end{equation}
Similarly, the left-hand side of \eqref{eqn:first tricky recurrence case}, up to a sign and a power of $q$, is equal to 
\[ (q-q^{-1}) [B+1] \prod_{i=2}^{B+1} (1+q^{Z-1-2i} x). \]

Having seen the kind of equality we want, we extract this as Theorem \ref{thm:reformedtelescope} below. The equality \eqref{eqn:first tricky recurrence case} follows from Theorem \ref{thm:reformedtelescope} with a little bookkeeping. Note that by replacing the variable $x$ with $q^Z x$, we can effectively ignore the variable $Z$. We have also divided both sides of the equation by $(q-q^{-1})$.
\end{proof}

\begin{thm} \label{thm:reformedtelescope} Fix an integer $B \ge 0$. We have
\begin{align} \nonumber \sum_{a=0}^B (-1)^a &  [2a+1] q^{2 \binom{a+1}{2}} \prod_{i=1}^{a} (1+q^{-1-2i} x) \prod_{i=a}^{B-1} (1+q^{+1+2i} x) = \\ & (-1)^B q^{-2B}[B+1] \prod_{i=2}^{B+1} (q^{2i+1} + x). \end{align} \end{thm}

\begin{proof} Our proof is to describe the partial sums as a telescoping product. Let
	\begin{equation} f(a) := (-1)^a  [2a+1] q^{2 \binom{a+1}{2}} \prod_{i=1}^{a} (1+q^{-1-2i} x) \prod_{i=a}^{B-1} (1+q^{1+2i} x), \end{equation}
	\begin{equation} PS(k) = \sum_{a=0}^k f(a). \end{equation}
Multiplying each term in the first product by $q^{2i+1}$ would multiply $f(a)$ by $q^a q^{2 \binom{a+1}{2}}$, so we can rewrite $f(a)$ as
\begin{equation} \label{fa} f(a) = (-1)^a q^{-a} [2a+1] \prod_{i=1}^{a} (q^{1+2i}+x) \prod_{i=a}^{B-1} (1+q^{1+2i} x). \end{equation}

Then we prove by induction that
\begin{equation} \label{ftoPSratio} f(a)/PS(a-1) = -q^{a-2}\frac{[2a+1]}{[a]} \frac{q^3+x}{1+q^{2a-1} x}, \qquad a \ge 1, \end{equation}
\begin{equation} \label{PStoPSratio} PS(a)/PS(a-1) = -q^{-2} \frac{[a+1]}{[a]} \frac{q^{1+2(a+1)} + x}{1+q^{2a-1} x}, \qquad a \ge 1, \end{equation}
\begin{equation} \label{PSformula} PS(a) = (-1)^a q^{-2a} [a+1] \prod_{i=2}^{a+1} (q^{1+2i} + x) \prod_{i=a}^{B-1} (1+q^{1+2i} x) , \qquad a \ge 0. \end{equation}
The desired result is just \eqref{PSformula} for $a=B$. The base case is the formula for $PS(0)$, which follows directly from \eqref{fa}.

Let $a \ge 1$ and suppose that \eqref{PSformula} correctly describes $PS(a-1)$. It is straightforward to verify \eqref{ftoPSratio} from the formulas above; $f(a)$ has an extra factor of $q^{1+2i}+ x$ where $i=1$, and $PS(a-1)$ has an extra factor of $1+q^{1+2i} x$ where $i=a-1$. Continuing, we have
\begin{equation} \frac{PS(a)}{PS(a-1)} = \frac{PS(a-1)}{PS(a-1)} + \frac{f(a)}{PS(a-1)} = 1 + \frac{f(a)}{PS(a-1)}. \end{equation}
For the next calculation we multiply numerator and denominator by $(q-q^{-1})$, replacing $[a]$ with $q^a - q^{-a}$, etcetera. Adding $1$ to \eqref{ftoPSratio} and putting both terms over a common denominator, the numerator has some nice cancellation:
\begin{align} \nonumber (q^a - q^{-a})&(1+q^{2a-1} x) - q^{a-2} (q^{2a+1}-q^{-(2a+1)})(q^3+x) \\ \nonumber & = (q^a - q^{-a} - q^{3a+2} + q^{-a}) + (q^{3a-1} - q^{a-1} - q^{3a-1} + q^{-a-3})x \\ \label{ratioachieved} & = (q^a - q^{3a+2}) + (q^{-a-3} - q^{a-1})x \\ \nonumber & = -q^{-2}(q^{a+1} - q^{-(a+1)})(q^{1+2(a+1)} + x). \end{align}
From this calculation \eqref{PStoPSratio} follows. It is easy to multiply the formula for $PS(a-1)$ from \eqref{PSformula} with the ratio $PS(a)/PS(a-1)$ from \eqref{PStoPSratio} and verify that the formula \eqref{PSformula} correctly describes $PS(a)$. \end{proof}

\begin{rem}
The computation of \eqref{PStoPSratio} (achieved in \eqref{ratioachieved}) was the only place where we did any ``serious'' algebra - and it is at least slightly surprising that it works. The right-hand side of \eqref{ftoPSratio} is a rational expression in $x$ with coefficients in $\mathbb{C}[q]$.  For a typical such expression $R$, one would not expect both $R$ and $R+1$ to factor so cleanly.
\end{rem}

%========================================================
\subsection{Variations on the telescope} \label{ssec:magictelescopevariations}

There are several variants on the previous results that we need for our recursive proofs. There are a few wrinkles that appear, but no significant new ideas. We provide most of the details, even when they become repetitive, because these details are very hard to reconstruct. We focus on the differences between the proofs below and the proofs of the previous section, and ignore a lot of the bookkeeping.

The following theorem arises in the proof of~\eqref{recursiveformula1} when both $b$ and $a$ are even.

\begin{thm}
\label{telescope a and b even}
Setting $\ell=2\nu+1$, when $\nu+1 \le k < \ell$ we have:
\begin{equation}
\begin{split}
(-1)^{k} (q^{2k} - q^{2\nu}) (q^{2\nu} - q^{2k-2}) \magic(\nu,k,\beta,+1) q^{k(k-\beta-\ell-2)}
 \\= \sum_{c = \ell-k}^{k-1} (-1)^c (q^{-2\nu}-q^{-2c}) \magic(\nu,c,\beta,0) q^{c(c-\beta-\ell+4)}. 
 \end{split}
 \end{equation}\end{thm}

\begin{proof}
The argument is very similar to the one in the previous theorem. Set $L = 2\nu$ and $d = L-c$. The line of symmetry is at $c = d = \nu$, but the corresponding term of the sum vanishes due to the factor $q^{-2\nu} - q^{-2c}$. The equivalent generating function identity is
\begin{equation}
\begin{split}
(-1)^k
(q^{k - \nu} - q^{\nu - k}) 
(q^{\nu-k+1} - q^{k-\nu-1}) 
q^{k(k - (2\nu+1))}
\prod_{\lambda \in \mathcal{L}} (1+q^{\lambda} x)
\\=
q^{1-3\nu}
\sum_{\nu+1 \leq c \leq k-1}
(-1)^c 
q^{-cd}
(q^{c-\nu} - q^{\nu-c})
(q^{2c}-q^{2d})
\prod_{\lambda \in \mathcal{R}(c)}
(1+q^{\lambda}x)
\end{split}
\end{equation}
where $\mathcal{L} = [[-2k+2, -2\nu-4]]$, $\mathcal{R}(c) = [[2-2c, -2\nu-2]] \cup [[2c-4\nu+2, 2k-4\nu-2]]$. As before, the top interval of $\mathcal{R}(c)$ has lower bound $2-2d$.

%We can rewrite $\mathcal{R}(c)$ as $[[2-2c,-2-L]] \cup [[2-2d,2k-2L-2]]$.

Letting $a = c-(\nu+1)$ and $B = k-\nu-2$ and $Z = 2-L$, we have
\begin{equation} \prod_{\lambda \in \mathcal{R}(c)} (1+q^{\lambda}x) = \prod_{i = 1}^{a} (1+q^{Z-2-2i} x) \prod_{i=a}^{B-1} (1+q^{Z+2+2i} x). \end{equation}
This is different from the previous case, as the bottom and top intervals of $\mathcal{R}(c)$ are further apart. We also have
\begin{equation} q^{2c} - q^{2d} = q^L(q^{2a+2} - q^{-(2a+2)}), \end{equation}
\begin{equation} q^{-cd} = q^{-(\nu +1 + a)(\nu-1 - a)} = q^{-\nu^2 + (a+1)^2} = q^{-\nu^2}q^{2\binom{a+1}{2} + a + 1}, \end{equation}
and an extra factor
\begin{equation} q^{c-\nu} - q^{\nu-c} = q^{a+1} - q^{-(a+1)}. \end{equation}

Now the result follows from the following theorem and some bookkeeping. This time we have divided both sides of the equation by $(q-q^{-1})^2$.
\end{proof}

\begin{thm} \label{reformed telescope a and b even} Fix an integer $B \ge 0$. We have
\begin{align} \nonumber \sum_{a=0}^B (-1)^a & [a+1][2a+2]q^{2 \binom{a+1}{2}+a} \prod_{i=1}^{a} (1+q^{-2-2i} x) \prod_{i=a}^{B-1} (1+q^{2+2i} x) = \\ & (-1)^B q^{-2B}[B+1][B+2] \prod_{i=2}^{B+1} (q^{2i+2} + x). \end{align} \end{thm}

\begin{proof} Let $f(a)$ be the $a$-th term in the sum, and $PS(a)$ be the $a$-th partial sum. Then we have
\begin{equation} f(a) = (-1)^a q^{-a} [a+1][2a+2]  \prod_{i=1}^{a} (q^{2i+2} + x) \prod_{i=a}^{B-1} (1+q^{2+2i} x), \end{equation}
\begin{equation} f(a)/PS(a-1) = (-1) q^{a-2} \frac{[2a+2](q^4+x)}{[a](1+q^{2i} x)}, \end{equation}
\begin{equation} PS(a)/PS(a-1) = (-1) q^{-2} \frac{[a+2](q^{2i+4}+x)}{[a](1+q^{2i} x)}, \end{equation}
\begin{equation} PS(a) = (-1)^a q^{-2a}[a+1][a+2] \prod_{i=2}^{a+1} (q^{2i+2} + x) \prod_{i=a}^{B-1} (1+q^{2+2i} x). \end{equation}
We leave the inductive proof to the reader; the outline is precisely as in Theorem \ref{thm:reformedtelescope}.
\end{proof}

%
% We prove this identity by explicitly telescoping the sum on the right hand side - the idea being that for $N$ between $\nu+1$ and $k+1$, the $N$th partial sums
% \begin{equation*}
% \text{RHS}(N)=
% \sum_{\nu+1 \leq c \leq N}
% (-1)^c
% q^{-cd}
% (q^{c-\nu} - q^{\nu-c})
% (q^{2c}-q^{2d})
% \prod_{\lambda \in \mathcal{R}(c)}
% (1+q^{\lambda}x)
% \end{equation*}
% themselves have a product form, namely:
% \begin{equation*}
% RHS(M)=RHS(\nu+1) \prod_{N=\nu+1}^M \frac{RHS(N+1)}{RHS(N)}
% \end{equation*}
% where the ratios in the product are given by the identity
% \begin{align*}
% \frac{RHS(N+1)}{RHS(N)}&=
% -\frac{ {\left(q^{4N + 4} - q^{4\nu}\right)} {\left(q^{2\nu + 2} + x\right)}}{{\left(q^{2N + 2} x + q^{4\nu}\right)} {\left(q^{2N} - q^{2\nu}\right)} q^{2}} + 1
% \\&=
%  -\frac{{\left(q^{2N + 2} + x\right)}{\left(q^{2N + 4} - q^{2\nu}\right)} q^{2\nu - 2}}{{\left(q^{2N + 2} x + q^{4\nu}\right)} {\left(q^{2N} - q^{2\nu}\right)}}. \\
% \end{align*}
% As before we show this inductively, the base case $N=\nu+1$ is vacuous, and the inductive step boils down to checking the equality of the two rational functions above. Multiplying out the factors when $N=k-1$ then yields the desired left-hand side.
% \end{proof}
%\begin{thm}

Here is the next variant we need to prove, corresponding to~\eqref{recursiveformula1} when both $b$ and $a$ are odd.

\begin{thm}
\label{telescope a and b odd}
Setting $\ell=2\nu-1$, when $\nu \le k < \ell$ we have:
\begin{equation}
\begin{split}
(-1)^{k} q^{\ell -1} (1 - q^{2 \beta}) \magic(\nu,k,\beta,-1) q^{k(k-\beta-\ell)}
 \\= \sum_{c = \ell-k}^{k-1} (-1)^c (1-q^{2c + 6 - 4 \nu}) \magic(\nu-1,c,\beta-1,-1) q^{c(c-\beta-\ell+3)}.
 \end{split}
 \end{equation}\end{thm}
 
\begin{rem} It seems like one might be able to prove Theorem \ref{telescope a and b odd} using Lemma \ref{lem:magicmanip1} and Theorem \ref{thm:telescoping sum}, but we were unable to make this approach work. \end{rem}
 
\begin{proof} We highlight only the differences with previous proofs. The factor of $(1-q^{2 \beta})$ on the left-hand side gives a difference between two generating functions, so that
\begin{equation} LHS(x) = (-1)^{k} q^{2\nu-2} q^{k(k-2\nu+1)} \left( \prod_{\lambda \in \mathcal{L}_1}(1+q^\lambda x) - \prod_{\lambda \in \mathcal{L}_1}(1+q^{\lambda + 2} x) \right), \end{equation}
where $\mathcal{L}_1 = [[2-2k,-2\nu]] \cup [[2k-4\nu+4,-2]]$. The two products have a common factor $\prod (1+q^{\lambda x})$ ranging over $\lambda \in [[4-2k, -2\nu]] \cup [[2k-4\nu+6,-2]]$. Dividing by this common factor, what remains in the parentheses is
\begin{equation}
\begin{split} (1+q^{2-2k} x)(1+q^{2k-4\nu+4} x) - (1+q^{2-2\nu} x)(1 + q^0 x)\\
	\label{thiscrazyguy} = (q^{2-2k} + q^{2k-4\nu+4} - q^{2-2\nu} - 1)x + (q^{6-4\nu} - q^{2-2\nu}) x^2 . \end{split} \end{equation}

Meanwhile, on the right-hand side, it is $\nu-1$ which plays the old role of $\nu$, which shifts our line of symmetry. Let $L = 2\nu-3$ and $d = L - c$. Our sum ranges over $\nu-1 \le c \le k-1$ and over $\nu-1 \le d \le k-2$; this time the two halves do not exactly match, so we keep track of the $c = k-1$ term separately! Also, using $\beta-1$ rather than $\beta$ will shift the generating function, multiplying the result by $x$. Using similar manipulations as before, we have
\begin{equation} RHS(x) = x(RHS_1(x) + RHS_2(x)) \end{equation}
where
\begin{equation} RHS_1(x) = \sum_{c = \nu-1}^{k-2} (-1)^c (q^{-2c} - q^{-2d}) q^{-cd} \prod_{\lambda \in \mathcal{R}_1(c)} (1+q^{\lambda} x), \end{equation}
\begin{equation} RHS_2(x) = \sum_{c = k-1}^{k-1} (-1)^c (1-q^{-2d}) q^{-cd} \prod_{\lambda \in \mathcal{R}_1(c)} (1+q^{\lambda} x). \end{equation}
Here $\mathcal{R}_1(c) = [[2-2c,-2(\nu-1)]] \cup [[2c-4(\nu-1)+4,-2]]$. The top interval is still $[[2-2d,-2]]$.

We remove the common interval $[[2k-4\nu+6,-2]]$ from $\mathcal{L}_1$ to obtain $\mathcal{L}$, and from $\mathcal{R}_1(c)$ to obtain $\mathcal{R}(c)$.

Let $Z = 2-L$ and $B = k-\nu$ and $a = c - \nu-1$. Then
\begin{equation} q^{-2c}-q^{-2d} = q^{-(L+1+2a)} - q^{-(L-2a-1)} = -q^{-L}(q^{2a+1}-q^{-(2a+1)}), \end{equation}
\begin{equation} q^{-cd} = q^{-(\nu-1+a)(\nu-2-a)} = q^{2 \binom{a+1}{2} - 2 \binom{\nu-1}{2}}. \end{equation}
We have:
\begin{multline*} \label{rewriteRHSoftelescopeODD} RHS_1(x) = (-1)^{\nu} q^{-2 \binom{\nu-1}{2}-L} (q-q^{-1}) \\ \times \sum_{a=0}^{B-1} (-1)^a  [2a+1] q^{2 \binom{a+1}{2}} \prod_{i=1}^{a} (1+q^{Z-1-2i} x) \prod_{i=a}^{B-1} (1+q^{Z+1+2i} x). \end{multline*}
Up to the leading scalar (and replacing $x$ with $q^Z x$), this is exactly the partial sum $PS(B-1)$ studied in Theorem \ref{thm:reformedtelescope}. Thus
\begin{multline*} RHS_1(x) = (-1)^{\nu} q^{-2 \binom{\nu-1}{2}-L} (q-q^{-1}) (-1)^{B-1} q^{-2(B-1)} \\ \times [B] \prod_{i=2}^{B} (q^{1+2i} + q^Z x) \prod_{i=B-1}^{B-1} (1+q^{Z+1+2i} x). \end{multline*} Meanwhile, $RHS_2(x)$ is just like the $a=B$ term except with $q^{-2c}-q^{-2d}$ replaced by $1-q^{-2d}$.
Since $d$ in this case is $L-k+1$, we have
\begin{multline*} RHS_2(x) = (-1)^{\nu-1} q^{-2 \binom{\nu-1}{2}} \\ \times \sum_{a=B}^{B} (-1)^a  (1-q^{-2L+2k-2}) q^{-a} \prod_{i=1}^{a} (q^{1+2i}+q^{Z} x) \prod_{i=a}^{B-1} (1+q^{Z+1+2i} x). \end{multline*}

Now let us add $RHS_1+RHS_2$. There is a common factor of 
\[ (-1)^{\nu+B-1} q^{-2 \binom{\nu-1}{2}} \prod_{i=2}^{B} (q^{1+2i} + q^Z x). \]
Ignoring this factor we compute that
\begin{align} \nonumber (q^B-q^{-B}) & q^{-2(B-1) - L} (1+q^{Z+1+2(B-1)} x) + (1-q^{-2L+2k-2}) q^{-B} (q^{1+2}+q^{Z} x) = \\ & (q^{5-\nu-k} - q^{5+\nu-3k} + q^{3+\nu-k} - q^{7-3\nu+k}) + (q^{5-\nu-k} - q^{9-3\nu-k}) x. \end{align}
There was a cancellation of $\pm q^{9-5\nu+k} x$. Multiplied by $x$, this agrees with the term in \eqref{thiscrazyguy} up to an overall factor of $-q^{3+\nu-k}$. After some bookkeeping, one obtains the desired result.
\end{proof}

% \BE{I did the above bookkeeping!}

Here is the final variant, corresponding to $b$ odd and $a$ even.

\begin{thm} \label{telescope b odd a even}
Setting $\ell = 2 \nu$, when $\nu \le k < \ell$ we have:
\begin{align}
(-1)^{k} q^{2 \ell - 3} (1-q^{2 \beta}) (q^{k-\nu} - q^{\nu-k}) \magic(\nu,k,\beta,0) q^{k(k-\beta-\ell)} =
 \\ \nonumber \sum_{c = \ell-k}^{k-1} (-1)^c (q^{c+1-\nu}-q^{\nu-c-1}) (q^{\ell-2-c} - q^{c+2-\ell}) \magic(\nu-1,c,\beta-1,0) q^{c(c-\beta-\ell+4)}.
\end{align}\end{thm}

The proof of this theorem merely combines the techniques above, and has the advantage of reusing Theorem \ref{reformed telescope a and b even}.

\begin{proof} First we evaluate the left-hand side. It has generating function
\begin{equation} LHS(x) = (-1)^k q^{4\nu-3} (q-q^{-1})[k-\nu] q^{k(k-2\nu)} \left(\prod (1+q^{\lambda} x) - \prod(1+q^{\lambda + 2} x)\right) \end{equation}
where both products range over the set $[[2-2k,-2\nu-2]] \cup [[2k-4\nu+2,-2]]$. Let
\begin{equation} H = (-1)^k q^{k(k-2\nu)} (q-q^{-1})[k-\nu] \prod (1+q^{\lambda} x) \end{equation}
where the product ranges over $Y = [[4-2k,-2\nu-2]] \cup [2k-4\nu+4,-2]]$. Then
\begin{equation} LHS(x) = H q^{4\nu-3} \left( (1+q^{2-2k} x)(1+q^{2k-4\nu+2} x) - (1+q^{-2\nu} x)(1+q^{0} x) \right). \end{equation}
Multiplying this out, we get
\begin{equation} \label{LHSfinal} LHS(x) = H x \left( (q^{4\nu-2k-1} + q^{2k-1} - q^{2\nu-3} - q^{4\nu-3}) + (q^{1} - q^{2\nu-3})x \right). \end{equation}

Now consider the right-hand side. Since $\nu-1$ is the parameter for $\magic$ and $\epsilon =0$, our symmetry involves $L = 2\nu-2$. Let $d = L-c$. Plugging in $c = \nu-1$ yields zero. The sum ranges over $\nu \le c \le k-1$ and $\nu \le d \le k-2$, so again there is a mismatch. Each summand can be rewritten as
\begin{equation} (-1)^c (q-q^{-1})^2 [c+1-\nu][d] q^c q^{-cd} q^{-c(\beta-1)} \magic(\nu-1,c,\beta-1,0). \end{equation}
Note that $(-1)^c = (-1)^d$. Thus we can write
\begin{equation} RHS(x) = x(RHS_1(x) + RHS_2(x)) \end{equation}
where
\begin{equation} RHS_1(x) = \sum_{c = \nu}^{k-2} (-1)^c (q-q^{-1})^2 (q^c[c+1-\nu][d] + q^d[d+1-\nu][c]) q^{-cd} \prod_{\lambda \in \mathcal{R}_1(c)} (1+q^{\lambda} x), \end{equation}
\begin{equation} RHS_2(x) = \sum_{c = k-1}^{k-1} (-1)^c (q-q^{-1})^2 q^c[c+1-\nu][d] q^{-cd} \prod_{\lambda \in \mathcal{R}_1(c)} (1+q^{\lambda} x). \end{equation}
Here $\mathcal{R}_1(c) = [[2-2c,-2(\nu-1)-2]] \cup [[2c-4(\nu-1)+2,-2]]$. The top interval is still $[[2-2d,-2]]$. Both $RHS_1(x)$ and $RHS_2(x)$ share a common factor of $\prod(1+q^{\lambda} x)$ ranging over $[[2(k-1)-4(\nu-1)+2,-2]] = [[2k-4\nu+4,-2]]$.

Let $Z = -2\nu+4 = 2-L$ and $B = k-\nu-2$ and let $a = c - \nu$. Then $d = \nu-2-a$ and $c = \nu+a$, so $c-d = 2(a+1)$. Expanding $[c]$ and $[d]$ we have
\begin{align} (q-q^{-1})&(q^c[c+1-\nu][d] + q^d[d+1-\nu][c]) \\ \nonumber  & =
	  (q^{c+d} - q^{c-d})[a+1] + (q^{c+d} - q^{d-c})[-1-a] \\\nonumber   &  = -[a+1](q^{c-d} - q^{d-c}) \\ \nonumber  & = -(q-q^{-1})[a+1][2(a+1)]. \end{align}
Meanwhile
\begin{equation} q^{-cd} = q^{(a+\nu)(a+2-\nu)} = q^{2\binom{a+1}{2} + a - \nu^2 + 2\nu}. \end{equation}
Thus we can rewrite $RHS_1(x)$ as
\begin{align} \nonumber RHS_1(x) = (-1)^{\nu+1}& q^{2\nu-\nu^2} (q-q^{-1})^2 \sum_{a=0}^{B} (-1)^a [a+1][2a+2] q^{2\binom{a+1}{2} + a}\\ & \prod_{i=1}^{a} (1+q^{Z-2-2i} x) \prod_{i=a}^{B-1} (1+q^{Z+2+2i}x) \prod_{[[2k-4\nu+2,-2]]} (1+q^{\lambda} x). \end{align}

Now we can apply Theorem \ref{reformed telescope a and b even} to deduce that
\begin{align} \nonumber RHS_1(x) = (-1)^{\nu+1} q^{2\nu-\nu^2} & (q-q^{-1})^2 (-1)^B q^{-2B} [B+1][B+2] \\ & \prod_{i=2}^{B+1} (q^{2+2i} + q^Z x) \prod_{[[2k-4\nu+2,-2]]} (1+q^{\lambda} x). \end{align}
For comparison we prefer $(1+q^{Z-2-2i}x)$ to $(q^{2+2i} + q^Z x)$, so we renormalize:
\begin{align}\nonumber RHS_1(x) = (-1)^{\nu+1+B} q^{2\nu-\nu^2} & (q-q^{-1})^2 q^{2 \binom{B+2}{2} - 2} [k-\nu-1][k-\nu] \\ & \prod_{i=2}^{B+1} (1+ q^{Z-2-2i} x) \prod_{[[2k-4\nu+2,-2]]} (1+q^{\lambda} x). \end{align}
The penultimate product ranges over the interval $[[4-2k,-2\nu-2]]$. Evaluating $B$ we conclude that
\begin{equation} RHS_1(x) = -H (q-q^{-1}) q^{-k(k-2\nu)}q^{2\nu-\nu^2+2 \binom{k-\nu}{2} - 2} [k-\nu-1] (1+q^{2k-4\nu+2} x). \end{equation}
Simplifying the exponent of $q$ we have
%\begin{equation} -k(k-2\nu) + 2\nu-\nu^2 + 2 \binom{k-\nu}{2} - 2 = -k + 3 \nu - 2 \end{equation}
\begin{equation} \label{RHS1final} RHS_1(x) = -H (q-q^{-1}) q^{-k + 3 \nu - 2} [k-\nu-1] (1+q^{2k-4\nu+2} x). \end{equation}
	
Meanwhile,
\begin{equation} RHS_2(x) = (-1)^{k-1} (q-q^{-1})^2 q^{k-1}[k-\nu][2\nu-k-1] q^{-(k-1)(2\nu-1-k)} \prod_{\lambda \in \mathcal{R}_1(k-1)} (1+q^{\lambda} x), \end{equation}
where the product ranges over $[[4-2k,-2\nu]] \cup [[2k-4\nu+4,-2]]$. We conclude that
\begin{equation} RHS_2(x) = - H (q-q^{-1}) q^{-k(k-2\nu)} q^{(k-1)(k+2-2\nu)} [2\nu-k-1](1+q^{-2\nu} x). \end{equation}
Simplifying the exponent of $q$ we have
% Continuing to simplify we have
% \begin{equation} -k(k-2\nu) + (k-1)(k+2-2\nu) = k - 2 + 2 \nu, \end{equation}
% so
\begin{equation} \label{RHS2final} RHS_2(x) = - H (q-q^{-1}) q^{k-2+2\nu} [2\nu-k-1] (1+q^{-2\nu} x). \end{equation}

Now it is simple arithmetic, using \eqref{LHSfinal}, \eqref{RHS1final}, and \eqref{RHS2final}, to verify that $LHS(x) = x (RHS_1(x) + RHS_2(x))$. \end{proof}

% Finally, we compute $RHS_1(x) + RHS_2(x)$ divided by $H$, and we get
% \begin{align} \nonumber
% 	-q^{-k + 3 \nu - 2}(q^{k-\nu-1} - q^{\nu+1-k})  (1+q^{2k-4\nu+2} x)
% 	- q^{k-2+2\nu}  (q^{2\nu-k-1} - q^{k+1-2\nu}) (1+q^{-2\nu x})\\
% &	= \left(-q^{2\nu-3} + q^{-2k+4\nu-1} - q^{4\nu-3} + q^{2k-1}\right) + \left(-q^{2k-2\nu-1} + q^{1} - q^{2\nu-3} + q^{2k-2\nu-1} \right)x.  \end{align}
% This exactly agrees with $LHS(x)/Hx$ as desired, see \eqref{LHSfinal}.

%% file: Preliminaries.tex
%!TEX root = QFrob2.tex

%========================================================
\section{Preliminaries} \label{sec:prelims}
%========================================================

For the rest of the paper we fix $n=3$, but when possible we make our definitions for general $n$, to give the reader some context.

\begin{notation} Let $\Om = \ZZ/n\ZZ$ be the vertices in the affine Dynkin diagram in type $\widetilde{A}_{n-1}$. Let $W_{\aff}$ be the affine Weyl group, with simple
reflections $S = \{s_i\}_{i \in \Om}$. Let $W_m$ be the quotient of $W_{\aff}$ by the normal subgroup $m \cdot \La_{\rt}$. Concretely, $m \cdot \La_{\rt}$ is the smallest normal subgroup containing $(s_0 w_0)^m$, where $w_0$ is the longest element of the finite Weyl group generated by $\{s_i\}_{1 \le i \le n-1}$. \end{notation}

\begin{defn} Let $z$ be a formal variable, and let $V_{z}$ be the free $\ZZ[z,z^{-1}]$-module with basis $\{x_i\}_{i \in \Om}$. It has an action of $W_{\aff}$ defined as
follows: \begin{equation} s_i(x_i) = z x_{i+1}, \quad s_i(x_{i+1}) = z^{-1} x_i, \quad s_i(x_j) = x_j \text{ if } j \notin \{i,i+1\}. \end{equation} For any
$m \ge 2$ let $V_m$ be the $\CC$-vector space obtained by specializing $z$ to a primitive $(nm)$-th root of unity $\ze \in \CC$. \end{defn}

This representation $V_z$ is a deformation of the reflection representation of $W_{\aff}$. The specialization $V_m$ is a faithful representation of the quotient group $W_m$.
% \footnote{Rather, it deforms the action of $W_{\aff}$ on the torus of affine $\gl_n$ rather than
% that of affine $\sl_n$. This $\gl_n$ variant is new in this paper, and is significantly easier to work with than the $\sl_n$ variant.}

\begin{notation} Let $R_z$ be the polynomial ring $\Sym(V_z)$ over the base ring $\ZZ[z,z^{-1}]$, and $R_m$ be the polynomial ring $\Sym(V_m)$ over $\CC$. Both are graded so that $V_z$ (resp. $V_m$) appears in degree $1$. \end{notation}
	
\begin{defn} For each $i \in \Om$ define
certain linear maps $R_z \to R_z$ of degree $-1$, called \emph{divided difference operators} or \emph{Demazure operators}, as follows: \begin{equation} \pa_i(f) = \frac{f - s_i f}{x_i - z x_{i+1}}. \end{equation}
These maps descend to $R_m$. \end{defn}

Demazure operators satisfy a twisted Leibniz rule:
\begin{equation} \pa_i(fg) = \pa_i(f) g + s_i(f) \pa_i(g). \end{equation}

\begin{defn}
The subalgebra of $\End_{R_z^{W_{\aff}}}(R_z)$ generated by Demazure operators $\pa_i$ for $i \in \Om$ is called the \emph{deformed affine nilCoxeter algebra} $\NC(z,n)$.
The subalgebra of $\End_{R_m^{W_m}}(R_m)$ generated by Demazure operators $\pa_i$ for $i \in \Om$ is called the \emph{exotic nilCoxeter algebra} $\NC(m,m,n)$. \end{defn}
	
The following quadratic and braid relations hold within $\NC(z,n)$ and hence within $\NC(m,m,n)$, see \cite[\S 4.1]{EJY1}:
\begin{equation} \pa_i \circ \pa_i = 0, \end{equation}
\begin{equation} \label{zeR3} z \pa_i \pa_{i+1} \pa_i = \pa_{i+1} \pa_i \pa_{i+1}. \end{equation}
\begin{equation} \pa_i \pa_j = \pa_j \pa_i \quad \text{ if } j \ne i, i \pm 1. \end{equation}
The third relation does not occur for $n=3$. Note the appearance of $z$ in \eqref{zeR3}, which should be specialized to $\ze$ when discussing operators on $R_m$.

To any word $\un{w} = (i_1, i_2, \ldots, i_{\ell})$ of length $\ell$ in the alphabet $\Om$ (such a word is usually called an \emph{expression}), we can associate the corresponding Demazure operator
\begin{equation} \pa_{\un{w}} := \pa_{i_1} \circ \pa_{i_2} \circ \cdots \circ \pa_{i_{\ell}}. \end{equation}
Using the braid relations, any two reduced expressions for the same element of $W_{\aff}$ give rise to Demazure operators which agree up to power of $z$ (which measures how many times \eqref{zeR3} must be applied). Using the quadratic relations, any non-reduced expression gives rise to the zero operator. When studying Demazure operators, we may as well choose one reduced expression for each element of $W_{\aff}$, which we do below.
	
A word is \emph{clockwise cyclic} if it has the form $(i+1,i+2,i+3,\ldots, i+\ell)$ for some $i \in \Om$ and $\ell \ge 1$ (with all indices considered modulo $n$). It is
\emph{widdershins cyclic} if it has the form $(i-1, i-2, \ldots, i-\ell)$. For example, when $n=3$, the word $(3,1,2,3,1,2,3)$ is clockwise cyclic, and $(2,1,3,2,1,3)$ is widdershins
cyclic.

\begin{rem} \label{rmk:roundabout} There are a number of interesting features of Demazure operators associated to cyclic words, see e.g. \cite[Theorem 4.25]{EJY1}. These lead to new relations in $\NC(m,m,n)$ called the \emph{roundabout relations}, see \cite[Theorem 4.35]{EJY1}. However, we were unable to use these results to simplify the proofs in this paper, so we do not recall the details. \end{rem}

\begin{defn} Fix $a, b \ge 0$, and $i \in \Om$. The word $\un{w}(a,b,i)$ has the form
	\begin{equation} \un{w}(a,b,i) = (k, k+1, \ldots, j, j+1, j, \ldots, i+1, i). \end{equation}
Here, the subword $(k,k+1, \ldots, j)$ is clockwise cyclic of length $a$, and the subword $(j, \ldots, i+1, i)$ is widdershins cyclic of length $b$. This word has length $\ell = a+b+1$, and ends in $i$. \end{defn}

\begin{ex} We have $\un{w}(3, 5, 2) = (1,2,3,1,3,2,1,3,2)$. One can view this as a clockwise word of length $4 = a+1$ and a counterclockwise word of length $6 = b+1$ which overlap in the middle index $1 = j+1$. The last index is $i = 2$. \end{ex}

\begin{ex} Here are some edge cases. We have $\un{w}(0,0,i) = (i)$. When $a = 0$, $\un{w}(0,b,i)$ is the widdershins cyclic word of length $b+1$ ending in $i$. When $b=0$, $\un{w}(a,0,i)$ is the clockwise cyclic word of length $a+1$ ending in $i$. \end{ex}
	
\begin{lem} \label{lem:abiparametrize} Every non-identity element of $W_{\aff}$ has a reduced expression $\un{w}(a,b,i)$ for a unique triple $(a,b,i)$. We write this element as $w(a,b,i)$. \end{lem}

The proof can be found in \cite[Lemma 5.4]{EJY1}. We restrict our attention to words of the form $\un{w}(a,b,i)$ henceforth.

\begin{notation} Let $\pa_{(a,b,i)} := \pa_{\un{w}(a,b,i)}$. \end{notation}

Now we discuss which polynomials need to be examined. The following easy result is \cite[Lemma 5.7]{EJY1}.

\begin{lem} \label{lem:killx123} Let $\un{w}$ be a word of length $\ell$ and $f \in R_z$ be a monomial of degree $\ell$. If $x_1 x_2 x_3$ divides $f$ then $\pa_{\un{w}}(f) = 0$.
\end{lem}

Because of this lemma, we need only examine monomials $x_1^k x_2^l x_3^m$ where at least one exponent in $\{k,l,m\}$ is zero. However, we can use symmetry to assume that $m = 0$.

\begin{defn} \label{defn:symmetries} Let $\si$ denote the rotation operator on $\Om$, for which $\si(i) = i+1$. Then $\si$ acts on $W_{\aff}$, and $\ZZ[z,z^{-1}]$-linearly on $V_z$ and $R$, by permuting indices. Thus $\si(s_i) = s_{i+1}$ and $\si(x_i) = x_{i+1}$. We set $\si(z) = z$.

Let $\tau$ be the automorphism of $W_{\aff}$ defined by $\tau(s_i) = s_{-i}$. Let $\tau$ be the $\ZZ$-linear automorphism of $V_z$ and $R$ given by $\tau(x_i) = x_{1-i}$, and $\tau(z) = z^{-1}$. \end{defn}

\begin{lem}(\cite[Lemma 4.8]{EJY1}) One has
	\begin{equation} \label{siondem} \si(\pa_i(f)) = \pa_{i+1}(\si(f)), \end{equation}
	\begin{equation} \label{tauondem} \tau(\pa_i(f)) = (-z) \pa_{-i}(\tau(f)). \end{equation}
\end{lem}

\begin{cor} We have
	\begin{equation} \label{sionabi} \si(\pa_{(a,b,i)}(f)) = \pa_{(a,b,i+1)}(\si(f)), \end{equation}
	\begin{equation} \label{tauonabi} \tau(\pa_{(a,0,i)}(f)) = (-z)^{a+1} \pa_{(0,a,-i)}(\tau(f)), \qquad \tau(\pa_{(0,a,i)}(f)) = (-z)^{a+1} \pa_{(a,0,-i)}(\tau(f)). \end{equation}
\end{cor}

\begin{proof} These equations follow from iterative applications of \eqref{siondem} and \eqref{tauondem}, after noting that $\si(\un{w}(a,b,i)) = \un{w}(a,b,i+1)$ and $\tau(\un{w}(a,0,i)) = \un{w}(0,a,-i)$. \end{proof}

\begin{rem} Meanwhile, $\tau(\un{w}(a,b,i))$ is a concatenation of a widdershins word with a clockwise word, and is not a reduced expression of the form $\un{w}(a',b',i')$ when both $a$ and $b$ are nonzero. \end{rem}

Suppose that $f = x_1^k x_2^l x_3^m$ and $k+l+m = a+b+1$, so that $\pa_{(a,b,i)}(f)$ is a scalar for degree reasons. Then $\pa_{(a,b,i)}(f)$ is fixed by $\si$, so is equal to $\pa_{(a,b,i+1)}(\si(f))$. If at least one of $\{k,l,m\}$ is zero, then up to the application of $\si$, we can assume $m=0$.

The scalars (elements of $\ZZ[z,z^{-1}]$) we aim to compute are below.

\begin{notation} Let $\Xi(a,b,i,k) = \pa_{\un{w}(a,b,i)}(x_1^k x_2^{a+b+1-k})$. \end{notation}

Our main result in \S\ref{sec:closedformula} is a closed formula for $\Xi(a,b,i,k)$. This formula depends in annoying ways on: \begin{itemize}
\item whether $a = 0$ or $a>0$ is even or $a>0$ is odd,
\item whether $b = 0$ or $b>0$ is even or $b>0$ is odd,
\item whether $k=0$ or $0 < k < a+b+1$ or $k = a+b+1$,
\item whether $i=1$ or $i=2$ or $i=3$.
\end{itemize}
One should expect on the order of $81$ separate formulas for the different cases. Setting $a=0$ or $b=0$ or $k=0$ or $k=a+b+1$ will change the formula dramatically, while switching parity of nonzero values of $a$ and $b$, or switching the index $i$, will often just result in a small adjustment. More precisely, from $a$ and $b$ we extract variables $\alpha$ and $\beta$ (depending on parity), and we describe our formula in terms of $\alpha$ and $\beta$, so that the difference between cases appears minimized.

Our formulas will involve quantum numbers. Up to an invertible scalar, these are quantum numbers in $q$, not in $z$! Remember that $z^3 = q^{-2}$. Thus the scalar $1 + z^3$ agrees with
the unbalanced quantum number $(2)_{q^{-1}} = 1 + q^{-2}$, and $1 + z^3 + z^6$ agrees with $(3)_{q^{-1}} = 1 + q^{-2} + q^{-4}$, etcetera.
%Powers of $z^3$ appear so frequently that we use the notation $y = z^3$.
For reasons of familiarity we wished to express our formulas using balanced quantum numbers like $[2]_q = q + q^{-1}$ and $[3]_q = q^{-2} + 1 + q^2$, but these can only be expressed using half-powers
of $z$. To this end, we introduce a square root of $z$, which we call $p$. So set \begin{equation} z = p^2, \qquad q = p^{-3}, \qquad p^6 = z^3 = q^{-2}. \end{equation} Henceforth
all quantum numbers and quantum binomial coefficients will be balanced with respect to the variable $q$. We give our formulas as expressions involving powers of $p$, $q$, and $z$,
but the result ultimately lives inside the ring $\ZZ[z,z^{-1}]$.

Balanced quantum binomial coefficients with respect to the variable $q$ will be denoted ${k \brack c}$. In various exponents we also use the ordinary binomial coefficients
$\binom{d}{2}$ for various integers $d$. Ultimately, the exponent of $p$ will be a quadratic polynomial in our various inputs ($a$, $b$, $k$) which has half-integer coefficients,
but nonetheless evaluates to an integer on all integer inputs. For example, $q^{\binom{d+1}{2}} = p^{\frac{-3}{2}(d^2 + d)}$. These quadratic polynomials are best written using
binomial coefficients when writing computer algebra code, because exponents are required to be integers.

Finally, in \S\ref{sec:closedformularou} we treat the case where $z$ is specialized to a primitive $3m$-th root of unity $\ze$. For a Laurent polynomial $g$ in $z$ we write $g(\ze)$ for
its specialization. We focus on the case where $a+b+1 = 3m$, because the Frobenius trace has degree\footnote{In general the Frobenius trace has degree $-m \binom{n}{2}$. One can see that the staircase monomial has degree $m \binom{n}{2}$. It is a
misleading coincidence that $\binom{n}{2} = n$ when $n=3$.} $-3m$. We also focus on the case $k=2m$. This is the action of Demazure operators on the staircase monomial $\stair =
x_1^{2m} x_2^m x_3^0$.

\begin{notation} Let $\xi_m(a,i) := \Xi(a,3m-a-1,i,2m)(\ze)$. \end{notation}

Incredibly, the messy formulas for $\Xi(a,b,i,k)$ simplify dramatically after specialization to a root of unity. Loosely speaking, here is the formula for $\xi_m(a,i)$. Depending on the triple $(a,b,i)$ with $a+b+1 = 3m$, there are integers $c$ and $t$, with $t$ close to $\frac{m}{2}$, for which
\begin{equation} \label{eq:loose} \xi_m(a,b,i) = \ze^{\foobar} m^2 \binom{t}{c}_{\ze^3}.\end{equation}
Here $\binom{t}{c}_{\ze^3}$ represents an unbalanced quantum binomial coefficient in
the variable $\ze^3$. The precise formula for the integers $\foobar$ and $c$ and $t$ in terms of the triple $(a,b,i)$ is annoying and depends on various parity considerations.

Our formula \eqref{eq:loose} is zero for $c < 0$ or $c > t$. This corresponds to $a < m-1$ or $a > 2m$. Indeed, the vanishing of $\pa_{(a,3m-a-1,i)}$ in $\NC(m,m,3)$ when $a < m-1$ or
$a > 2m$ is already a consequence of the roundabout relations mentioned in Remark \ref{rmk:roundabout}.

\begin{ex} \label{ex:dangitscomplex} In this computer-based example, we illustrate just how surprising it is that $\xi_m$ admits a simple formula. Begin with $\stair$, and apply Demazure operators one at a time, to see if there are any patterns that emerge. It would take pages to print this example, so we instead point the reader to our code \cite[\texttt{ExploringDemazureOperators.m}]{EJYcode}, which can be executed on MAGMA's free online interface\footnote{Again, it is found here: \texttt{http://magma.maths.usyd.edu.au/calc/}}. This program will print the polynomial one obtains after applying the first $j$ Demazure operators, though for reasons of clarity it is printed as a list of terms. For reasons of space, we have ignored any terms in the polynomial which are divisible by $x_1 x_2 x_3$; these do not affect the final answer thanks to Lemma \ref{lem:killx123}.

Also for sanity, we provide a table of how MAGMA prints the powers of $\ze$, since these powers are easily obfuscated by the way MAGMA treats cyclotomic fields. For small $m$ we note that the appearance of e.g. $4\ze^7 - \ldots$ in a coefficient implies that, when trying to write this coefficient as a linear combination of powers of $\ze$, at least four powers of $\ze$ are required.

We treat the case where $m = 5$ and $a = 7$ and $b = 7$ and $i=1$ in the code as written, noting that $a+b+1 = 3m$. We encourage the reader to vary $a$ or $m$ and see what happens.

For example, consider what happens after applying 8 or 9 Demazure operators from the word $\un{w}(a,b,i)$ to $\stair$. The coefficients in the resulting polynomial have the appearance of total chaos. However, after applying $2m$ Demazure operators, structure crystalizes seemingly from nowhere: all coefficients are multiples of $m$, and are powers of $z$ times quantum numbers. Things become more chaotic until finally, at the final ($3m$-th) step, the result is suddenly a multiple of $m^2$.

Similar behavior occurs for other values of $a$ and $b$, though the moment of first crystalization is after either $2m$ or $2m-1$ steps (depending on parity, it appears). The result depends strongly on $i$ until the final step! When $m$ is not prime, the various prime factors of $m$ start to divide coefficients at different times, which is quite mysterious!
\end{ex}

% \begin{defn} A \emph{cyclic word} is a word in $\Om$ of the form $(i, i+1, i+2, \ldots, i+(d-1))$ (clockwise) or of the form $(i,i-1, i-2, \ldots, i-(d-1))$ (widdershins) for some $i \in \Om$ and $d \ge 1$. We write $(i, i+1, \ldots, i+(d-1) = j)$ as $\cw_{i_L, d}$ when we wish to emphasize that it starts with $i$, and as $\cw_{j_R, d}$ when we wish to emphasize that it ends in $j$. Similarly, we write $(i,i-1, i-2, \ldots, i-(d-1) = j)$ as $\ws_{i_L, d}$ or $\ws_{j_R, d}$. \end{defn}

% \begin{ex} When $n=3$, $\cw_{2_L,5} = \cw_{3_R,5} = (2,3,1,2,3)$ is a clockwise cyclic word of length $5$. \end{ex}

%% file: Symmetry.tex
%!TEX root = QFrob2.tex

%%%%%%%%%%%%%%%%%%%%%%%%%%%%%%%%%%%%%%%%%%%%%%%%%%%%%%%%%
%========================================================
\section{Symmetries and recursions} \label{sec:symmetry}
%========================================================
%%%%%%%%%%%%%%%%%%%%%%%%%%%%%%%%%%%%%%%%%%%%%%%%%%%%%%%%%

Throughout this chapter, $\ell$ will always equal $a+b+1$ unless otherwise stated.

The goal of this chapter is to reduce the amount of work we need to do to both state and prove a formula for the scalars $\Xi(a,b,i,k)$. Ultimately, we hope to explain the following table.
\begin{equation} \label{benstable} \begin{tikzpicture}
  \matrix [column sep={3cm,between origins}, row sep = {1cm,between origins}, text width = 2cm, text height=0.6 cm]
  {
                        & \node {$i=1$};  & \node {$i=2$};& \node{$i=3$};\\
    \node {$k=0$};      & \node(c10){};   & \node(c20){}; & \node(c30){};\\
    \node {$k=1$};      & \node(c11){};   & \node(c21){}; & \node(c31){};\\
                        & \node(c12){};   & \node(c22){}; & \node(c32){};\\
    \node {\vdots};     & \node(c13){};   & \node(c23){}; & \node(c33){};\\
                        & \node(c14){};   & \node(c24){}; & \node(c34){};\\
                        & \node(c15){};   & \node(c25){}; & \node(c35){};\\
    \node {$k=\ell-1$}; & \node(c16){};   & \node(c26){}; & \node(c36){};\\
    \node {$k=\ell$};   & \node(c17){};   & \node(c27){}; & \node(c37){};\\
  };
  \node[fit=(c10), draw,inner sep=0pt](a1){edge case};
  \node[fit=(c11)(c13), draw,inner sep=0pt](c1){do it};
  \node[fit=(c14)(c16), draw,inner sep=0pt,fill=white!50!gray](c2){don't};
  \node[fit=(c17),inner sep=0pt, fill=white!50!gray,draw](a2){edge case};

  \node[fit=(c20),draw, fill=white!50!gray,inner sep=0pt](z1) {edge case};
  \node[fit=(c21),draw, inner sep=0pt](z3){special};                             
  \node[fit=(c22)(c26),draw, inner sep=0pt](z5){bulk}; 
  \node[fit=(c27),draw, fill=white!50!gray, inner sep=0pt] (z7) {zero};       

  \node[fit=(c30),inner sep=0pt,draw, fill=white!50!gray] (z8) {zero};
  \node[fit=(c31)(c35),inner sep=0pt, draw, fill=white!50!gray](z6){bulk};
  \node[fit=(c36),inner sep=0pt,draw, fill=white!50!gray](z4)(z4){special};
  \node[fit=(c37),inner sep=0pt,draw, fill=white!50!gray](a4)(z2){edge case};

  \draw [<->,red,thick] (a1.south west) to [bend right = 20] node [midway] {} (a2.north west);
  \draw [<->,orange,thick] (a2.north east) to node [midway] {} (z1.south west);
 % \draw [<->,red,thick] (a1.south east) to [bend left = 10] node [midway] {} (a4.north west);
  \draw[<->,blue,thick] (z1.east) to node [midway]   {} (z2.west);
  \draw[<->,blue,thick] (z3.east) to node [midway]   {} (z4.west);
  \draw[<->,blue,thick] (z5.east) to node [midway]   {} (z6.west);
  \draw[<->,green,thick] (z7.east) to node [midway]   {} (z8.west);
  \draw[<->,red,thick] (c1.west) to [bend right=10] node [midway]   {} (c2.west);
\end{tikzpicture}
\end{equation}
The arrows in this table represent symmetries that we define in \S\ref{ssec:exploit}, which relate $\Xi(a,b,i,k)$ to $\Xi(a,b,i',k')$, for a fixed value of $a$ and $b$. The color coding of arrows is explained in Remark \ref{rmk:colors}. The gray-shaded entries can all be determined from the white-shaded entries, and we verify that our symmetries do not provide any additional constraints between white-shaded entries.

There is one additional symmetry \eqref{sitauonXi} not pictured above, because it changes the values of $a$ and $b$. More precisely, it relates $\Xi(c,0,i,k)$ with
$\Xi(0,c,i,\ell-k)$. Due to this symmetry we need only examine the case when $a, b > 0$, the case when $a>0$ and $b=0$, and the base case where $a=b=0$. See Lemma
\ref{lem:aftersymmetry} for details.

Our formulas for $\Xi(a,b,i,k)$ take one form when $k \notin \{0,\ell\}$ and a different form when $k \in \{0,\ell\}$, which is why those are labeled as edge cases. Note that $\Xi(a,b,2,\ell) = \Xi(a,b,3,0) = 0$.

We also prove recursive formulas which relate $\Xi(a,b,i,k)$ with $\Xi(a',b',i',k')$ where $a' + b' = a + b - 1$. That is, the length $\ell' = a' + b' + 1$ is decreased
by one. The recursive formulas for $\Xi(a,b,i,k)$ are compatible with symmetry, so that proving the recursive formulas for the white-shaded entries is sufficient to imply them all.
These recursive formulas have a different form when $k \in \{0,\ell\}$. More interestingly, they also have a different form in the boxes labeled ``special.'' These are cases where
$k \notin \{0,\ell\}$ but $k' \in \{0,\ell'\}$. See Lemma \ref{lem:recursive} for details.

%========================================================
\subsection{Exploiting Symmetries} \label{ssec:exploit}
%========================================================

In Definition \ref{defn:symmetries} we introduced two symmetries $\si$ and $\tau$. The goal of this section is to prove four key consequences.

\begin{notation} For $f \in \ZZ[z,z^{-1}]$, let $\overline{f}$ be obtained by the $\ZZ$-linear automorphism which swaps $z$ and $z^{-1}$. \end{notation}

\begin{prop} For all $a, b \ge 0$ and all $0 \le k \le \ell = a+b+1$, we have
\begin{subequations} \label{allsymmetries}
\begin{equation} \label{Xi1symmetry} \Xi(a,b,1,k) = -z^{2k-\ell} \Xi(a,b,1,\ell-k), \end{equation}
\begin{equation} \label{Xi23symmetry} \Xi(a,b,2,k) = -z^{\ell-k} \Xi(a,b,3,\ell-k). \end{equation}
For all $i \in \Om$ we also have
\begin{equation} \label{sionXi} \Xi(a,b,i,\ell) = \Xi(a,b,i+1,0), \end{equation}
\begin{equation} \label{sitauonXi} \overline{\Xi(c,0,i,k)} = (-z)^{\ell} \Xi(0,c,-i-1,\ell-k). \end{equation}
\end{subequations}
In the last equation, $\ell = c+1$.
\end{prop}

\begin{proof}
From the definition of the divided difference operator we have
\begin{equation} \label{paantiinv} \pa_i(f) = - \pa_i(s_i(f)). \end{equation}
Let us note that
\begin{equation} \label{actionofs1} s_i(x_i^k x_{i+1}^{d - k}) = z^{2k-d}x_i^{d-k} x_{i+1}^k. \end{equation}	
Using \eqref{paantiinv} and \eqref{actionofs1} for $i=1$ we have
\begin{equation} \Xi(a,b,1,k) = \pa_{(a,b,1)}(x_1^k x_2^{\ell-k}) = - \pa_{(a,b,1)}(z^{2k-\ell} x_1^{\ell-k} x_2^k) = -z^{2k-\ell} \Xi(a,b,1,\ell-k), \end{equation}
which is \eqref{Xi1symmetry}.

Recall that $\si$ fixes all scalars in $\ZZ[z,z^{-1}]$, and that all values of $\Xi$ are scalars. Similar to the above we have 
\begin{align} \nonumber \Xi(a,b,2,k) = \pa_{(a,b,2)}&(x_1^k x_2^{\ell - k}) = \pa_{(a,b,3)}(x_2^k x_3^{\ell-k}) \\ & = -z^{\ell-k} \pa_{(a,b,3)}(x_2^k x_1^{\ell-k}) = -z^{\ell-k} \Xi(a,b,3,\ell-k). \end{align}
The first and last equalities are by definition, the second equality uses \eqref{sionabi}, and the third equality uses \eqref{paantiinv} and \eqref{actionofs1} for $i=3$. This proves \eqref{Xi23symmetry}.

Another application of \eqref{sionabi} is 
\begin{equation} \Xi(a,b,i,\ell) = \pa_{(a,b,i)}(x_1^{\ell}) = \si(\pa_{(a,b,i)}(x_1^{\ell})) = \pa_{(a,b,i+1)}(x_2^{\ell}) = \Xi(a,b,i+1,0), \end{equation}
which proves \eqref{sionXi}.

For the final symmetry \eqref{sitauonXi} we use the symmetry $\si^{-1} \tau$. Recall that $\tau$ acts on scalars as $\tau(f) = \overline{f}$. One can also compute that $\tau(x_1^k x_2^{\ell - k}) = x_3^k x_2^{\ell - k}$. Thus $\si^{-1} \tau$ acts on scalars as $f \mapsto \overline{f}$, and sends $x_1^k x_2^{\ell -k}$ to $x_1^{\ell - k} x_2^k$. Using this and \eqref{sionabi} and  \eqref{tauonabi}, we get
\begin{align}\nonumber \overline{\Xi(c,0,i,k)} & = \si^{-1} \tau (\pa_{(c,0,i)}(x_1^k x_2^{\ell-k})) \\ & = (-z)^{\ell} \pa_{(0,c,-i-1)}(x_1^{\ell - k} x_2^k) = (-z)^{\ell} \Xi(0,c,-i-1,\ell-k), \end{align}
which proves \eqref{sitauonXi}. \end{proof}

\begin{cor} If $\ell$ is even then
\begin{equation} \Xi(a,b,1,\ell/2) = 0. \end{equation}
\end{cor}

\begin{proof} This follows from \eqref{Xi1symmetry}, or from $s_1$-invariance: $\pa_1(x_1^{\ell/2} x_2^{\ell/2}) = 0$. \end{proof}

In \S\ref{sec:closedformula} we will define formulas for $\Xi(a,b,i,k)$. We check that they satisfy \eqref{allsymmetries}, which reduces the amount of work we need to do to verify their correctness. Indeed, we even use \eqref{allsymmetries} to define our formulas in many cases. For example, we use \eqref{sitauonXi} to derive the
formula for $\Xi(a,b,i,k)$ when $a=0$ from the case when $b=0$. However, the fact that the resulting formula satisfies \eqref{allsymmetries} is not completely tautological. In theory, the symmetries might have ``monodromy,'' that is, we might be able to follow a chain of symmetries in a loop and deduce an unexpected consequence. We now argue that this is not the case.

The symmetries \eqref{Xi23symmetry} and \eqref{sionXi} seem to provide two different constants of proportionality between $\Xi(a,b,2,\ell)$ and $\Xi(a,b,3,0)$. There is no contradiction
here, as $\Xi(a,b, 2, \ell) = \Xi(a,b,3,0) = 0$, see Example \ref{ex:recursion2edgecases}. Similarly, \eqref{Xi1symmetry} implies that $\Xi(a,b,1,\ell/2)=0$ when $\ell$ is even, as discussed in the previous corollary.

The symmetries \eqref{Xi1symmetry} and \eqref{sionXi} and \eqref{Xi23symmetry} produce a consistent loop \begin{equation} \label{secretsymmetry} \Xi(a,b,2,0) = \Xi(a,b,1,\ell) =
-z^{\ell} \Xi(a,b,1,0) = -z^{\ell} \Xi(a,b,3,\ell) = \Xi(a,b,2,0). \end{equation} One can confirm that \eqref{secretsymmetry} is also consistent with \eqref{sitauonXi} when $a = 0$ or
$b=0$. This handles all cases where $k \in \{0,\ell\}$.

When $k \notin \{0,\ell\}$ the symmetry \eqref{sionXi} is not present, and applying a combination of \eqref{sitauonXi} and either \eqref{Xi1symmetry} or \eqref{Xi23symmetry}, one can only derive a consistent loop of equalities. For example
\begin{align} \nonumber \Xi(c,0,1,k) = -z^{2k-\ell} & \Xi(c,0,1,\ell-k) = \overline{z^{-2k} \Xi(0,c,1,k)} \\ & = \overline{-z^{-\ell} \Xi(0,c,1,\ell-k)} = \Xi(c,0,1,k). \end{align}

\begin{rem} \label{rmk:colors} In the table of \eqref{benstable}, \eqref{Xi1symmetry} is drawn in red, and \eqref{Xi23symmetry} is drawn in blue. When $i=1$, \eqref{sionXi} is drawn in orange. Green represents the overlap of both \eqref{Xi23symmetry} and \eqref{sionXi} for $i=2$, which is consistent since both values of $\Xi$ are zero. When $i=3$, \eqref{sionXi} not pictured, but it agrees with the composition of a red, orange, and blue arrow by \eqref{secretsymmetry}. The fourth symmetry \eqref{sitauonXi} changes the values of $a$ and $b$, so is not described by \eqref{benstable}. \end{rem}

Our previous discussion is summarized in the following lemma, c.f. the white-shaded regions in \eqref{benstable}.

\begin{lem} \label{lem:aftersymmetry} The following computations determine all values of $\Xi(a,b,i,k)$ after using the symmetries \eqref{sionXi}, \eqref{sitauonXi}, \eqref{Xi1symmetry}, and \eqref{Xi23symmetry}. \begin{itemize}
\item The \emph{base cases}: $\Xi(0,0,i,k)$ for $0 \le k \le 1$,
\item The cases $\Xi(a,b,1,k)$ for $a > 0$ and $b \ge 0$ and $\ell/2 \le k$,
\item The cases $\Xi(a,b,2,k)$ for $a > 0$ and $b \ge 0$ and $k > 0$.
\end{itemize}
Moreover, suppose one specifies these values of $\Xi(a,b,i,k)$, and uses various symmetries to determine the other values of $\Xi$. The result is consistent, and satisfies the four symmetries, if and only if $\Xi(a,b,2,\ell) = 0$, and $\Xi(a,b,1,\ell/2)=0$ when $\ell$ is even. \end{lem}

\subsection{Base cases} \label{ssec:basecase}
%========================================================

We now discuss the case when $a = b = 0$. Either $k=0$ or $k = 1 = \ell$. In this case, $\pa_{(0,0,i)} = \pa_i$, and we already know the values of $\pa_i(x_j)$. From this we deduce the following theorem, an immediate calculation using the definition of the divided difference operators.

\begin{thm} \label{thm:basecase} We have
\begin{subequations} \label{recursiveformulabase}
\begin{equation} \label{recursiveformula2base} \Xi(0,0,2,0) = \pa_2(x_2) = 1, \qquad \Xi(0,0,2,1) = \pa_2(x_1) = 0. \end{equation}
\begin{equation} \label{recursiveformula3base} \Xi(0,0,3,0) = \pa_3(x_2) = 0, \qquad \Xi(0,0,3,1) = \pa_3(x_1) = -z. \end{equation}	
\begin{equation} \label{recursiveformula1base} \Xi(0,0,1,0) = \pa_1(x_2) = -z, \qquad \Xi(0,0,1,1) = \pa_1(x_1) = 1. \end{equation}
\end{subequations}
\end{thm}

%========================================================
\subsection{Recursive formulas} \label{ssec:recursive}
%========================================================

By applying Demazure operators one at a time, one can derive recursive formulas for $\Xi(a,b,i,k)$, which reduce the total length $a+b+1$ by one. These recursive formulas have many edge cases, which become evident when trying to derive them. Let us explore the recursive formulas when $i=2$. 

\begin{ex} \label{ex:recursion2} Suppose that $b > 0$. Then $\un{w}(a,b,2) = \un{w}(a,b-1,3) \circ (2)$. Note that 
\begin{equation} \label{manysurvivors} \pa_2(x_1^k x_2^{\ell-k}) = x_1^k \pa_2(x_2^{\ell-k}) =  x_1^k \sum_{j = 0}^{\ell-k-1} x_2^j (z x_3)^{\ell - k - 1 - j}. \end{equation}
When examining this innocent-seeming formula there are a surprising number of edge cases, which we leave for later. For now, assume that $0 < k < \ell - 1$.

A great simplification occurs when working modulo $x_1 x_2 x_3$, which we can do thanks to Lemma \ref{lem:killx123}. Only two terms in \eqref{manysurvivors} survive:
\begin{equation} \label{twosurvivors} \pa_2(x_1^k x_2^{\ell - k}) \equiv x_1^k x_2^{\ell - k - 1} + z^{\ell - k - 1} x_1^k x_3^{\ell - k -1}. \end{equation}
So
\begin{equation} \Xi(a,b,2,k)  = \pa_{\un{w}(a,b-1,3)}(\pa_2(x_1^k x_2^{\ell - k})) =  \pa_{\un{w}(a,b-1,3)}(x_1^k x_2^{\ell - 1 - k} + z^{\ell - k - 1} x_1^k x_3^{\ell - k -1}). \end{equation}
Now $\pa_{\un{w}(a,b-1,3)}(x_1^k x_2^{\ell - 1 - k}) = \Xi(a,b-1,3,k)$ by definition. To analyze the other term, we apply $\si$. Since $\si$ fixes scalars, we deduce that
\begin{equation} \pa_{\un{w}(a,b-1,3)}(x_1^k x_3^{\ell - k -1}) = \pa_{\un{w}(a,b-1,1)}(x_2^k x_1^{\ell-1-k}). \end{equation}
Thus
\begin{equation} \Xi(a,b,2,k) = \Xi(a,b-1,3,k) + z^{\ell - k - 1} \Xi(a,b-1,1,\ell-k-1). \end{equation}
Using \eqref{Xi1symmetry} we have our preferred recursive formula 
\begin{equation} \label{itsarecursionyay} \Xi(a,b,2,k) = \Xi(a,b-1,3,k) - z^{2\ell - 3k - 2} \Xi(a,b-1,1,k), \quad \text{if } a \ge 0, b > 0, 0 < k < \ell-1. \end{equation}
\end{ex}
	
\begin{ex} \label{ex:recursion2edgecases} Now let us consider some of the edge cases when examining $\Xi(a,b,2,k)$ for $b > 0$. If $k = \ell$ then $\pa_2(x_1^k x_2^{\ell - k}) = \pa_2(x_1^\ell) = 0$. Indeed, more generally we have
\begin{equation} \label{kislitszero} \Xi(a,b,2,\ell) = 0, \quad \text{for all } a, b \ge 0. \end{equation}

When $k = \ell - 1$, we continue as in Example \ref{ex:recursion2}, working modulo $x_1 x_2 x_3$, except instead of two surviving terms as in \eqref{twosurvivors} we only have one:
\begin{equation} \pa_2(x_1^{\ell-1} x_2) = x_1^{\ell-1}. \end{equation}
So
\begin{equation} \label{Xi2specialist} \Xi(a,b,2,\ell-1) = \Xi(a,b-1,3,\ell-1). \end{equation}

When $k = 0$, all the terms from \eqref{manysurvivors} survive modulo $x_1 x_2 x_3$, and there is no great simplification. Following the remaining arguments of Example
\ref{ex:recursion2} we deduce that \begin{equation} \label{notaneasyrecursion} \Xi(a,b,2,0) = \sum_{j=0}^{\ell - 1} z^{\ell -j - 1} \Xi(a,b-1,1,\ell-j-1) = -\sum_{j=0}^{\ell - 1} z^{2\ell -3j - 2}
\Xi(a,b-1,1,j). \end{equation}
However, we can alternatively apply \eqref{sionXi} to deduce that
\begin{equation} \Xi(a,b,2,0) = \Xi(a,b,1,\ell). \end{equation}
Thus turns \eqref{notaneasyrecursion} into a formula for $\Xi(a,b,1,\ell)$ in terms of various $\Xi(a,b-1,1,j)$. Our general recursive formula for $\Xi(a,b,1,k)$ will be very similar in form to \eqref{notaneasyrecursion}. \end{ex}

% \begin{rem} The edge cases $k=\ell$ and $k= 0$ correspond to monomials $x_1^{\ell}$ and $x_2^{\ell}$, which have additional symmetry over the generic monomial. It is not so surprising that the formula for $\Xi(a,b,i,k)$ obeys one pattern when $0 < k < \ell$, and a separate pattern when $k \in \{0,\ell\}$. Nor is it surprising that the recursive formulas which govern $\Xi(a,b,i,k)$ when $k \in \{0,\ell\}$ would be different either.
%
% It may seem more surprising that the case $k = \ell-1$ gave rise to a special recursive formula, though this also makes perfect sense. While $k = \ell-1$ is not an edge case for
% $\Xi(a,b,i,k)$, it is an edge case for $\Xi(a,b-1,i,k)$. Thus our recursive formula compares a generic case to an edge case, and can be
% expected to have a special form. \end{rem}

We continue our exploration of recursive formulas by considering the $b=0$ case.
	
\begin{ex} \label{ex:recursion2b0} When $b=0$, one can derive recursive formulas for $\Xi(a,0,2,k)$ in almost exactly the same way. The key difference is that $\un{w}(a,0,2) = \un{w}(a-1,0,1) \circ (2)$. When simplifying
\[ \pa_{\un{w}(a-1,0,1)}(z^{\ell - k - 1} x_1^k x_3^{\ell - k -1}) \]
we choose to first apply the formula $\pa_1(s_1(f)) = - \pa_1(f)$ to obtain
\[ -\pa_{\un{w}(a-1,0,1)}(z^{\ell - 1} x_2^k x_3^{\ell - k -1}), \]
and then apply $\si^{-1}$ to obtain
\[ -\pa_{\un{w}(a-1,0,3)}(z^{\ell - 1} x_1^k x_2^{\ell - k -1}). \]
Ultimately, instead of \eqref{itsarecursionyay} one has
\begin{equation} \Xi(a,0,2,k) = \Xi(a-1,0,1,k) - z^{\ell-1} \Xi(a-1,0,3,k), \quad \text{if } a > 0, 0 < k < \ell-1. \end{equation}
	
The cases $k \in \{0, \ell-1, \ell\}$ have their own recursive formulas. The case $k = \ell$ is covered by \eqref{kislitszero}, and the case $k=0$ will be treated by symmetry. When $k = \ell-1$ we have $\pa_2(x_1^{\ell-1} x_2) = x_1^{\ell-1}$ as before, so
\begin{equation} \label{Xi2b0specialist} \Xi(a,0,2,\ell-1) = \Xi(a-1,0,1,\ell-1). \end{equation}
\end{ex}

For pedagogical reasons, let us briefly discuss recursive formulas for $\Xi(a,b,3,k)$. For brevity, we focus on a special case.

\begin{ex} Consider $\Xi(a,b,3,1)$ when $b > 0$, which is $\pa_{(a,b,3)}(x_1 x_2^{\ell-1})$. Note that $\pa_3(x_1 x_2^{\ell-1}) = -z^{-1} x_2^{\ell-1}$. Since $\un{w}(a,b,3) = \un{w}(a,b-1,1) \circ (3)$, we deduce that
\begin{equation} \label{Xi3specialist} \Xi(a,b,3,1) = -z^{-1} \Xi(a,b-1,1,0). \end{equation}
	
On the other hand, we could have used the symmetry \eqref{Xi23symmetry} instead, to write
\[ \Xi(a,b,3,1) = -z^{-1} \Xi(a,b,2,\ell-1),\]
relating this special case for $i=3$ to the earlier special case for $i=2$. Continuing with \eqref{Xi2specialist} and \eqref{sionXi} we get
\begin{equation} \Xi(a,b,3,1) = -z^{-1} \Xi(a,b,2,\ell-1) = -z^{-1} \Xi(a,b-1,3,\ell-1) = -z^{-1} \Xi(a,b-1,1,0). \end{equation}
Thus we end up deducing \eqref{Xi3specialist} from \eqref{Xi2specialist} and symmetry.  \end{ex}

The point we wish to make by the previous example is this. There are many recursive formulas for the values of $\Xi(a,b,i,k)$, which one can compute directly. There are also
symmetries which relate various values of $\Xi(a,b,i,k)$. The recursive formulas and the symmetries are compatible: applying the symmetries to one recursive formula, one obtains the
other recursive formula. This is ultimately a tautological statement because the way in which the formulas are derived is compatible with the symmetries; the proof is unenlightening
and is left as a tedious exercise. This is the intuitive answer, though it takes some thought to puzzle through why, and we prefered to convince the reader with example rather than proof.

Rather than list all the recursive formulas, we provide a minimal list that determines the remaining formulas via symmetry, c.f. \eqref{benstable}.

\begin{lem} \label{lem:recursive} The following recursive formulas hold for $\Xi(a,b,i,k)$. Moreover, $\Xi(a,b,i,k)$ is the unique function satisfying these recursive formulas, the base cases \eqref{recursiveformulabase}, and the symmetries \eqref{sionXi}, \eqref{sitauonXi}, \eqref{Xi1symmetry}, and \eqref{Xi23symmetry}.

\begin{subequations} \label{recursiveformulamain}
\begin{equation} \label{recursiveformula2} \Xi(a,b,2,k) = \Xi(a,b-1,3,k) - z^{2\ell - 3k - 2} \Xi(a,b-1,1,k) \quad \text{ if } a,b > 0, 0 < k < \ell-1, \end{equation}
\begin{equation} \label{recursiveformula1} \Xi(a,b,1,k) = \sum_{c = \ell - k}^{k-1} z^{k-1-c} \Xi(a,b-1,2,c) \quad \text{ if } a,b > 0, \ell/2 \le k \le \ell, \end{equation}
\begin{equation} \label{recursiveformula2b} \Xi(a,0,2,k) = \Xi(a-1,0,1,k) - z^{\ell-1} \Xi(a-1,0,3,k) \quad \text{ if } a > 0, 0 < k < \ell-1, \end{equation}
\begin{equation} \label{recursiveformula1b} \Xi(a,0,1,k) = \sum_{c = \ell - k}^{k-1} z^{k-1-c} \Xi(a-1,0,3,c) \quad \text{ if } a > 0, \ell/2 \le k \le \ell, \end{equation}
\begin{equation} \label{recursiveformula2l} \Xi(a,b,2,\ell) = 0, \end{equation}
\begin{equation} \label{recursiveformula2special} \Xi(a,b,2,\ell-1) = \Xi(a,b-1,3,\ell-1) \quad \text{ if } a,b > 0, \end{equation}
\begin{equation} \label{recursiveformula2bspecial} \Xi(a,0,2,\ell-1) = \Xi(a-1,0,1,\ell-1) \quad \text{ if } a>0.  \end{equation}
\end{subequations}

Note that \eqref{recursiveformula2} and \eqref{recursiveformula1} and \eqref{recursiveformula2special} hold when $a = 0$ as well, but these formulas are not needed for uniqueness.
\end{lem}

% \begin{rem} \BE{I say we comment this remark out} Here are some additional recursive formulas, for posterity.
% \begin{subequations}
% \begin{equation} \label{recursiveformula3} \Xi(a,b,3,k) = z^{\ell - 2k} \Xi(a,b-1,2,k-1) - z^{-1} \Xi(a,b-1,1,k-1), \quad \text{ if } b > 0, 1 < k < \ell, \end{equation}
% \begin{equation} \label{recursiveformula1prime} \Xi(a,b,1,k) = -\sum_{c = k}^{\ell-1-k} z^{k-1-c} \Xi(a,b-1,2,c), \quad \text{ if } b > 0, 0 \le k \le \ell/2, \end{equation}
% \begin{equation} \label{recursiveformula3b} \Xi(a,0,3,k) = z^{\ell - 3k + 1}\Xi(a-1,0,1,k-1) - z^{-1} \Xi(a-1,0,2,k-1), \quad \text{ if } b = 0, a > 0, k > 0, \end{equation}
% \begin{equation} \label{recursiveformula1bprime} \Xi(a,0,1,k) = - \sum_{c = k}^{\ell - 1 - k} z^{k-1-c} \Xi(a-1,0,3,c), \quad \text{ if } b = 0, a > 0, 0 \le k \le \ell/2, \end{equation}
% \begin{equation} \label{recursiveformula3c} \Xi(a,b,3,0) = 0, \end{equation}
% \begin{equation} \Xi(a,b,3,1) = -z^{-1} \Xi(a,b-1,1,0), \quad \text{ if } b > 0. \end{equation}
% \begin{equation} \Xi(a,0,3,1) = -z^{-1} \Xi(a-1,0,2,0), \quad \text{ if } a > 0. \end{equation}
% \end{subequations}
%
% To give another example of the compatibility of recursive formulas and symmetry: \eqref{recursiveformula2} says
% \[
% \Xi(0,b,2,k) = \Xi(0,b-1, 3, k) - z^{2(b+1)-3k-2}\Xi(0,b-1,1,k).
% \]
% Applying \eqref{sitauonXi} and applying the bar involution, this is equivalent to
% \[
% (-z)^{b+1} \Xi(b,0,3,b+1-k) = (-z)^b\Xi(b-1,0,2,b-k) - (-z)^bz^{-2(b+1)+3k+2}\Xi(b-1,0,1,b-k)
% \]
% which is equivalent to \eqref{recursiveformula3b}. \BE{comment until here}
% \end{rem}

\begin{proof} In Examples \ref{ex:recursion2} and \ref{ex:recursion2edgecases} and \ref{ex:recursion2b0}, we already proved all the recursive formulas in \eqref{recursiveformulamain} except \eqref{recursiveformula1} and \eqref{recursiveformula1b}. Both of these follow immediately from
\begin{equation} \label{itspa1} \pa_1(x_1^k x_2^{\ell-k}) = \sum_{c = \ell - k}^{k-1} z^{k-1-c} x_1^c x_2^{\ell-1-c}, \qquad \text{ if } \ell/2 \le k \le \ell. \end{equation} Note that $\ell/2 \le k$ is the same as $\ell - k \le k$. The formula \eqref{itspa1} is easier to prove after dividing by the $s_1$-invariant monomial $x_1^{\ell - k} x_2^{\ell - k}$.

The uniqueness of a solution to these equations is clear. We use induction on the length $a+b+1$, and compute $\Xi(a,b,i,k)$. Lemma \ref{lem:aftersymmetry} implies that, after applying symmetry, we are in one of the cases handled by \eqref{recursiveformulamain}, so $\Xi(a,b,i,k)$ is determined by values of $\Xi$ with smaller lengths.
\end{proof}

When we prove the correctness of our formulas in the next chapter, we will do so by verifying that they meet the conditions of Lemma \ref{lem:recursive}. Note that we have checked all cases of our formulas for small values of $\ell$ by computer, which obviates the need to check the edgiest of edge cases, like when $a=1$ and $b=0$.

%% file: NastyFour.tex
%!TEX root = QFrob2.tex

%%%%%%%%%%%%%%%%%%%%%%%%%%%%%%%%%%%%%%%%%%%%%%%%%%%%%%%%%
%========================================================
\section{A closed formula for scalars in the deformed affine nilHecke algebra} \label{sec:closedformula}
%========================================================
%%%%%%%%%%%%%%%%%%%%%%%%%%%%%%%%%%%%%%%%%%%%%%%%%%%%%%%%%

Our goal in this chapter is to state an explicit formula for the scalars $\Xi(a,b,i,k)$, which were introduced in \S\ref{sec:prelims}. Throughout this chapter, $\ell = a+ b + 1$
unless otherwise stated.

Thanks to the work done in \S\ref{sec:symmetry}, we need only specify the value of $\Xi(a,b,i,k)$ for certain quadruples $(a,b,i,k)$, see Lemma \ref{lem:aftersymmetry}. The
remaining values are determined by symmetry. However, because they obey nearly the same formulas and it makes the result easier to reference, we chose to provide formulas in some of
the redundant cases as well. For example, our formulas for $\Xi(a,b,1,k)$ when $k \notin \{0,\ell\}$ do not distinguish between $k > \ell/2$ and $k < \ell/2$, because the same
formula works in both cases.

The scalars $\Xi(a,b,i,k)$ follow a general formula whenever $a > 0$, $b > 0$, and $0 < k < \ell$. We call this the \emph{standard regime}, and we call the other possibilities \emph{edge cases}. The formula for edge cases is a variant on the formula in the standard regime, but is much simpler.

Sections \S\ref{ssec:formula} and \S\ref{ssec:formulatrue} are dedicated to the statement of the results, which are spread over several theorems which treat different cases. The
rest of the chapter is the proof of these theorems. In \S\ref{ssec:formulassymmetric} we confirm that the formulas we provide do indeed satisfy the four symmetries from
\S\ref{ssec:exploit}, a fact which is not trivial. We check that the formulas satisfy recursive formulas in \S\ref{ssec:recursionb0} and \S\ref{ssec:recursionmain}. As a result,
Lemma \ref{lem:recursive} will imply that our formulas are correct.

\subsection{Variations on quantum factorial} \label{ssec:rho}
%========================================================

\begin{defn} Let 
\begin{equation} \label{defrho} \rho(d) := (q-q^{-1})^d [d]! = \prod_{c=1}^d (q^c - q^{-c}). \end{equation}
By the standard convention for empty products, $\rho(0) = 1$. 

Let 
\begin{equation} \label{defrhoprime} \rho'(d) := q^{-\binom{d+1}{2}} \rho(d) = \prod_{c=1}^d (1 - q^{-2c}). \end{equation}
\end{defn}

\begin{rem} Up to a power of $q$, $\rho(d)$ is the size of $GL_n(\mathbb{F}_{q^2})$. \end{rem}

\begin{lem} \label{lem:diffofy} For any value of $c$ and $d$ we have $q^{-2c} - q^{-2d} = q^{-(c+d)} (q-q^{-1}) [d-c]$. \end{lem}

\begin{proof} Left to the reader. \end{proof}

%========================================================
\subsection{The formula in the standard regime} \label{ssec:formula}
%========================================================

From the numbers $a$ and $b$ we will extract dependent variables $\alpha$ and $\beta$ and $\nu$, and we state portions of the formula in terms of these new variables. We've also normalized in order to shunt the weirdest parts of the formula to $i=3$. Here is the crazy formula for $\Xi$ within the standard regime.

\begin{thm} \label{thm:mainnasty} Fix $a, b > 0$ and $i \in \Om$ and $0 < k < a+b+1 =: \ell$. Define $\alpha$ and $\beta$ by 
\begin{equation} \label{alphabetadef} a = 2\alpha + 2 \;\; \text{or} \;\; 2 \alpha + 1, \qquad b = 2 \beta + 2 \;\; \text{or} \;\; 2 \beta + 1, \end{equation}
depending on parity. Let 
\begin{equation} \label{nudef} \nu = \alpha + \beta + 2. \end{equation}
Then
\begin{equation} \label{eq:mainnasty} \Xi(a,b,i,k) = \thesign \cdot \gamma_1 \cdot \gamma_2 \cdot \gamma_3 \cdot \kappa_1 \cdot \kappa_2 \cdot \kappa_3 \cdot \lambda_1 \cdot \lambda_2 \cdot \lambda_3 \cdot \lambda_4 \cdot \lambda_5. \end{equation}
The three most interesting factors $\gamma$ are defined by 
\begin{subequations}
\begin{equation} \gamma_1 = \rho'(\alpha) \rho'(\beta), \end{equation}
\begin{equation} \gamma_2 = \begin{cases}
	1 & \text{ if $a$ and $b$ are odd} \\
	q^{-(\nu)} (q-q^{-1}) [k-\nu] & \text{ if $a$ is even and $b$ is odd } \\
	q^{-(\nu)} (q-q^{-1}) [\nu-k] & \text{ if $a$ is odd and $b$ is even and } i=1 \\ 
	q^{-(\ell-1)} (q-q^{-1}) [\ell-1-k] & \text{ if $a$ is odd and $b$ is even and } i=2 \\
	q^{-(\ell - 1)} (q-q^{-1})[1-k] & \text{ if $a$ is odd and $b$ is even and } i=3 \\
	q^{-(\ell)} (q-q^{-1})^2 [k-\nu][\nu+1-k] & \text{ if $a$ and $b$ are even and } i=1 \\
	q^{-(\ell+\nu-1)} (q-q^{-1})^2 [k-\nu][\ell-1-k] & \text{ if $a$ and $b$ are even and } i=2 \\
	q^{-(\ell+\nu-1)} (q-q^{-1})^2 [k-1][\nu+1 - k] & \text{ if $a$ and $b$ are even and } i=3 	
	 \end{cases} \end{equation}
\begin{equation} \gamma_3 = \begin{cases}
	\magic(\nu,k,\beta,-1) & \text{ if $a$ and $b$ are odd} \\
	\magic(\nu,k,\beta,0) & \text{ if $a$ is even and $b$ is odd } \\
	\magic(\nu,k,\beta,0) & \text{ if $a$ is odd and $b$ is even and } i=1 \\ 
	\magic(\nu,k,\beta,-1) & \text{ if $a$ is odd and $b$ is even and } i=2 \\
	q^{\beta} \magic(\nu,k-1,\beta,-1) & \text{ if $a$ is odd and $b$ is even and } i=3 \\
	\magic(\nu,k,\beta,+1) & \text{ if $a$ and $b$ are even and } i=1 \\
	\magic(\nu,k,\beta,0) & \text{ if $a$ and $b$ are even and } i=2 \\
	q^{\beta} \magic(\nu,k-1,\beta,0) & \text{ if $a$ and $b$ are even and } i=3 	
	 \end{cases}	 
\end{equation}
\end{subequations}
We have
\begin{equation} \thesign = (-1)^{\beta + k}. \end{equation}
The remaining factors are powers of $p$. The factors $\kappa$ depend on the value of $k$.
\begin{subequations}
\begin{eqnarray} 
	\kappa_1 & = & z^k q^{k(k - \beta - \ell)}, \\
	\kappa_2 & = & \begin{cases} q^{2k} & \text{ if } i = 2, \\ 1 & \text{ else}. \end{cases}
\end{eqnarray}
\end{subequations}
The remaining factors $\lambda$ only depend on $(a,b,i)$. Let $\phi$ denote the number of even elements in the pair $(a,b)$. Then $\ell = 2 \alpha + 2 \beta + 3 + \phi = 2 \nu - 1 + \phi$.
\begin{subequations}
\begin{eqnarray} 
	\lambda_1 & = & z^{\binom{\beta}{2}} z^{-\binom{\ell+1}{2}} p^{-(\beta+1)(\ell+3\beta)}, \\
	\lambda_2 & = & \begin{cases}
	z^{\beta} & \text{ if $a$ is odd }  \\
	z^{-\beta-3} & \text{ if $a$ is even, }
	\end{cases} \\
	\lambda_3 & = & \begin{cases}
	z^{\ell+1} & \text{ if $b$ is odd }  \\
	1 & \text{ if $b$ is even, }
	\end{cases} \\
	\lambda_4 & = & p^{(\beta+3)(\phi-1)}, \\
	\lambda_5 & = & \begin{cases}
	1 & \text{ if } i=1, \\
	z^{\ell} & \text{ if } i=2, \\
	z^{-\ell} & \text{ if } i =3, \end{cases}
\end{eqnarray}
\end{subequations}
\end{thm}

\begin{rem} All these scalars, like $\gamma_2$ or $\lambda_5$, are functions of the quadruple $(a,b,i,k)$, and we will write $\gamma_2(a,b,i,k)$ when the quadruple varies or is not yet fixed by the context. \end{rem}

\begin{rem} \label{rem:whenqis1} Suppose that $q$ is specialized to $1$. Then $\gamma_1 = 0$ unless $\alpha = \beta = 0$, and $\gamma_2 = 0$ unless $a$ and $b$ are odd, so the only
nonzero possibility in the standard regime is $a = b = 1$. Up to symmetry, $\un{w}(1,1)$ is the longest element of the finite Weyl group.
\end{rem}

%========================================================
\subsection{The formula for edge cases} \label{ssec:formulatrue}
%========================================================

We also need to give formulas when $a=0$, when $b=0$, when $k=0$, and when $k = a+b+1$, and for combinations of these conditions. We've already treated the base cases $a=b=0$ in Theorem \ref{thm:basecase}. Many of the remaining cases come from symmetry. According to \eqref{benstable} and Lemma \ref{lem:aftersymmetry}, we need to treat the following three cases:
\begin{itemize} \item $a, b > 0$ and $k = \ell$,
	\item $a>0$, $b = 0$, and $k = \ell$. 
	\item $a>0$, $b = 0$, and $k \notin \{0,\ell\}$,
\end{itemize}
As it turns out, one formula suffices for $a>0$ and $b \ge 0$ and $k = \ell$.

For purposes of comparison it helps to use some of the same notation as in the standard regime, i.e. the variables $\alpha$, $\beta$, and $\nu$. However, it sometimes helps to adjust the definition of $\alpha$ or $\beta$, as parity considerations apply differently in edge cases.

\textbf{Case I:} $a > 0$ and $b \ge 0$ and $k= \ell$. The case $k=0$ is included via symmetry.

We are computing $\Xi(a,b,i,\ell) = \pa_{(a,b,i)}(x_1^{\ell})$. This is zero when $i=2$, since $x_1^{\ell}$ is invariant under $s_2$. What is less obvious is that the result is also zero when $b$ is odd. This is because $\un{w}(a,b,i)$ has right descent set $\{i,i+1\}$. Just as the operator $\pa_i$ kills any polynomial invariant under the reflection $s_i$, the operator $\pa_i \pa_{i+1} \pa_i$ kills any polynomial invariant under any of the reflections $s_i$ or $s_{i+1}$ or $s_i s_{i+1} s_i$.

For use in this case and other cases with $k \in \{0,\ell\}$ we define
\begin{equation} \nabla = \begin{cases}
	1 & \text{ if $i = 1$,} \\
	0 & \text{ if $i = 2$,} \\
	-z^{-\ell} & \text{ if $i = 3$,} \end{cases} \end{equation}

In this case we will adjust the definition of $\alpha$.

\begin{thm} \label{thm:klencase}  Fix $a > 0$ and $b \ge 0$ and $i \in \Om$ and let $\ell := a+b+1$. Define $\alpha$ (anew) and $\beta$ by 
\begin{equation} \label{alphabetadefklen} a = 2\alpha \;\; \text{or} \;\; 2 \alpha + 1, \qquad b = 2 \beta + 2 \;\; \text{or} \;\; 2 \beta + 1. \end{equation}
Then $\Xi(a,b,i,\ell) = \Xi(a,b,i+1,0)$ satisfies
\begin{equation} \label{kleneasyformula} \Xi(a,b,i,\ell) = \begin{cases} 0 & \text{if } b \text{ is odd,} \\
(-1)^{\beta+\ell} \nabla z^{-\binom{\ell}{2} + \binom{\beta+1}{2} + \ell(\beta+1)} \rho'(\alpha + \beta + 1) & \text{if } b \text{ is even.} \end{cases} \end{equation}
% Then $\Xi(a,b,i,\ell) = \Xi(a,b,i+1,0) = 0$ if $b$ is odd or $i=2$. Otherwise,
% \begin{equation} \label{kleneasyformula} \Xi(a,b,i,\ell) = \Xi(a,b,i+1,0) = (-1)^{\beta+\ell} \nabla z^{-\binom{\ell}{2}} \rho'(\alpha + \beta + 1) p^{(2 \ell + \beta)(\beta + 1)}. \end{equation}
When $b=0$ so that $\beta = -1$, this formula simplifies to
\begin{equation} \label{eq:bzerokleneasy} \Xi(a,0,i,\ell) = \Xi(a,0,i+1,0) = (-1)^{\ell-1} \nabla z^{-\binom{\ell}{2}} \rho'(\alpha). \end{equation}
\end{thm}

\begin{rem} When $i=3$, there is a simplification which absorbs $\nabla$ appropriately:
\begin{equation} \label{klen3} \Xi(a,b,3,\ell) = \Xi(a,b,1,0) = (-1)^{\beta+\ell+1} z^{-\binom{\ell+1}{2}} \rho'(\alpha + \beta + 1) p^{(2 \ell + \beta)(\beta + 1)}. \end{equation}
\end{rem}

\begin{rem} One might ask whether Theorem \ref{thm:klencase} is somehow a special case of Theorem \ref{thm:mainnasty}, or what the relationship might be. Indeed, there is a relationship between them, but it is quite technical. See \S\ref{ssec:recursion1l} for an example where this edge case is related to a case in the standard regime. \end{rem}

\textbf{Case II:} $a>0$ and $b=0$ and $0 < k < \ell$.

\begin{thm} \label{thm:bzerokmain} Fix $a > 0$ and $i \in \Om$. Let $b = 0$ and $0 < k < \ell := a+b+1$. Let
\begin{equation} \label{alphadefcaseII} a = 2\alpha + 2 \;\; \text{or} \;\; 2\alpha+1. \end{equation}
Then
\begin{equation} \label{bzerokmaineasyformula} \Xi(a,0,i,k) = (-1)^{k+1} q^{2\binom{k}{2}} z^{-\binom{\ell}{2}}  \rho'(\alpha) p^{k(3\ell - 1)} \whyme \end{equation}
where
\begin{equation} \whyme = \begin{cases}
	0 & \text{ if $i = 1$ and $a$ is odd,} \\
	p^{- 2\ell } & \text{ if $i = 1$ and $a$ is even,} \\
	p^{- 3k} & \text{ if $i = 2$ and $a$ is odd,} \\
	p^{3\ell-6k-3} & \text{ if $i = 2$ and $a$ is even,} \\
	-p^{-\ell-3k} & \text{ if $i = 3$ and $a$ is odd,} \\
	p^{-\ell-3} & \text{ if $i = 3$ and $a$ is even.} \end{cases} \end{equation}
\end{thm}

For sake of completeness, we spell out the remaining cases.

\textbf{Case III:} $a = 0$ and $b > 0$ and $k = \ell$. The case $k=0$ is included via symmetry. We define this formula using \eqref{sitauonXi}.

\begin{thm} Fix $b > 0$ and $i \in \Om$ and let $\ell = b+1$. Then
\begin{equation} \label{klenazero} \Xi(0,b,i,\ell) = \Xi(0,b,i+1,0) = (-z)^{-\ell} \overline{\Xi(b,0,-i-1,0)} = (-z)^{-\ell} \overline{\Xi(b,0,-i-2, \ell)}, \end{equation}
where $\Xi(b,0,-i-2,\ell)$ is given in Theorem \ref{thm:klencase}.
%  More explicitly, define $\beta$ by
% \begin{equation} \label{alphadefklen} b = 2\beta \quad \text{or} \quad b = 2 \beta + 1. \end{equation}
% Then $\Xi(0,b,2,\ell) = \Xi(0,b,3,0) = 0$, and for $i \ne 2$ one has
% \begin{equation} \Xi(0,b,i,\ell) = (-z)^{-\ell} (-1)^{\ell-1} q^{\binom{\beta + 1}{2}} z^{\binom{\ell}{2}} (-1)^{\alpha} \rho(\alpha) \overline{\nabla}
% 	= (z)^{-\ell} q^{\binom{\beta + 1}{2}} z^{\binom{\ell}{2}} (-1)^{\beta-1} \rho(\beta) \overline{\nabla}. \end{equation}
\end{thm}

\textbf{Case IV:} $a=0$ and $b>0$ and $0 < k < \ell$. We define this formula using \eqref{sitauonXi}.

\begin{thm} Fix $b > 0$ and $i \in \Om$ and let $\ell = b+1$. Then
\begin{equation} \label{kmainazero} \Xi(0,b,i,k) = (-z)^{-\ell} \overline{\Xi(b,0,-i-1,\ell-k)}, \end{equation}
where $\Xi(b,0,-i-1,\ell-k)$ is given in Theorem \ref{thm:bzerokmain}.
\end{thm}

%========================================================
\subsection{Proof of the formulas: they satisfy symmetry} \label{ssec:formulassymmetric}
%========================================================

We need to check that our formulas satisfy the symmetries \eqref{allsymmetries}, in order to apply Lemma \ref{lem:recursive}. 
% Indeed, \eqref{sionXi} was used to define $\Xi$ when $k = 0$, and \eqref{sitauonXi} was used to define $\Xi$ when $a = 0$. Thus these symmetries hold automatically, and we need only check  \eqref{Xi1symmetry}, and \eqref{Xi23symmetry}.

\begin{lem} When $k \in \{0,\ell\}$, the formulas from \S\ref{ssec:formulatrue} satisfy \eqref{allsymmetries}. \end{lem}

\begin{proof} As discussed in \S\ref{ssec:exploit}, \eqref{secretsymmetry} encapsulates all the symmetries between nonzero values of $\Xi$ for $k \in \{0,\ell\}$. We reprint it here.
\[ \Xi(a,b,2,0) = \Xi(a,b,1,\ell) = -z^{\ell} \Xi(a,b,1,0) = -z^{\ell} \Xi(a,b,3,\ell). \]
Formulas for $k = 0$ are derived from formulas for $k = \ell$ using \eqref{sionXi}, so two of these equalities hold by construction. We need only check that $\Xi(a,b,1,\ell) = -z^{\ell} \Xi(a,b,3,\ell)$.

When $a > 0$, consider \eqref{kleneasyformula}. All the dependency on $i$ is within the factor $\nabla$, which satisfies the desired ratio of $-z^{\ell}$ between the cases of $i=1$ and $i=3$. 

When $a = 0$, \eqref{klenazero} implies that
\[ \frac{\Xi(0,b,1,\ell)}{\Xi(0,b,3,\ell)} = \overline{\left(\frac{\Xi(b,0,3,\ell)}{\Xi(b,0,1,\ell)}\right)}.\]
Now the result holds since $-z^{\ell} = \overline{-z^{-\ell}}$.
\end{proof}

\begin{lem} When $0 < k < \ell$ and $b = 0$ or $a = 0$, the formulas from \S\ref{ssec:formulatrue} satisfy \eqref{allsymmetries}. \end{lem}

\begin{proof} We defined our formulas in the case $a=0$ by using \eqref{sitauonXi} from the case $b=0$, so this symmetry is manifestly satisfied. There are no surprising
consequences of these symmetries when $k \notin \{0,\ell\}$, see the discussion around \eqref{secretsymmetry}.

It remains to confirm (when $b = 0$) that \eqref{bzerokmaineasyformula} is compatible with \eqref{Xi1symmetry} and \eqref{Xi23symmetry}. The proof is bookkeeping. Since this is the first such bookkeeping proof in the chapter, we have written it out in full.

For either equation, both sides have a constant factor of $q^{-\binom{\alpha+1}{2}} z^{-\binom{\ell}{2}} \rho(\alpha)$.  We have
\begin{equation} \frac{\Xi(a,0,i,k)}{\Xi(a,0,i',\ell-k)} = (-1)^{\ell} q^{2\binom{k}{2} -2 \binom{\ell-k}{2}} p^{(2k - \ell)(3 \ell - 1)} \frac{\whyme(a,0,i,k)}{\whyme(a,0,i',\ell-k)}. \end{equation}
It is a straightforward computation that
\begin{equation} \binom{\ell-k}{2} - \binom{k}{2} = \frac{1}{2} (-\ell+1)(2k-\ell) \end{equation}
and hence
\begin{equation} \label{ybinomcrap} q^{2\binom{\ell-k}{2} -2 \binom{k}{2}} = p^{3(\ell-1)(2k-\ell)}. \end{equation}
Thus
\begin{equation} \frac{\Xi(a,0,i,k)}{\Xi(a,0,i',\ell-k)} = (-1)^{\ell} z^{2k-\ell} \frac{\whyme(a,0,i,k)}{\whyme(a,0,i',\ell-k)}. \end{equation}

The case of \eqref{Xi1symmetry} when $a$ is odd is trivial, as both sides are zero.

For the case of \eqref{Xi1symmetry} when $a$ is even, the factor $\whyme$ is the same for the numerator and denominator. Since $\ell$ is odd we have
\begin{equation} \frac{\Xi(a,0,1,k)}{\Xi(a,0,1,\ell-k)} = -z^{2k-\ell}, \end{equation}
as desired.

For the case of \eqref{Xi23symmetry} when $a$ is odd, we have
\begin{equation} \frac{\Xi(a,0,2,k)}{\Xi(a,0,3,\ell-k)} = (-1)^{\ell} z^{2k-\ell} \frac{p^{-3k}}{-p^{-\ell-3(\ell-k)}}. \end{equation}
Note that $\ell$ is even, but an extra sign appears in the ratio of the two values of $\whyme$. So
\begin{equation} \frac{\Xi(a,0,2,k)}{\Xi(a,0,3,\ell-k)} = -z^{2k-\ell} p^{-3k+4\ell-3k} = -z^{\ell-k}, \end{equation}
as desired.

For the case of \eqref{Xi23symmetry} when $a$ is even, we have
\begin{equation} \frac{\Xi(a,0,2,k)}{\Xi(a,0,3,\ell-k)} = (-1)^{\ell} z^{2k-\ell} \frac{p^{3 \ell - 6k - 3}}{p^{-\ell-3}}. \end{equation}
Since $\ell$ is odd we have
\begin{equation}  \frac{\Xi(a,0,2,k)}{\Xi(a,0,3,\ell-k)} = -z^{2k-\ell} p^{4 \ell - 6k} = -z^{\ell - k}, \end{equation}
as desired. \end{proof}

\begin{lem} In the standard regime, the formulas of Theorem \ref{thm:mainnasty} satisfy the symmetries \eqref{Xi1symmetry} and \eqref{Xi23symmetry}. \end{lem}

\begin{proof} We need to compute the ratio
\[ \frac{\Xi(a,b,i,k)}{\Xi(a,b,i',\ell-k)} \]
when either $i=i'=1$, or $i=2$ and $i'=3$. As the numerator and denominator have the same value of $a$ and $b$, $\alpha$ and $\beta$ and $\ell$ and $\nu$ are unchanged. No change is therefore made to $\gamma_1$ or $\lambda_r$ for $r \in \{1, 2, 3, 4\}$.

The most interesting part of the computation was already done in Corollaries \ref{cor:gammasym1} and \ref{cor:gammasym2}, where we proved that
\begin{equation} \frac{\gamma_3(a,b,i,k)}{\gamma_3(a,b,i',\ell-k)} = q^{\beta(2k-\ell)}. \end{equation}
Regardless of $i$ and $i'$ we have
\begin{equation} \frac{\kappa_1(a,b,i,k)}{\kappa_1(a,b,i',\ell-k)} = z^{k-(\ell-k)} q^{k(k-\beta-\ell) - (\ell-k)(-k-\beta)} = z^{2k-\ell} q^{-\beta(2k-\ell)}. \end{equation}
When $i = i' = 1$ we have $\kappa_2 \lambda_5 = 1$, while when $i= 2$ and $i'=3$ we have
\begin{equation} \frac{\kappa_2 \lambda_5(a,b,2,k)}{\kappa_2 \lambda_5(a,b,3,\ell-k)} = q^{2k} z^{2\ell}. \end{equation}
Altogether, we have
\begin{equation} \frac{\gamma_3\kappa_1 \kappa_2 \lambda_5(a,b,i,k)}{\gamma_3\kappa_1 \kappa_2 \lambda_5(a,b,i',\ell-k)} = \begin{cases} z^{2k-\ell} & \text{ when } i=i'=1, \\ z^{\ell-k} & \text{ when } i=2, i'=3. \end{cases} \end{equation}

Thus in both cases, the symmetries hold so long as
\begin{equation} \label{desired} \frac{\gamma_2 \mu(a,b,i,k)}{\gamma_2 \mu(a,b,i',\ell-k)} = -1. \end{equation}
When $\ell$ is odd, $\mu$ differs by a sign on top and bottom, so we will show that $\gamma_2$ is the same on top and bottom. When $\ell$ is even, $\mu$ is the same on top and bottom, so we will show that $\gamma_2$ differs by a sign.

When $a$ and $b$ are odd, $\ell$ is odd. Also, $\gamma_2 = 1$ on both top and bottom.

When $a$ and $b$ are even, $\ell = 2\nu + 1$ is odd. We have
\begin{subequations}
\begin{equation} \frac{\gamma_2(a,b,1,k)}{\gamma_2(a,b,1,\ell-k)} = \frac{[k-\nu][\nu+1-k]}{[\ell-k-\nu][\nu+1-\ell+k]} =  \frac{[k-\nu][\nu+1-k]}{[\nu+1-k][k-\nu]} = 1. \end{equation}
\begin{equation} \frac{\gamma_2(a,b,2,k)}{\gamma_2(a,b,3,\ell-k)} = \frac{[k-\nu][\ell-1-k]}{[\ell-k-1][(\ell-\nu)-(\ell-k)]} = 1. \end{equation}	

When $a$ is even and $b$ odd, $\ell = 2\nu$ is even. We have
\begin{equation} \frac{\gamma_2(a,b,i,k)}{\gamma_2(a,b,i',\ell-k)} = \frac{[k-\nu]}{[\ell-k-\nu]} = \frac{[k-\nu]}{[\nu-k]} = -1. \end{equation}

When $a$ is odd and $b$ is even, $\ell = 2\nu$ is even. We have
\begin{equation} \frac{\gamma_2(a,b,1,k)}{\gamma_2(a,b,1,\ell-k)} = \frac{-[k-\nu]}{-[\ell-k-\nu]} = \frac{[k-\nu]}{[\nu-k]} = -1. \end{equation}
\begin{equation} \frac{\gamma_2(a,b,2,k)}{\gamma_2(a,b,3,\ell-k)} = \frac{[\ell-1-k]}{[1-(\ell-k)]} = -1. \end{equation}
\end{subequations}

In all cases, \eqref{desired} holds.
\end{proof}

%========================================================
\subsection{Proof of the recursive formula: when $b=0$.}\label{ssec:recursionb0}
%========================================================

When $b=0$ and the length is sufficiently small, all cases of the recursive formula are elementary and have been checked by computer. We ignore several edge cases below.

We first prove \eqref{recursiveformula2b}, which states that
\begin{equation} \label{proveme1}
\Xi(a,0,2,k) = \Xi(a-1,0,1,k) -p^{2a} \Xi(a-1,0,3,k)
\end{equation}
when $a>1, k<\ell-1 = a$.

We split into two cases, based on the parity of $a$. We use \eqref{bzerokmaineasyformula} to define all relevant values of $\Xi$.

Suppose $a=2 \alpha + 2$ is even. Then $a-1 = 2 \alpha + 1$, so the same value of $\alpha$ applies to both sides of \eqref{proveme1}. All three terms in \eqref{proveme1} have a common factor of $(-1)^{k+1}  y^{-\binom{k}{2}} \rho'(\alpha)$. The remainder of the proposed identity says
\[
z^{-\binom{\ell}{2}} p^{k(3\ell - 1)} p^{3 \ell - 6k - 3} = 0 - p^{2(\ell-1)} z^{-\binom{\ell-1}{2}} p^{k(3(\ell-1) - 1)} (-p^{-(\ell-1) - 3k}).
\]
We leave the reader to verify this equality, using the fact that
\[ \binom{\ell}{2} - \binom{\ell-1}{2} = \ell-1. \]

Suppose that $a = 2 \alpha + 1$ is odd. Then $a-1 = 2(\alpha-1) + 2$,  so the right-hand side of \eqref{proveme1} uses $\alpha-1$ instead of $\alpha$. All three terms have a common factor of $(-1)^{k+1} y^{-\binom{k}{2}}\rho'(\alpha-1)$. Dividing by this and by $z^{-\binom{-\ell-1}{2}}p^{3k(\ell-1)-1}$, the remainder of the proposed identity says
\[
z^{-\ell+1} (1 - q^{-2\alpha}) p^{3k} p^{-3k} = p^{-2(\ell-1)} - p^{2(\ell-1)} p^{-(\ell-1)-3}.
\]
Again, we leave this simple verification to the reader.

Now we prove \eqref{recursiveformula2bspecial}, which states that
\begin{equation} \label{proveme2}
\Xi(a,0,2,\ell-1) = \Xi(a-1,0,1,\ell-1)
\end{equation}
when $a > 1$ and $\ell = a+1$.  The left-hand side is defined using \eqref{bzerokmaineasyformula}, while the right-hand side is defined using \eqref{eq:bzerokleneasy} with a different formula for $\alpha$.

When $a$ is even, $a = 2 \alpha + 2$ and $a-1 = 2 \alpha + 1$, so the same value of $\alpha$ is used for both terms. The signs match: $(-1)^{\ell}$ versus $(-1)^{\ell-2}$. There is a common factor of $\rho'(\alpha)$.  Dividing by this common factor and by $z^{-\binom{\ell-1}{2}}$, the remainder of the proposed equality is
\begin{equation} q^{2\binom{\ell-1}{2}} z^{-\ell+1} p^{(\ell-1)(3\ell - 1)} p^{3\ell-6(\ell-1)-3} = 1, \end{equation}
which is easily verified.

When $a$ is odd, $a = 2\alpha + 1$ and $a-1 = 2 \alpha$, so again the same value of $\alpha$ is used for both terms. Ignoring the same common factor and the same signs, what remains is
\[q^{2\binom{\ell-1}{2}} z^{-\ell+1} p^{(\ell-1)(3\ell - 1)} p^{- 3(\ell-1)} = 1. \]
This is the same computation after observing that $p^{3\ell-6(\ell-1)-3} = p^{- 3(\ell-1)}$, an equality of two different entries in the table for $\whyme$ when $k = \ell-1$.

Now we prove \eqref{recursiveformula1b}, which states that
\begin{equation} \label{proveme3}
\Xi(a,0,1,k) = \sum_{c = \ell - k}^{k-1} z^{k-1-c} \Xi(a-1,0,3,c)
\end{equation}
when $a > 0$ and $k > \ell-k$ (the case $k = \ell-k$ is zero by symmetry). This is a very different computation. We begin with the case $k < \ell$, so that all terms are defined using \eqref{bzerokmaineasyformula}.

When $a$ is odd, the left side is zero, and we will argue that the terms in the right side cancel in pairs. Dividing by the factors in \eqref{bzerokmaineasyformula} which are constant as $c$ varies, the right side of \eqref{proveme3} simplifies to
\begin{equation} \sum_{c=\ell-k}^{k-1} z^{-c} (-1)^c  q^{2\binom{c}{2}} p^{c(3 (\ell-1) - 1)} = \sum_{c=\ell-k}^{k-1} (-1)^c p^{3c(\ell-1-c)}. \end{equation}
The exponent of $p$ is a quadratic polynomial in $c$, which is symmetric around the half-integer $\frac{\ell-1}{2}$. This is also the average of the two endpoints $k-1$ and $\ell - k$ of the sum. There are an even number of summands, which cancel in pairs due to the sign $(-1)^c$.

Similarly, in the case where $a$ is even, the right hand side is a constant factor times the sum
\begin{equation} \sum_{c = \ell-k}^{k-1} (-1)^c z^{-c} q^{2\binom{c}{2}} p^{c(3 (\ell-1) - 1)}p^{-3c} = \sum_{c = \ell-k}^{k-1} (-1)^c p^{3c(\ell-2-c)}, \end{equation}
with the extra factor of $p^{-3c}$ on the left side coming from $\whyme$. 	The exponent of $p$ is a quadratic polynomial in $c$, which is symmetric around the half-integer $\frac{\ell-2}{2}$. All the terms for $c$ between $\ell-k$ and $k-2$ will cancel in pairs, so the only surviving term has $c = k-1$, and we wish to prove that
\[ \Xi(a,0,1,k) = \Xi(a-1,0,3,k-1). \]
Note that $a = 2 \alpha + 2$, so $a-1 = 2 \alpha + 1$, and the same value of $\alpha$ is used on both sides. Dividing by the typical common factors, our remaining equality is
\[ (-1)^{k+1} q^{2(k-1)} z^{-(\ell-1)} p^{k(3\ell - 1)} p^{-2\ell} = (-1)^{k} p^{(k-1)(3(\ell-1)-1)} (-p^{-(\ell-1) - 3(k-1)}). \]
The signs clearly match. Making the exponents match is an exercise.

Finally, we need to check \eqref{proveme3} when $k = \ell$. Now the left side, and the $c = 0$ and $c = \ell-1$ terms on the right side, are governed by \eqref{eq:bzerokleneasy} instead of \eqref{bzerokmaineasyformula}. Also, the $c=0$ term on the right side is zero, see \eqref{benstable}.

Suppose that $a$ is odd. By the same computation as above, all the terms on the right side except for ($c=0$ and) $c = \ell-1$ will cancel in pairs. We need to confirm that
\begin{equation} \Xi(a,0,1,\ell) = \Xi(a-1,0,3,\ell-1). \end{equation}
Note that $\alpha$ is being defined by \eqref{alphabetadefklen}, and $a = 2 \alpha+1$ and $a-1 = 2 \alpha$, so the same value of $\alpha$ is being used on both sides of the equation. Dividing by the common factor $\rho'(\alpha)$ and by $z^{-\binom{\ell-1}{2}}$, the remainder of the proposed equality is
\begin{equation} (-1)^{\ell-1} z^{-\ell+1} = (-1)^{\ell-2} (-z^{-(\ell-1)}), \end{equation}
where $\nabla$ was nontrivial on the right side. This equality obviously holds.

Suppose that $a$ is even. By the same computation as above, all the terms on the right side except for ($c=0$ and) $c = \ell-2$ and $c = \ell-1$ will cancel in pairs. We need to confirm that
\begin{equation} \Xi(a,0,1,\ell) = z \Xi(a-1,0,3,\ell-2) + \Xi(a-1,0,3,\ell-1). \end{equation}
This time, $a = 2 \alpha$ and $a-1 = 2(\alpha-1) + 1$, so a different value of $\alpha$ is being used on the right hand side. Dividing both sides by $\rho'(\alpha-1) z^{-\binom{\ell-1}{2}}$, the remainder of the proposed equality is
\begin{align} \nonumber (-1)^{\ell-1} & z^{-\ell+1} (1 - q^{-2\alpha}) = \\ & (-1)^{\ell-1} z q^{2\binom{\ell-2}{2}}p^{(\ell-2)(3(\ell-1) - 1)}(-p^{-(\ell-1)-3(\ell-2)}) + (-1)^{\ell-2}(-z^{-(\ell-1)}). \end{align}
Again, the verification is an exercise.

%\BE{I have done all exercises in this chapter!!}

%========================================================
\subsection{Proof of the recursive formula: when $b=1$.}\label{ssec:recursionb1}
%========================================================

The formulas \eqref{recursiveformula1} and \eqref{recursiveformula2} when $b=1$ also involve the edge case formulas as $b-1 = 0$.

Consider the case of \eqref{recursiveformula1} when $b = 1$ and $k = \ell$. It states that
\[ \Xi(a,1,1,\ell) = \sum_{c = 0}^{\ell-1} z^{\ell-1-c} \Xi(a,0,2,c). \]
However, by Theorem \ref{thm:klencase}, $\Xi(a,1,1,\ell) = 0$ since $b$ is odd. Letting $\ell = a+1$ to easier study the right side, we need to prove that
\begin{equation} \label{proveme4} \sum_{c = 0}^{\ell} z^{\ell - c} \Xi(a,0,2,c) = 0. \end{equation}
The $c = \ell$ term is zero. Meanwhile, $\Xi(a,0,2,0) = \Xi(a,0,1,\ell)$. So \eqref{proveme4} is equivalent to
\begin{equation} \label{proveme5} \Xi(a,0,1,\ell) = -z^{-\ell} \sum_{c = 1}^{\ell-1} z^{\ell-c} \Xi(a,0,2,c) = - \sum_{c = 1}^{\ell-1} z^{-c} \Xi(a,0,2,c). \end{equation}	
Using \eqref{bzerokmaineasyformula}, the right side of \eqref{proveme5}, up to a constant factor independent of $c$, is
\begin{equation} \sum_{c=1}^{\ell-1} (-1)^{c} z^{-c} q^{2\binom{c}{2}} p^{c(3 \ell - 1)} \begin{cases} p^{-3c} & \text{ $a$ is odd} \\ p^{-6c} & \text{ $a$ is even} \end{cases}. \end{equation}

When $a$ is odd, the summands have the form $(-1)^c p^{3c(c + \ell - 1)}$. The exponent of $p$ is symmetric around the half-integer $\frac{\ell-1}{2}$. Thus all terms from $c = 1$ to $c = \ell-2$ cancel out, and only the $c = \ell- 1$ term survives. It remains to confirm
\begin{equation} \Xi(a,0,1,\ell) = - z^{-\ell+1} \Xi(a,0,2,\ell-1), \quad \text{ when $a$ is odd.} \end{equation}
Using \eqref{eq:bzerokleneasy} and \eqref{bzerokmaineasyformula}, and noting that both sides use the same value of $\alpha$, this amounts to showing that
\begin{equation} (-1)^{\ell-1} = -z^{-\ell+1} (-1)^{\ell} q^{2\binom{\ell-1}{2}} p^{(\ell-1)(3\ell - 1)} p^{-3(\ell-1)}. \end{equation}
This is an elementary exericse.

When $a$ is even, the summands have the form $(-1)^c p^{3c(c+\ell - 2)}$. The exponent of $p$ is symmetric around the half-integer $\frac{\ell-2}{2}$. Thus all terms from $c=1$ to $c = \ell - 3$ cancel out, and only the $c = \ell-2$ and $c = \ell-1$ terms survive. It remains to confirm
\begin{equation} \Xi(a,0,1,\ell) = -z^{-\ell+1} \Xi(a,0,2,\ell-1) - z^{-\ell+2} \Xi(a,0,2,\ell-2), \quad \text{ when $a$ is even}. \end{equation}
Using \eqref{eq:bzerokleneasy} and \eqref{bzerokmaineasyformula}, this time the values of $\alpha$ differ. Letting $a = 2 \alpha$, and dividing by $\rho'(\alpha-1)$, we need to show that
\begin{align} \nonumber (-1)^{\ell-1} (1 - q^{-2\alpha}) =  
-z^{-\ell+1} & (-1)^{\ell} q^{2\binom{\ell-1}{2}} p^{(\ell-1)(3\ell - 1)} p^{3\ell-6(\ell-1)-3}  \\
& - z^{-\ell+2} (-1)^{\ell-1} q^{2\binom{\ell-2}{2}} p^{(\ell-2)(3\ell - 1)} p^{3\ell-6(\ell-2)-3}. \end{align}
This is another elementary exercise.

Now consider the case of \eqref{recursiveformula1} when $b=1$ and $\ell/2 \le k < \ell = a+2$. It states that
\[ \Xi(a,1,1,k) = \sum_{c = \ell-k}^{k-1} z^{k-1-c} \Xi(a,0,2,c). \]
This time the left hand side is governed by Theorem \ref{thm:mainnasty}, with $\beta = 0$, while the right hand side is governed by Theorem \ref{thm:bzerokmain}. Both sides have a factor of $\rho'(\alpha)$, which we divide by henceforth.

Suppose $a = 2 \alpha + 1$ is odd. Then 
\begin{align} LHS = (-1)^k z^k q^{k(k-\ell)} z^{-\binom{\ell+1}{2}} p^{-\ell}z^{\ell+1}p^{-3} =
	(-1)^k q^{k(k-\ell)} p^{- \ell^2 + 2k - 1}. \end{align}
\begin{align} \nonumber RHS = \sum_{c = \ell-k}^{k-1} z^{k-1-c} & (-1)^{c+1} q^{2\binom{c}{2}} z^{-\binom{\ell-1}{2}} p^{c(3\ell-4)}p^{-3c}
	\\ & = \sum_{c = \ell-k}^{k-1} (-1)^{c+1} p^{2k-4-6c - \ell^2 + 3 \ell + 3 \ell c - 3 c^2} . \end{align}
Dividing both sides by $(-1)^k p^{2k - 1 - \ell^2}$, we need to prove that
\begin{equation} q^{k(k-\ell)} = \sum_{c = \ell-k}^{k-1} (-1)^{k+c+1} q^{1 + 2c - \ell - \ell c + c^2}. \end{equation}
When $c = k-1$, the summand in the RHS agrees with the LHS. The remaining terms cancel using the quadratic argument above; the exponent is symmetric around $\frac{\ell-2}{2}$, which is halfway between $\ell-k$ and $k-2$.

Suppose $a = 2 \alpha + 2$ is even. Then
\begin{align} \nonumber LHS = (-1)^k & q^{-\nu}(q-q^{-1})[k-\nu] z^k q^{k(k-\ell)} z^{-\binom{\ell+1}{2}}p^{-\ell}z^{-3} z^{\ell+1}
	\\ & = (-1)^k q^{k - 2 \nu + k^2 - k\ell} p^{2k-\ell^2 -4} + (-1)^{k+1} q^{-k+k^2 - k \ell} p^{2k-\ell^2 -4} . \end{align}
\begin{align} \nonumber RHS = \sum_{c = \ell-k}^{k-1} & z^{k-1-c} (-1)^{c+1} q^{2\binom{c}{2}} z^{-\binom{\ell-1}{2}}p^{c(3\ell-4)}p^{3 \ell -6c-6}
	\\ & = \sum_{c = \ell-k}^{k-1} (-1)^{c+1} p^{2k -3c^2  - \ell^2 + 6 \ell  + 3 c \ell - 9c - 10}. \end{align}
Dividing both sides by $(-1)^k p^{2k - \ell^2  - 4}$ and noting that $2 \nu = \ell$, we need to prove that
\begin{equation} q^{k - \ell + k^2 - k\ell} - q^{-k+k^2 - k \ell} = \sum_{c = \ell-k}^{k-1} (-1)^{k+c+1} q^{ c^2  -2 \ell  - c \ell +3c +2}. \end{equation}
The $c = k-1$ and $c = k-2$ terms on the RHS match the two terms on the LHS. The remaining terms cancel using the quadratic argument above; the exponent is symmetric around $\frac{\ell-3}{2}$, which is halfway between $\ell-k$ and $k-3$.

Now we treat the $b=1$ case of \eqref{recursiveformula2} when $0 < k < \ell-1$, which states that
\[ \Xi(a,1,2,k) = \Xi(a,0,3,k) - z^{2\ell - 3k - 2} \Xi(a,0,1,k). \]
Again, the left hand side is governed by Theorem \ref{thm:mainnasty}, with $\beta = 0$, while the right hand side is governed by Theorem \ref{thm:bzerokmain}. Both sides have a factor of $\rho'(\alpha)$, which we divide by henceforth.

Suppose $a = 2 \alpha + 1$ is odd. Then 
\begin{align} LHS = (-1)^k z^k q^{k(k-\ell)} q^{2k} z^{-\binom{\ell+1}{2}} p^{-\ell}z^{\ell+1}p^{-3} z^{\ell} =
	(-1)^k p^{-3k^2 + 3 k \ell - \ell^2 + 2 \ell  - 4k - 1}. \end{align}
\begin{align} RHS = (-1)^{k+1} q^{2\binom{k}{2}}z^{-\binom{\ell-1}{2}}p^{k(3\ell-4)}(-p^{-\ell+1-3k}) + 0 = (-1)^k p^{-3 k^2 - \ell^2 + 2 \ell + 3 k \ell - 4k -1}. \end{align}
The two sides clearly agree.

Suppose $a = 2 \alpha + 2$ is even. Then
\begin{align} \nonumber LHS = (-1)^k & q^{-\nu}(q-q^{-1})[k-\nu] z^k q^{k(k-\ell)} q^{2k} z^{-\binom{\ell+1}{2}}p^{-\ell}z^{-3} z^{\ell+1}z^{\ell}
	\\ &= (-1)^k q^{3k - 2 \nu + k^2 - k\ell} p^{2k-\ell^2 -4 + 2 \ell} + (-1)^{k+1} q^{k+k^2 - k \ell} p^{2k-\ell^2 -4 + 2 \ell} . \end{align}
\begin{align} \nonumber RHS = (-1)^{k+1} & q^{2\binom{k}{2}}z^{-\binom{\ell-1}{2}}p^{k(3\ell-4)} \left( p^{-\ell+1-3}) - z^{2\ell - 3k - 2} p^{-2(\ell-1)} \right)
	\\ &= (-1)^k p^{-3k^2  - \ell^2 + 5 \ell + 3 \ell k - 7k - 4} + (-1)^{k+1} p^{-3k^2 -k - \ell^2 + 2 \ell - 4 + 3 \ell k}. \end{align}
Making these match is another exercise.

%\BE{I did all exercises in this chapter.}

%========================================================
\subsection{Proof of \eqref{recursiveformula2special}}\label{ssec:recursionspecial}
%========================================================

In this section we check \eqref{recursiveformula2special} which mixes the standard regime with an edge case. It states that
\[\Xi(a,b,2,\ell-1) = \Xi(a,b-1,3,\ell-1) \quad \text{ if } a,b > 0. \]
Note that $\Xi(a,b-1,3,\ell-1)$ is defined using \eqref{klen3}, regardless of whether or not $b-1=0$, so that $b=1$ is not a special edge case of the formula. One should also recall that \eqref{klen3} uses the definition of $\alpha$ from \eqref{alphabetadefklen}. Meanwhile, $\Xi(a,b,2,\ell-1)$ is defined using \eqref{eq:mainnasty}.

Recall from Theorem \ref{thm:klencase} that $\Xi(a,b-1,3,\ell-1) = 0$ if $b-1$ is odd. Meanwhile, if $b$ is even and $i=2$, then $[\ell-1-k]$ divides $\gamma_2$, so when $k = \ell-1$ we have $\gamma_2 = 0$. Thus both sides of \eqref{recursiveformula2special} are zero when $b$ is even.

Henceforth $b = 2 \beta + 1$ is odd. Note that $b-1 = 2(\beta-1) + 2$, so the two sides of \eqref{recursiveformula2special} use different values of $\beta$. Let us simplify the formula from \eqref{eq:mainnasty} under the assumption $k = \ell-1$ and $i=2$.
\begin{equation} \kappa_1 \kappa_2 = p^{(\ell-1)(3 \beta -1)}. \end{equation}

We begin with the case where $a = 2 \alpha + 1$ is odd, so $\ell = 2 \nu - 1$. Both sides of \eqref{recursiveformula2special} use the same value of $\alpha$. Indeed, \eqref{klen3} (with length $\ell-1$ and values $\alpha$ and $\beta - 1$) becomes
\begin{equation} \label{klen3modifiedaodd} \Xi(a,b-1,3,\ell-1) = (-1)^{\beta+\ell-1} z^{-\binom{\ell}{2}} \rho'(\alpha + \beta) p^{(2 \ell + \beta - 3)(\beta)}. \end{equation}

The sign clearly matches with $\thesign$. The factor of $z^{-\binom{\ell}{2}}$ in \eqref{klen3modifiedaodd} matches $\lambda_5$ times the factor of $z^{-\binom{\ell+1}{2}}$ in $\lambda_1$.
Let us note that
\begin{equation} \frac{\rho'(\alpha+\beta)}{\rho'(\alpha)\rho'(\beta)} = q^{\binom{\alpha+1}{2} + \binom{\beta+1}{2} - \binom{\alpha+\beta+1}{2}} {\alpha+\beta \brack \beta} = q^{-\alpha \beta} {\alpha+\beta \brack \beta}. \end{equation}
Thus
\begin{equation} \frac{\Xi(a,b,2,\ell-1)}{\Xi(a,b-1,3,\ell-1)} = \frac{\gamma_3}{{\alpha+\beta \brack \beta}} q^{\alpha \beta} \frac{p^{(\ell-1)(3 \beta -1)}}{p^{(2 \ell + \beta - 3)(\beta)}} z^{\binom{\beta}{2}}  p^{-(\beta+1)(\ell+3\beta)} z^{\beta} z^{\ell+1} p^{-(\beta+3)}, \end{equation}
which simplifies to
\begin{equation} \frac{\gamma_3}{{\alpha+\beta \brack \beta}} q^{\alpha \beta} p^{-3 \beta^2 - 3 \beta} = \frac{\gamma_3}{{\alpha+\beta \brack \beta}} q^{\beta(\alpha + \beta + 1)}. \end{equation}
	
It remains to prove that
\begin{equation} \label{specialmagic1} q^{\beta(\nu-1)} \magic(\nu,2 \nu - 2,\beta,-1) = {\alpha + \beta \brack \beta}. \end{equation}
This was proven in \eqref{magicchuvan}.
	
Now consider the case when $a = 2 \alpha + 2 = 2(\alpha+1)$ is even, and hence $\ell = 2 \nu$. Now the two sides of \eqref{recursiveformula2special} use different values of $\alpha$.
Indeed, \eqref{klen3} (with length $\ell-1$ and values $\alpha+1$ and $\beta - 1$) becomes
\begin{equation} \label{klen3modifiedaeven} \Xi(a,b-1,3,\ell-1) = (-1)^{\beta+\ell-1} z^{-\binom{\ell}{2}} \rho'(\alpha + \beta+1) p^{(2 \ell + \beta - 3)(\beta)}. \end{equation}

We leave the reader to perform a similar analysis to the above. The slight complication is that $\gamma_2$ is nontrivial, but $\gamma_2 \rho'(\alpha + \beta)$ is equal to $\rho'(\alpha+\beta+1)$ up to a power of $q$. Eventually one computes that
\begin{equation} \frac{\Xi(a,b,2,\ell-1)}{\Xi(a,b-1,3,\ell-1)} = q^{\beta\nu} \frac{\gamma_3}{{\alpha + \beta \brack \beta}}. \end{equation}
It remains to prove that
\begin{equation} \label{specialmagic2} q^{\beta\nu} \magic(\nu,2 \nu -1,\beta,0) =  {\alpha + \beta \brack \beta}. \end{equation}
This was also proven in \eqref{magicchuvan}.

%\BE{full Check is in nastytwo.tex, but in previous notation though}

%========================================================
\subsection{Proof of \eqref{recursiveformula1} in the standard regime}\label{ssec:recursion1}
%========================================================

We now prove \eqref{recursiveformula1} when $a > 0$, $b>  1$, and $\ell/2 \le k < \ell$. It states  
\begin{equation} \label{thebigguy} \Xi(a,b,1,k) = \sum_{c = \ell-k}^{k-1} z^{k-1-c} \Xi(a,b-1,2,c), \end{equation}
and all terms $\Xi$ are in the standard regime. (We treated the $b=1$ case earlier in \S\ref{ssec:recursionb1}.) 

We begin with the case when $b = 2 \beta + 2$ is even, so that $b-1 = 2 \beta + 1$ and both sides use the same value of $\beta$ (and $\alpha$ and $\nu$). Thus both sides have the
same value for $\gamma_1$ and $\lambda_2$. Moreover, we have
\begin{equation} \frac{\lambda_1 \lambda_3 \lambda_4 \lambda_5(a,b,1)}{\lambda_1 \lambda_3 \lambda_4 \lambda_5(a,b-1,2)} = z^{-\ell}p^{-(\beta+1)}z^{-\ell}p^{\beta+3}z^{-(\ell-1)} = z^{-3\ell+2}. \end{equation}

We begin with the case where $b$ is even and $a$ is odd, so that $\ell = 2 \nu$. Dividing \eqref{thebigguy} by $\gamma_1$ and all $\lambda_r(a,b-1,2)$ and by $(-1)^{\beta}$, what remains of the left side is
\begin{equation} LHS = (-1)^{k} q^{-\nu} (q-q^{-1}) [\nu-k] \magic(\nu,k,\beta,0) z^{-3\ell+2} z^k q^{k(k-\beta-\ell)}. \end{equation}
What remains of the right side is
\begin{equation} RHS = \sum_{c = \ell-k}^{k-1} z^{k-1-c} (-1)^c \magic(\nu,c,\beta,-1) z^c q^{c(c-\beta-\ell+3)}. \end{equation}
Multiplying both sides by $z^{1-k}$ and manipulating, we get
\begin{equation} newLHS = (-1)^k (q^{-2k} - q^{-2\nu}) \magic(\nu,k,\beta,0) q^{2 \ell-2} q^{k(k-\beta-\ell+1)}, \end{equation}
\begin{equation} newRHS = \sum_{c = \ell-k}^{k-1} (-1)^c \magic(\nu,c,\beta,-1)q^{c(c-\beta-\ell+3)}. \end{equation}

The equality of $newLHS$ and $newRHS$ is an immediate consequence of Theorem~\ref{thm:telescoping sum}.

Now we consider the case where $b$ is even and $a$ is even, so that $\ell = 2 \nu + 1$. Dividing \eqref{thebigguy} by $\gamma_1$ and all $\lambda_r(a,b-1,2)$ and by $(-1)^{\beta}$, what remains of the left side is
\begin{equation} LHS = (-1)^{k} q^{-\ell} (q-q^{-1})^2 [k-\nu][\nu+1-k] \magic(\nu,k,\beta,+1) z^{-3\ell+2} z^k q^{k(k-\beta-\ell)}. \end{equation}
What remains of the right side is
\begin{equation} RHS = \sum_{c = \ell-k}^{k-1} z^{k-1-c} (-1)^c q^{-\nu} (q-q^{-1}) [c - \nu] \magic(\nu,c,\beta,0) z^c q^{c(c-\beta-\ell+3)}. \end{equation}
Multiplying both sides by $z^{1-k}$ and manipulating, we get
\begin{equation} newLHS = (-1)^{k} (q^{2k} - q^{2\nu}) (q^{2\nu} - q^{2k-2}) \magic(\nu,k,\beta,+1) q^{k(k-\beta-\ell-2)}. \end{equation}
\begin{equation} newRHS = \sum_{c = \ell-k}^{k-1} (-1)^c (q^{-2\nu}-q^{-2c}) \magic(\nu,c,\beta,0) q^{c(c-\beta-\ell+4)}. \end{equation}
The equality of these two expressions is, by some cosmic coincidence, precisely what we proved in 
Theorem~\ref{telescope a and b even}.

Now we consider the case where $b = 2 \beta+1$ is odd. Now $b-1 = 2 (\beta-1) + 2$, so the RHS uses $\beta-1$ and $\nu-1$ in its formulas. This time we have
\begin{equation} \frac{\gamma_1(a,b)}{\gamma_1(a,b-1)} = 1 - q^{-2 \beta}, \end{equation}
and
\begin{align} \frac{\lambda_1 \lambda_3 \lambda_4 \lambda_5(a,b,1)}{\lambda_1 \lambda_3 \lambda_4 \lambda_5(a,b-1,2)} =
%	z^{\beta-1} z^{-\ell} p^{-(\beta+1)(\ell+3\beta)}p^{(\beta)(\ell+3 \beta - 4)} z^{\ell+1}p^{something that was wrong, lambda4 misstep} z^{-\ell+1} \\ & =
	p^{\phi-1} q^{\ell + 2 \beta}. \end{align}

Consider the case where $b$ is odd and $a$ is odd, so that $\ell = 2 \nu - 1$. Dividing \eqref{thebigguy} by $\gamma_1(a,b-1)$ and all $\lambda_r(a,b-1,2)$ and by $(-1)^{\beta-1}$, what remains of the left side is
\begin{equation} LHS = (-1)^{k+1} (1-q^{-2 \beta}) \magic(\nu,k,\beta,-1) z p^{-1} q^{\ell + 2 \beta} z^k q^{k(k-\beta-\ell)}. \end{equation}
What remains of the right side is
\begin{align} & RHS = \\\nonumber \sum_{c = \ell-k}^{k-1} z^{k-1-c} (-1)^c q^{-(\ell-2)} & (q-q^{-1}) [\ell-2-c] \magic(\nu-1,c,\beta-1,-1) z^c q^{c(c-\beta-\ell+4)}. \end{align}
Multiplying both sides by $z^{1-k}$ and manipulating we get
\begin{equation} newLHS = (-1)^{k} q^{\ell -1} (1 - q^{2 \beta}) \magic(\nu,k,\beta,-1) q^{k(k-\beta-\ell)}, \end{equation}
\begin{equation} newRHS = \sum_{c = \ell-k}^{k-1} (-1)^c (1-q^{2c + 6 - 4 \nu}) \magic(\nu-1,c,\beta-1,-1) q^{c(c-\beta-\ell+3)}. \end{equation}
The equality $newLHS = newRHS$ is Theorem \ref{telescope a and b odd}.

Finally, consider the case when $b$ is odd and $a$ is even, so that $\ell = 2 \nu$. Dividing \eqref{thebigguy} by $\gamma_1(a,b-1)$ and all $\lambda_r(a,b-1,2)$ and by $(-1)^{\beta-1}$, then multiplying by $z^{1-k}q^{\ell+\nu-3}$ and manipulating as before, we get
\begin{equation} LHS = (-1)^{k} q^{2 \ell - 3} (1-q^{2 \beta}) (q^{k-\nu} - q^{\nu-k}) \magic(\nu,k,\beta,0) q^{k(k-\beta-\ell)}, \end{equation}
\begin{align} & RHS = \\ \nonumber \sum_{c = \ell-k}^{k-1} (-1)^c (q^{c+1-\nu}-q^{\nu-c-1}) & (q^{\ell-2-c} - q^{c+2-\ell}) \magic(\nu-1,c,\beta-1,0) q^{c(c-\beta-\ell+4)}. \end{align}
The equality $LHS = RHS$ is Theorem \ref{telescope b odd a even}.	

%  what remains of the left side is
% \begin{equation} LHS = (-1)^{k+1} (1-q^{-2 \beta}) q^{-\nu} (q^{k-\nu} - q^{\nu-k}) \magic(\nu,k,\beta,0) z^{-1} q^{\ell + 2 \beta} z^k q^{k(k-\beta-\ell)}. \end{equation}
% What remains of the right side is
% \begin{equation} RHS = \sum_{c = \ell-k}^{k-1} z^{k-1-c} (-1)^c q^{-(\ell + \nu -3)} (q-q^{-1})^2 [c+1-\nu][\ell-2-c] \magic(\nu-1,c,\beta-1,0) z^c q^{c(c-\beta-\ell+4)}. \end{equation}
% Multiplying both sides by $z^{1-k}q^{\ell+\nu-3}$ and manipulating we get
% \begin{equation} newLHS = (-1)^{k} q^{2 \ell - 3} (1-q^{2 \beta}) (q^{k-\nu} - q^{\nu-k}) \magic(\nu,k,\beta,0) q^{k(k-\beta-\ell)}, \end{equation}
% \begin{equation} newRHS = \sum_{c = \ell-k}^{k-1} (-1)^c (q^{c+1-\nu}-q^{\nu-c-1}) (q^{\ell-2-c} - q^{c+2-\ell}) \magic(\nu-1,c,\beta-1,0) q^{c(c-\beta-\ell+4)}. \end{equation}
% The equality $newLHS = newRHS$ is Theorem \ref{telescope b odd a even}.

%========================================================
\subsection{Proof of \eqref{recursiveformula1} when $k = \ell$}\label{ssec:recursion1l}
%========================================================

We now prove \eqref{recursiveformula1} when $a,b > 0$ and $k = \ell$, another case which mixes the standard regime with an edge case. It states
\[ \Xi(a,b,1,\ell) = \sum_{c = 0}^{\ell-1} z^{\ell-1-c} \Xi(a,b-1,2,c). \]
We already treated the case $b=1$ in \S\ref{ssec:recursionb0}, so we can assume $b > 1$. The $c = \ell-1$ summand on the right side is zero. Meanwhile, the $c=0$ summand is defined using \eqref{sionXi}, via $\Xi(a,b-1,2,0) = \Xi(a,b-1,1,\ell-1)$. So we need to show
\begin{equation} \label{proveme10} \Xi(a,b,1,\ell) - z^{\ell-1} \Xi(a,b-1,1,\ell-1) = \sum_{c=1}^{\ell-2} z^{\ell-1-c} \Xi(a,b-1,2,c) = z \Xi(a,b,1,\ell-1). \end{equation}
The last equality arose from \eqref{recursiveformula1} when $k = \ell-1$, as proven in the last section. By Theorem \ref{thm:klencase}, $\Xi(a,b,1,\ell) = 0$ when $b$ is odd, and so $\Xi(a,b-1,1,\ell-1)=0$ if $b$ is even. Thus, regardless of the parity of $b$, only one of the two terms on the left-hand side of \eqref{proveme10} is nonzero.

Let us examine $\Xi(a,b,1,\ell-1)$. Define $\alpha$ and $\beta$ using \eqref{alphabetadef}. By examination of $\gamma_3$ when $i=1$, we see that
\begin{equation} \gamma_3(a,b,1,k) = \magic(\nu,k,\beta,\epsilon), \qquad \ell = 2 \nu + \epsilon. \end{equation}
In particular, when $k = \ell-1$ the conditions of Lemma \ref{lem:magicchuvan} are satisfied, and it implies that
\begin{equation} \gamma_3(a,b,1,\ell-1) = q^{\beta(\nu-\ell)} {\nu-2 \brack \beta}. \end{equation}
Since $\nu - 2 = \alpha + \beta$, and $\gamma_1$ is equal to $[\alpha]![\beta]!$ up to a power of $q$ and a power of $(q-q^{-1})$, we have
\begin{equation} \label{gamma1gamma3inproof} \gamma_1 \gamma_3(a,b,1,\ell-1) = (q-q^{-1})^{\alpha + \beta} q^{-\binom{\alpha+1}{2}} q^{-\binom{\beta+1}{2}} q^{\beta(\nu-\ell)} [\alpha+\beta]! = q^{\alpha \beta} q^{\beta(\nu-\ell)}  \rho'(\alpha+\beta). \end{equation}
This last equality followed since
\[ \binom{\alpha+\beta+1}{2} - \binom{\alpha+1}{2} = \binom{\beta+1}{2} = \alpha \beta. \]

Note that the left side of \eqref{proveme10} is defined using \eqref{kleneasyformula}, and is a sign and a power of $p$ times $\rho'(d)$ for some $d$. More precisely, the formula
\eqref{kleneasyformula} states that $d = \alpha+\beta+1$, but this is for a value of $\alpha$ and $\beta$ which need not agree with the values just used above. This discrepancy is
exactly accounted for by $\gamma_2(a,b,1,\ell-1)$, as we now justify case by case.

When $a$ and $b$ are odd, both sides of \eqref{proveme10} use the same value of $\alpha$. In the left side of \eqref{proveme10} only $\Xi(a,b-1,1,\ell-1)$ is nonzero, and it uses
$\beta-1$ instead of $\beta$. So $d = \alpha + \beta$, which matches the factor of $\rho'(\alpha+\beta)$ in \eqref{gamma1gamma3inproof}. Meanwhile, $\gamma_2(a,b,1,\ell-1) = 1$.

% Powers of $q$ and signs
% On RHS, we have
% \begin{equation} (-1)^{\ell-1+\beta} z q^{\beta(\nu - \ell + \alpha)} z^{\ell-1} q^{(\ell-1)(-\beta-1)} z^{\binom{\beta}{2}} z^{-\binom{\ell+1}{2}} p^{-(\beta+1)(\ell+3\beta)} z^{\beta} z^{\ell+1}p^{-(\beta+3)}\end{equation}
% On LHS we have
% \begin{equation} -z^{\ell-1} (-1)^{\beta+\ell-2} z^{-\binom{\ell-1}{2}} p^{(2 \ell + \beta - 3)(\beta)} \end{equation}
% signs match, now ignoring signs. simplifying
% \begin{equation} RHS/LHS =  q^{\beta(\nu - \ell + \alpha)}  q^{(\ell-1)(-\beta-1)} q^{\beta^2  + \ell  -1 + \ell \beta}. \end{equation}
% it works, I did exercise.   

When $a$ is even and $b$ is odd, \eqref{alphabetadefklen} implies that $\Xi(a,b-1,1,\ell-1)$ uses $\alpha+1$ instead of $\alpha$. It also uses $\beta-1$ instead of $\beta$. So $d =
\alpha+\beta+1$. Meanwhile, $\gamma_2(a,b,1,\ell-1)$ has a factor of $(q-q^{-1}) [\ell-1-\nu]$, and $\ell = 2\nu$, so $\ell-1-\nu = \alpha+\beta+1$. Multiplying this by the factor of
$\rho'(\alpha+\beta)$ in \eqref{gamma1gamma3inproof}, we get $\rho'(\alpha+\beta+1)$ up to powers of $q$, as desired.

When $a$ is odd and $b$ is even, both sides of \eqref{proveme10} use the same value of $\alpha$. In the left side of \eqref{proveme10} only $\Xi(a,b,1,\ell)$ is nonzero, and it uses the
same value of $\beta$. So $d = \alpha + \beta + 1$. Once again, $\gamma_2(a,b,1,\ell-1)$ has a factor of $-(q-q^{-1})[\ell-1-\nu]$, and $\ell-1-\nu = \alpha+\beta+1$.

Finally, when $a$ and $b$ are both even, we have $d = \alpha+\beta+2$. Meanwhile, $\gamma_2(a,b,1,\ell-1)$ has a factor of $-(q-q^{-1})^2 [\ell-1-\nu][\ell-2-\nu]$. Since $\ell = 2 \nu
+ 1$, we have $\ell - 1 - \nu = \nu = \alpha+\beta+2$, and $\ell-2-\nu = \alpha+\beta+1$. Combining these factors with $\rho'(\alpha+\beta)$ from \eqref{gamma1gamma3inproof}, we get
$\rho'(\alpha+\beta+2)$ up to powers of $q$, as desired.
	
Matching up the signs and powers of $q$ is a tedious exercise which we leave to the reader.

%\BE{I only did this tedious exercise when $a$ and $b$ are odd... computer verified however}

%========================================================
\subsection{Proof of \eqref{recursiveformula2}}\label{ssec:recursionmain}
%========================================================

%\BE{THIS SECTION COMPLETE!!}

Fix $a > 0$ and $b > 1$ and $0 < k < \ell-1$. Consider \eqref{recursiveformula2}, which states that
\[ \Xi(a,b,2,k) = \Xi(a,b-1,3,k) - z^{2\ell - 3k - 2} \Xi(a,b-1,1,k).\]
Let us prove this formula, in which all factors $\Xi$ are governed by \eqref{eq:mainnasty}.

Whenever $b$ is even, all factors $\Xi$ use the same value of $\alpha$ and $\beta$ and $k$. Thus $\gamma_1$ and $\lambda_2$ and $\thesign$ are the same for all terms. One can also verify that $\gamma_3$ is the same for all three terms, regardless of the parity of $a$.  Also, note that
\begin{equation} \frac{\gamma_2(a,b,2,k)}{\gamma_2(a,b-1,3,k)} = \frac{\gamma_2(a,b,2,k)}{\gamma_2(a,b-1,1,k)} = q^{-(\ell-1)}(q-q^{-1})[\ell-1-k], \end{equation}
regardless of the parity of $a$. No other terms depend on the parity of $a$, so $a$ can be either even or odd below.

The product $\kappa_1 \kappa_2$ has an additional factor of $q^k$ on the left hand side (an extra $q^{2k}$ in $\kappa_2$ and an extra $q^{-1}$ in $\kappa_1$ from the difference in length).  Finally, note that
\begin{equation} \frac{\lambda_1 \lambda_4(a,b)}{\lambda_1 \lambda_4(a,b-1)} = z^{-\ell} p^{-(\beta+1)} p^{\beta+3} = z^{-\ell+1}. \end{equation}
So, dividing by 
\[ (\thesign \gamma_1 \gamma_2 \gamma_3 \kappa_1 \lambda_1 \lambda_2 \lambda_4)(a,b-1,k), \] what remains is
\begin{equation} LHS = q^{-(\ell-1)}(q-q^{-1})[\ell-1-k] q^{k} z^{-\ell+1} z^{\ell} = z q^{k+1-\ell} (q^{\ell-1-k} - q^{k+1-\ell}). \end{equation}
\begin{equation} RHS = z^{\ell}z^{-(\ell-1)} - z^{2 \ell - 3k - 2} z^{\ell} = z(1 - z^{3 \ell - 3k - 3}) = z(1 - q^{2(k+1-\ell)}). \end{equation}
The equality of the two sides is now evident.

Now we consider the case when $b = 2 \beta + 1$ is odd. The $\Xi$ values on the right side use the value $\beta-1$ instead of $\beta$, and we have a fair bit of work to do.

Let us first compare the values of some $\lambda$ variables. We have
\begin{equation} \frac{\lambda_1 \lambda_3 \lambda_4(a,b)}{\lambda_1 \lambda_3 \lambda_4(a,b-1)} = z^{\beta-1} z^{-\ell} \frac{p^{-(\beta+1)(\ell + 3 \beta)}}{p^{-\beta(\ell+3 \beta - 4)}} z^{\ell+1} \frac{p^{(\beta+3)(\phi-1)}}{p^{(\beta+2)(\phi)}}  =  p^{-6 \beta -\ell - 3 + \phi}. \end{equation}
Also, the ratio of $\lambda_2$ variables is $z$ if $a$ is odd and $z^{-1}$ if $a$ is even. Since $b$ is odd, $\phi$ indicates the parity of $a$, so that
\begin{equation} \frac{\lambda_2(a,b)}{\lambda_2(a,b-1)} = z^{1 - 2 \phi}. \end{equation}

Now compare the values of the $\kappa$ variables. We have
\begin{equation} \frac{\kappa_1 \kappa_2(a,b,2,k)}{\kappa_1 \kappa_2(a,b-1, \ne 2,k)} = \frac{q^{2k} q^{k(k-\beta-\ell)}}{q^{k(k-\beta-\ell+2)}} = 1. \end{equation}

Consider the case when $a$ is even and $b$ is odd. Note that the formula for $\gamma_2$ on the right hand side uses the value $\nu-1$ rather than $\nu$. So $[\nu-k]$ divides $\gamma_2(a,b-1,\ne 2,k)$. Dividing both sides by 
\[ (\thesign \gamma_1 \lambda_1 \lambda_2 \lambda_3 \lambda_4 \kappa_1 \kappa_2)(a,b-1,\ne 2,k) \cdot (q-q^{-1})[\nu-k],\] the remainder is
\begin{subequations} \label{LHSRHS1}
\begin{equation}  LHS = (1-q^{-2\beta}) q^{-\nu} \magic(\nu,k,\beta,0) p^{ -6  \beta - \ell - 4} z^{\ell}. \end{equation}
\begin{align}  \nonumber RHS = q^{-(\ell+\nu-3)} & (q-q^{-1}) [k-1] q^{\beta-1} \magic(\nu-1,k-1,\beta-1,0) z^{-(\ell-1)} \\ & - z^{2\ell - 3k - 2}q^{-(\ell-1)}(q-q^{-1})[k-\nu+1]\magic(\nu-1,k,\beta-1,+1). \end{align}
\end{subequations}
Note that the sign change in $\thesign$ is cancelled by the difference between $[\nu-k]$ and $[k-\nu]$. The extra factor of $(1-q^{-2\beta})$ on the LHS comes from $\gamma_1$.

Dividing both sides by $(q-q^{-1}) p^{-3\beta + 5 \nu - 4}$, some elementary manipulation (using $\ell = 2\nu$) yields
\begin{equation} LHS = [\beta] \magic(\nu,k,\beta,0). \end{equation}
\begin{align}  \nonumber RHS = [k-1] & \magic(\nu-1,k-1,\beta-1,0)\\ &  - q^{-3\nu +2k - \beta +1}[k-\nu+1]\magic(\nu-1,k,\beta-1,+1). \end{align}
The equality of these two sides is Lemma \ref{lem:magicmanip1}.

Finally, consider the case when $a$ and $b$ are both odd. Dividing both sides by
\[  (\thesign \gamma_1 \lambda_1 \lambda_2 \lambda_3 \lambda_4 \kappa_1 \kappa_2)(a,b-1,\ne 2,k),\] the remainder is
\begin{equation} LHS = (-1) (1-q^{-2 \beta}) \magic(\nu,k,\beta,-1) p^{-\ell -6  \beta -1}  z^{\ell}. \end{equation}
\begin{align} \nonumber RHS = q^{-(\ell-2)} & (q-q^{-1})[1-k] q^{\beta-1} \magic(\nu-1,k-1,\beta-1,-1)z^{-(\ell-1)}\\ & - z^{2\ell - 3k - 2} q^{-(\nu-1)} (q-q^{-1}) [\nu-1-k] \magic(\nu-1,k,\beta-1,0). \end{align}
The main difference between these equations and \eqref{LHSRHS1} is (a sign, and) that all offsets have been reduced by one. Dividing both sides by $-(q-q^{-1}) p^{-3\beta + 2 \nu - 2}$, some elementary manipulation (using $\ell = 2 \nu -1$) yields
\begin{equation} LHS = [\beta] \magic(\nu,k,\beta,-1). \end{equation}
\begin{align} \nonumber RHS =  [k-1] & \magic(\nu-1,k-1,\beta-1,-1) \\ & + q^{-3\nu +2k - \beta +3} [\nu-1-k] \magic(\nu-1,k,\beta-1,0). \end{align}
Again, the equality of these two sides is Lemma \ref{lem:magicmanip1}.

%% file: FormulaROUpapertwo.tex
%!TEX root = QFrob2.tex

%%%%%%%%%%%%%%%%%%%%%%%%%%%%%%%%%%%%%%%%%%%%%%%%%%%%%%%%%
%========================================================
\section{Evaluation at a root of unity} \label{sec:closedformularou}
%========================================================
%%%%%%%%%%%%%%%%%%%%%%%%%%%%%%%%%%%%%%%%%%%%%%%%%%%%%%%%%

Fix $m \ge 2$. Let $\ze$ be a primitive $3m$-th root of unity. We are interested in computing the scalar $\xi_m(a,i) \in \CC$, defined as
\begin{equation} \xi_m(a,i) := \Xi(a,3m-a-1,i,2m)(\ze) \end{equation} That is, $\xi_m(a,i)$ the scalar obtained by applying $\pa_{(a,3m-a-1,i)}$ to the \emph{staircase polynomial}
\begin{equation} \stair = x_1^{2m} x_2^m, \end{equation}
and then evaluating the result under the specialization $z \mapsto \ze$. A loose formula for $\xi(a,i)$ was already discussed in \eqref{eq:loose}.

Recall that the scalar $\Xi(a,3m-a-1,i,2m)$ lives in $\ZZ[z^{\pm 1}]$, though it is expressed as a product of factors e.g. $\gamma_j$ which live within the larger ring $\ZZ[p^{\pm 1}]$. When we specialize $z \mapsto \ze$, we also specialize $p$ to a primitive $6m$-th root of unity. We use frequently that $p^{3m} = -1 = q^m$.

\begin{lem} \label{lem:whenxinonzero} We have $\xi_m(a,i) = 0$ if either $a \ge 2m+1$ or $a \le m-2$. \end{lem}

\begin{proof}[Proof 1] Using \cite[Corollary 4.40]{EJY1}, $\pa_{(a,b,i)} = 0$ under these conditions because $\un{w}(a,b,i)$ contains a cyclic word of length at least $2m+2$. \end{proof}

\begin{proof}[Proof 2] If $a \ge 2m+1$ then $\alpha \ge m$ in either \eqref{alphabetadef} or \eqref{alphabetadefklen}. If $a \le m-2$ then $b \ge 2m+1$ and $\beta \ge m$ in \eqref{alphabetadef}.
Either way, in all of our formulas, $\rho(c)$ will divide $\Xi(a,b,i,2m)$ for some $c \ge m$. But $q^{m} - q^{-m} = 0$ divides $\rho(c)$. \end{proof}

A consequence of this lemma is that either $\Xi(a,3m-a-1,i,2m)$ is zero or $(a,3m-a-1,i,2m)$ is within the standard regime. Throughout this chapter, we define $\alpha$ and $\beta$
and $\nu$ as in Theorem \ref{thm:mainnasty}, and use the notation from that chapter. For brevity we write e.g. $\gamma_j$ to represent $\gamma_j(a,3m-a-1,i,2m) \in \ZZ[p^{\pm 1}]$
after specializing $p$ to the appropriate root of unity.

%========================================================
\subsection{Precise statement of the result} \label{ssec:evalatrou}
%========================================================

We introduce two more variables, $d$ and $\bottom$, which make it more convenient to state our results. Both are approximately half of $m$. The variable $\bottom$ represents the lower bound on $\alpha$ or $\beta$ for $\xi_m(a,i)$ to be nonzero.

\begin{notation} \label{notation:bottom} Let $m = 2d$ or $m = 2d+1$ depending on parity. Unless $m$ and $a$ and $b$ are odd, set $\bottom = d-1$. If $m$ and $a$ and $b$ are odd, set $\bottom = d$. \end{notation}

\begin{lem} One has $\alpha + \beta = m-1+\bottom$, or equivalently,
\begin{equation} \label{theyaddup} (\alpha - \bottom) + (\beta - \bottom) = m-1-\bottom. \end{equation}
The condition that $\alpha, \beta \le m-1$ is equivalent to the condition that $\bottom \le \alpha \le m-1$ or the condition that $\bottom \le \beta \le m-1$. \end{lem}

\begin{proof} A boring verification. \end{proof}

% 	If $m = 2d$ is even then either $a$ is odd and $b$ even, or $a$ is even and $b$ odd. Either way, $6d = 3m = a+b+1 = 2\alpha + 2 \beta + 4 = 2 \nu$, so $\nu = \frac{3m}{2} = 3d$. In particular, $\alpha + \beta = \nu - 2 = 3d-2$ so $\beta > m-1 = 2d-1$ if and only if $\alpha < d-1 = \bottom$. Also, $(\alpha - \bottom) + (\beta - \bottom) = 3d-2-2(d-1) = d = m - 1 - \bottom$.
%
% If $m = 2d+1$ is odd and $a$ and $b$ are even, then $6d+3 = 3m = a+b+1 = 2\alpha + 2 \beta + 5 = 2 \nu + 1$, so $\nu = \frac{3m-1}{2} = 3d+1$. In particular, $\alpha + \beta = \nu - 2 = 3d-1$, so $\beta > m-1 = 2d$ if and only if $\alpha < d-1 = \bottom$. Also, $(\alpha - \bottom) + (\beta - \bottom) = 3d-1-2(d-1) = d+1 = m - 1 - \bottom$.
%
% If $m = 2d+1$ is odd and $a$ and $b$ are odd, then $6d+3 = 3m = a+b+1 = 2\alpha + 2 \beta + 3 = 2 \nu - 1$, so $\nu = \frac{3m+1}{2} = 3d+2$. In particular, $\alpha + \beta = \nu-2 = 3d$, so $\beta > m-1 = 2d$ if and only if $\alpha < d = \bottom$. Also, $(\alpha - \bottom) + (\beta - \bottom) = 3d-2d = d = m - 1 - \bottom$. \end{proof}

\begin{rem} In the loose formula \eqref{eq:loose}, the integer $t$ is $m-1-\bottom$, and the integer $c$ is $\beta - \bottom$. That is, our formula for $\xi_m(a,i)$ will use the quantum number ${m-1-\bottom \brack \beta - \bottom}$. By \eqref{theyaddup} we have
\begin{equation} {m-1-\bottom \brack \alpha - \bottom} = {m-1-\bottom \brack \beta - \bottom} \end{equation}
which explains the symmetry in our results. \end{rem}

\begin{notation} We say that $\alpha$ (resp. $\beta$) is in the \emph{nonzero range} if $\bottom \le \alpha \le m-1$. \end{notation}

\begin{lem}\label{lem:whengamma1nonzero} The scalar $\gamma_1(\ze)$ is nonzero if and only if $\bottom \le \alpha \le m-1$. \end{lem}
	
\begin{proof} See proof 2 of Lemma \ref{lem:whenxinonzero}. \end{proof}
	
 % Recall that $\gamma_1 = \rho'(\alpha) \rho'(\beta)$, which equal to $[\alpha]! [\beta]!$ up to a power of $q$ and a power of $(q-q^{-1})$. After specialization, $q$ and $q - q^{-1}$ are sent to units, as is $[c]$ whenever $c < m$. However, $[m] \mapsto 0$. Thus $\gamma_1(\ze) = 0$ if and only if $\alpha \ge m$ or $\beta \ge m$. By the previous lemma, $\gamma_1 \ne 0$ if and only if $\bottom \le \alpha \le m-1$. \end{proof}

\begin{thm} \label{thm:whatisw0} Let $m = 2d$ or $m = 2d+1$ depending on parity. Fix $a$ and $i$, let $b = 3m-1-a$, let $\alpha$ and $\beta$ be as in Theorem \ref{thm:mainnasty}, and let $d$ and $\bottom$ be as in Notation \ref{notation:bottom}. Then
\begin{equation} \xi_m(a,i) = (-1)^{d+a+\beta} m^2 z^{2m} q^{\binom{d+1}{2}-\binom{\alpha+1}{2} -\binom{\beta+1}{2}} {m-1-\bottom \brack \beta - \bottom} \cdot \blahblah \end{equation}
where
\begin{equation} \blahblah = \begin{cases}
	p^{-2\beta^2 - 6 \beta - 4} & \text{ if $m$ is even and $a$ is even,} \\ 
	p^{-2\beta^2 - 2 \beta} & \text{ if $m$ is even and $a$ is odd,} \\
	p^{-2 \beta^2 -5\beta -3 + 3d} & \text{ if $m$ is odd and $a$ is even,} \\ 
	p^{-2 \beta^2  - 3 \beta -1} & \text{ if $m$ is odd and $a$ is odd.}
\end{cases} \end{equation}
\end{thm}

%\BE{CONFIRMED WITH CODE!}

The proof of this theorem will comprise the rest of this chapter. The interesting part will be to evaluate the product $\thesign \gamma_1 \gamma_2 \gamma_3$, see Lemma \ref{lem:gammafactorsrou}. 

\begin{rem} One interesting feature of Theorem \ref{thm:whatisw0}, when compared with Theorem \ref{thm:mainnasty}, is that the value of $\xi_m(a,i)$ is independent of $i$! This is yet another ``surprising'' symmetry, though see \cite[Corollary 3.20]{EJY1}. \end{rem}

% \begin{rem} For any given case (i.e. parity of $m$ and $a$ and value of $i$), the exponent of $p$ is a particular quadratic polynomial in the variables $\alpha$ and $d$. To verify
% that two quadratic polynomials agree, one need only check their equality on finitely many values of the inputs. As our computations have been verified by computer for very many
% cases, this is another way to verify Theorem \ref{thm:whatisw0} without doing the terrible exercise. Code can be found at \BE{cite}. \BE{Please make sure this isn't bullshit.} \end{rem}

The variables $\alpha$ and $\beta$ and $d$ in Theorem \ref{thm:whatisw0} are not independent. We have
\begin{subequations} \label{alphabetad}
\begin{eqnarray} \alpha + \beta = 3d, & \; & \text{ when $a$ and $b$ are odd}, \\
\alpha + \beta = 3d-2, & \; & \text{ when $a$ is even and $b$ is odd, or vice versa}, \\
\alpha + \beta = 3d-1, & \; & \text{ when $a$ and $b$ are even}. \end{eqnarray}
\end{subequations}
Thus we can rewrite the formula for $\xi_m$ using only $\beta$ and $d$. Given Theorem \ref{thm:whatisw0}, the proof of the following corollary is straightforward: one rewrites $\binom{\alpha+1}{2}$ as a function of $\beta$ and $d$, and removes extraneous terms using $p^{6m} = 1$. For example, when $m$ is even, $q^{4d} = 1$.

\begin{cor} For $\alpha$ in the nonzero range we have
\begin{equation} \xi_m(a,i) = (-1)^{d+\beta+1} m^2 {m-1-\bottom \brack \beta - \bottom} \cdot \blahblah \end{equation}
where
\begin{equation} \blahblah = \begin{cases}
	p^{\beta^2 + 3 \beta d - d - 1} & \text{ if $m$ is even and $a$ is even,} \\ 
	p^{\beta^2 + 3 \beta d + 4\beta - 7d + 3} & \text{ if $m$ is even and $a$ is odd,} \\
	p^{\beta^2 -9 \beta d - 2\beta - 7d - 2} & \text{ if $m$ is odd and $a$ is even,} \\ 
	p^{\beta^2 -9\beta d - 3 \beta - 7d - 3} & \text{ if $m$ is odd and $a$ is odd.} \end{cases}
\end{equation}
\end{cor}

%========================================================
\subsection{Quantum numbers at roots of unity} \label{ssec:qnumrou}
%========================================================

We discuss some algebraic properties of quantum numbers at a root of unity. Throughout this section, let $q$ be an arbitrary primitive $2m$-th root of unity. Quantum numbers will be
quantum numbers in $q$.

Whenever $q^{2m} = 1$ we have $q^m - q^{-m} = 0$, so 
\begin{equation} [m] = 0. \end{equation}
Because $[2][k] = [k+1] + [k-1]$ for all $k \in \ZZ$, we see that $0 = [m-1] + [m+1]$. Multiplying by $[2]$ again, one can deduce inductively that
\begin{equation} [m-k] = -[m+k]. \end{equation}
Combining this with the standard Clebsh-Gordan rule for multiplication of quantum numbers, one computes that
\begin{equation} [m-1]^2 = 1. \end{equation}
It is a good exercise to confirm that $q^m = -1$ implies that $[m-1]=1$. Again, multiplying by $[2]$ and using induction implies that
\begin{equation} \label{mirror} [m-k] = [k]. \end{equation}
Combining the above, one deduces the periodic formulas
\begin{equation} \label{period} [k+m] = -[k], \qquad [k+2m] = [k]. \end{equation}

\begin{lem} For any $0 \le j \le m$ we have
\begin{equation} \label{2mminus1} {2m-1 \brack j} = (-1)^j. \end{equation}
\end{lem}

\begin{proof} Since $[2m-j] = [-j] = - [j]$, we have
\begin{equation} {2m-1 \brack j} = \frac{[2m-1][2m-2]\cdots [2m-j]}{[1][2] \cdots [j]} = (-1)^j. \end{equation} Our bounds on $j$ remove the possibility of any issues involved with dividing zero by zero. \end{proof}

Multiplying formulas \eqref{mirror} and \eqref{period} by $q-q^{-1}$, we get analogous formulas for $q^k - q^{-k}$. For example,
\begin{equation} \label{mirror2} q^{m-k} - q^{k-m} = q^k - q^{-k}. \end{equation}
Recall the definition of $\rho$ from \eqref{defrho}.
 
\begin{lem} Let $m = 2d$ be even. We have
	\begin{equation} \label{rhom1trigeven} \rho(m-1) = q^{d(m-1)} m. \end{equation}
Since $q^d$ is a primitive $4$-th root of unity, this formula implies $\rho(m-1) = \pm \sqrt{-1} m$.

Now let $m = 2d+1$ be odd. We have
	\begin{equation} \label{rhom1trigodd} \rho(m-1) =(-1)^d m. \end{equation}
\end{lem}

\begin{rem} In either case, a consequence we will not use is $\rho'(m-1) = m$. \end{rem}
%\BE{THIS STATEMENT IS PROBABLY EASIER FROM STANDARD FORMULAS!!!}

\begin{proof} We first prove the result for a particular root of unity, and then address what happens for other roots of unity.

Let $\theta = \pi/m$. One choice for $q$ is $e^{\sqrt{-1} \theta}$. For this choice we have
\begin{equation} q + q^{-1} = 2 \cos(\theta), \qquad q - q^{-1} = 2\sqrt{-1} \sin(\theta). \end{equation}
Similarly, 
\begin{equation} q^k + q^{-k} = 2 \cos(k \theta), \qquad q^k - q^{-k} = 2\sqrt{-1} \sin(k \theta). \end{equation}
We write $\sqrt{-1}$ rather than the complex number $i$ because $i$ is already used for an element of $\Om$.

A standard formula from trigonometry\footnote{We call this standard because it appears on the Wikipedia page listing trigonometric identities, under ``Finite products of trigonometric functions.'' There is a more well-known identity, listed right below it, which states that
\[ \sin(nx) = 2^{n-1} \prod_{c=0}^{n-1} \sin(c \theta + x). \]
Dividing by the $c=0$ term, and then taking the limit as $x$ goes to zero, one obtains the desired result.},
says that
\begin{equation} \prod_{k=1}^{m-1} \sin(k \theta) = \frac{m}{2^{m-1}}. \end{equation}
Multiplying both sides by $(2\sqrt{-1})^{m-1}$ we get
\begin{equation} \label{rhom1trigabort} \rho(m-1) = \prod_{k=1}^{m-1} (q^k - q^{-k}) = (\sqrt{-1})^{m-1} m. \end{equation}

When $m = 2d$ is even, $q^d = \sqrt{-1}$, so we can rewrite this equation as \eqref{rhom1trigeven}. When $m = 2d-1$ is odd, then $(\sqrt{-1})^{m-1} = (-1)^d$ is a sign which does not
depend on the choice of $\sqrt{-1}$. We obtain \eqref{rhom1trigodd}.

To argue that \eqref{rhom1trigeven} and \eqref{rhom1trigodd} hold for arbitrary primitive $2m$-th roots of unity $q$, we can multiply both sides by $q^{\binom{m}{2}}$ to get polynomials in $q$ rather than Laurent polynomials. The result is a polynomial with integral coefficients. Because the Galois group acts transitively on primitive roots of unity, the same polynomial vanishes regardless of the choice of $q$.
\end{proof}

\begin{rem} Note that \eqref{rhom1trigabort} is not independent of the choice of root of unity when $m$ is even. Replacing $q$ with $q^{-1}$ would multiply $\rho(m-1)$ by a sign. \end{rem}
%	, but not the RHS of the equation. The RHS is a Laurent polynomial with complex coefficients, so the above argument involving Galois groups would not apply.  \end{rem}
	
It is well known that primes dividing $2m$ will ramify in the cyclotomic field $\QQ(\ze)$. We can explicitly construct square roots of $m$ (or $m/2$) using values of $\rho$.

\begin{lem} Let $m = 2d$ be even. Then 
\begin{equation} \label{rhom1even} \rho(m-1) = \rho(d) \rho(d-1). \end{equation}
Also, 
\begin{equation} \label{rhom1evenagain} \rho(d) = 2 q^{d} \rho(d-1), \end{equation}
where $q^d = \pm \sqrt{-1}$ is a primitive 4-th root of unity.

Let $m = 2d+1$ be odd. Then 
\begin{equation} \label{rhom1odd} \rho(m-1) = \rho(d)^2. \end{equation} \end{lem}
	
\begin{proof} When $m = 2d$ is even, \eqref{mirror} says that
\begin{equation} \label{mirror3}[d+k] = [d-k], \end{equation}
and \eqref{mirror2} says that
\begin{equation}\label{mirror4} q^{d+k} - q^{-(d+k)} = q^{d-k} - q^{-(d-k)}. \end{equation}
Now \eqref{rhom1even} is an immediate consequence of \eqref{mirror4}. Meanwhile, $\rho(d)/\rho(d-1) = q^d - q^{-d} = 2 q^d$, since $q^{-2d} = q^m = -1$.	

When $m = 2d+1$ is odd, \eqref{mirror2} immediately implies \eqref{rhom1odd}. \end{proof}

\begin{lem} After specialization to a root of unity we have
\begin{equation} \label{gamma1ze} \rho(\alpha)\rho(\beta) = \rho(\bottom) \rho(m-1) {m-1-\bottom \brack \alpha - \bottom}. \end{equation}
\end{lem}

% \begin{lem} Let $q$ be a primitive $2m$-th root of unity. The following equation holds whenever $m = 2d$ is even and $d-1 \le \alpha \le m-1$. Letting $\beta = 3d-2-\alpha$ we have
% \begin{equation}\label{gamma1ze} \rho(\alpha) \rho(\beta) = \rho(d)\rho(d-1)^2 {d \brack \alpha+1-d}. \end{equation}
% \end{lem}
%
% Up to units, this says that $\gamma_1$ is a quantum binomial coefficient times $m \sqrt{\frac{m}{2}}$.

\begin{proof} If $\alpha$ is not in the nonzero range then both sides will vanish. For the rest of this proof we assume that $\alpha$ is in the nonzero range.
	
When $\alpha = \bottom$ then $\beta = m-1$, and the formula \eqref{gamma1ze} is immediate. For the general case, we divide by the known equality for the special case $\alpha = \bottom$. When $\alpha = \bottom + c$ and $\beta = m-1-c$, then
\begin{equation} \label{helloworld}\frac{\rho(\alpha) \rho(\beta)}{\rho(\bottom) \rho(m-1)} = \frac{\prod_{k=1}^{c} (q^{\bottom+k} - q^{-(\bottom+k)})}{\prod_{k=1}^{c} (q^{m-k} - q^{-(m-k)})} = \prod_{k=1}^{c} \frac{[\bottom+k]}{[m-k]}. \end{equation}
Meanwhile,
\begin{equation} \label{helloworld2} {m-1-\bottom \brack \alpha-\bottom} = {m-1-\bottom \brack c} = \prod_{k=1}^{c} \frac{[m-\bottom-k]}{[k]} = \prod_{k=1}^{c} \frac{[\bottom+k]}{[m-k]}. \end{equation}
The last equality holds by applying \eqref{mirror} to both the numerator and denominator. Since \eqref{helloworld} and \eqref{helloworld2} match, we deduce the result. \end{proof}

By combining the results in this section, one can prove that $\gamma_1$ times either $\rho(d)$ or something similar is equal to $m^2 {m-1-\bottom \brack \alpha - \bottom}$ times a power of $q$. The details depend on various parity considerations, and will appear in the proofs below.

%========================================================
\subsection{Magic at a root of unity} \label{ssec:magicrou}
%========================================================

In this section we study the evaluation of $\gamma_3$ at a root of unity. We recall again the $q$-binomial theorem:
\begin{equation} \label{qbinomialthmredux}
\sum_{j=0}^n
{n \brack j} q^{-j(n-1)}t^j
=
\prod_{c=0}^{n-1}\left(1+q^{-2c}t\right).
\end{equation}

We first consider when $m=2d$ is even. Then $k = 4d$ and $\nu = 3d$ (c.f. \eqref{alphabetad}).

Recall that
\[ \term(\nu,k,\beta,\epsilon,j) = {k-1 \brack \beta-j} {\nu-k-1 \brack j} q^{j(-3\nu- 2 \epsilon+2 k)}, \]
and
\[ \magic(\nu, k,\beta,\epsilon) = \sum_{j = 0}^{\beta} \term(\nu,k,\beta,\epsilon,j). \]
Below we assume that $0 \le j \le \beta$, which (when $\beta$ is in the nonzero range) also implies that $j < m$ and $\beta - j < m$.

\begin{lem}  Let $m = 2d$ be even. For $\beta$ in the nonzero range, we have
\begin{equation} \label{magic0} \magic(3d, 4d,\beta,0)(\ze) = (-1)^{\beta +d-1} q^{\binom{d}{2}} \rho(d-1), \end{equation}
\begin{equation} \label{magic1} \magic(3d, 4d,\beta,-1)(\ze) = (-1)^{\beta + d} q^{\binom{d+1}{2}} \frac{\rho(d)}{1-q^2}, \end{equation}
\begin{equation} \label{magic3} q^{\beta} \magic(3d, 4d-1,\beta,-1)(\ze) = (-1)^{\beta + d} q^{\binom{d+1}{2}} \frac{\rho(d)}{1-q^2}. \end{equation}
\end{lem}

\begin{proof} Let us focus on \eqref{magic0} as the other cases are similar. We lay out several steps in the proof, which will be imitated in other cases. Specializing appropriately, we get
\begin{equation} \term(3d,4d,\beta,\epsilon,j) = {2m-1 \brack \beta-j} {-d-1 \brack j} q^{-dj} q^{-2 j\epsilon}. \end{equation}
As an aside, note that $q^d$ is a primitive fourth root of unity, so $q^{-dj}$ is just a fourth root of unity. 

\begin{enumerate} \item By \eqref{2mminus1}, ${2m-1 \brack \beta-j} = (-1)^{\beta -j}$.

\item The usual rules for binomial coefficients with negative values state that
\begin{equation} {-d-1 \brack j} = (-1)^j {d+j \brack j}. \end{equation}

Combining the previous steps we deduce that
\begin{equation} \term(3d,4d,\beta,\epsilon,j) = (-1)^{\beta} {d+j \brack j} q^{-dj} q^{-2 j\epsilon}. \end{equation}

\item Using \eqref{mirror3} we have  ${d+j \brack j} = {d-1 \brack j}$.

\item As a consequence of the previous step, $\term(3d,4d,\beta,\epsilon,j)$ vanishes unless $0 \le j \le d-1$. Note that $d-1 = \bottom \le \beta$, and the condition $j \le \beta$ is now redundant. So, regardless of the value of $\beta$ in the nonzero range, we have
\begin{equation} \label{magicsimplified} \magic(3d, 4d,\beta,\epsilon) =\sum_{j = 0}^{d-1} \term(3d,4d,\beta,\epsilon,j) = (-1)^{\beta} \sum_{j=0}^{d-1} {d-1 \brack j} q^{-dj} q^{-2j\epsilon} \end{equation}

\item Note that $q^{-dj} q^{-2j\epsilon} = q^{-(d-2)j} (q^{-2(1+\epsilon)})^j$. Now we are in a position to apply the $q$-binomial theorem \eqref{qbinomialthmredux} with $n=d-1$ and $t= q^{-2(1+\epsilon)}$. It states that
\begin{equation} \label{refertome}
\sum_{j=0}^{d-1} {d-1 \brack j} q^{-dj} q^{-2j\epsilon} = \prod_{c=0}^{d-2}\left(1+q^{-2(c+1+\epsilon)}\right) = \prod_{c=1+\epsilon}^{d-1+\epsilon}\left(1+q^{-2c}\right).
\end{equation}

\item Recall that
\[ \rho(d-1) = \prod_{c=1}^{d-1} (q^c - q^{-c}). \]
We wish to add signs to the right side of \eqref{refertome} to obtain an expression more like $\rho$. We do this by using the fact that $q^{2d} = -1$, and by reindexing the sum.
\begin{equation} \label{reindexproduct} \prod_{c=1+\epsilon}^{d-1+\epsilon}\left(1+q^{-2c}\right) = \prod_{c=1+\epsilon}^{d-1+\epsilon}\left(1-q^{2(d-c)}\right) = \prod_{c=1-\epsilon}^{d-1-\epsilon}  \left(1-q^{2c}\right). \end{equation}
When $\epsilon = 0$, we get
\begin{equation} \label{epsilon0laststep} \prod_{c=1}^{d-1}  \left(1-q^{2c}\right) = (-1)^{d-1} q^{\binom{d}{2}} \rho(d-1).\end{equation}

\item This step is meant to encompass the algebra needed to reach the final desired form, which in this case is almost nothing. Combining \eqref{epsilon0laststep} with \eqref{refertome} and \eqref{magicsimplified}, we deduce \eqref{magic0}.
\end{enumerate}

To prove \eqref{magic1} we follow the same steps. In the last step, $\epsilon = -1$, and the product ranges from $c=2$ to $d$. Multiplying and dividing by the $c=1$ term, the right
side of \eqref{reindexproduct} is $(-1)^{d} q^{\binom{d+1}{2}} \frac{\rho(d)}{1-q^2}$.

For \eqref{magic3} we instead need to examine
\begin{equation} \term(3d,4d-1,\beta,-1,j) = {2m-2 \brack \beta-j} {-d \brack j} q^{-dj}. \end{equation}
Minor modifications of the first three steps yield
\begin{equation} {2m-2 \brack \beta-j} = (-1)^{\beta-j} [\beta-j+1], \qquad {-d \brack j} = (-1)^j {d \brack j}, \end{equation}
so that	
\begin{equation} \term(3d,4d-1,\beta,-1,j) = (-1)^{\beta} [\beta+1-j] {d \brack j} q^{-dj}. \end{equation}
	
The fourth step is the analysis of which terms (indexed by $j$) are nonzero. This time the result is evidently zero for $j > d$, but the $j=d$ term may be nonzero. Now our sum will range from $0$ to $d$, regardless of the value of $\beta$. A small wrinkle appears and is immediately smoothed: we also assumed $j \le \beta$, and $\beta = d-1$ is permitted, which would rule out the $j=d$ term. Thankfully, in the special case when $\beta = d-1$ and $j=d$, $[\beta + 1-j] = 0$ and this term vanishes. Ultimately, our replacement for \eqref{magicsimplified} is
\begin{equation} \label{alongtheway} (-q)^{\beta} \magic(3d,4d-1,\beta,-1) = q^{\beta} \sum_{j=0}^{d} [\beta+1-j] {d \brack j} q^{-dj}, \end{equation}
for all $\beta$ in the nonzero range.

Now a new argument appears. We claim that the right side of \eqref{alongtheway} is actually independent of $\beta$. The difference between this formula for $\beta$ and for
$\beta-1$, after some manipulation, is
\begin{equation} \sum_{j=0}^d {d \brack j} q^{-dj} (q^{\beta} [\beta+1-j] - q^{\beta-1}[\beta-j]) = q^{2 \beta} \sum_{j=0}^d {d \brack j} q^{-(d+1)j}. \end{equation}
Using the $q$-binomial theorem \eqref{qbinomialthmredux} with $n=d$ and $t = q^{-2}$, we get the product
\[ \prod_{c=0}^{d-1} \left( 1 + q^{-2c-2} \right). \]
When $c = d-1$, the factor is $1+q^{-2d} = 0$.

To evaluate the right side of \eqref{alongtheway}, we may choose our favorite value of $\beta$, which will be $\beta = 2m-1$. Then $[\beta+1-j] = [2m-j] = -[j]$, and
$[j]{d \brack j} = [d]{d-1 \brack j-1}$. So with a little reindexing we have
\begin{align}\nonumber (-q)^{\beta} \magic(3d,4d-1,\beta,-1) & = -q^{2m-1}[d] \sum_{j=0}^{d} {d-1 \brack j-1} q^{-dj}  \\ & = -q^{2m-1}[d] \sum_{j=0}^{d-1} {d-1 \brack j} q^{-dj}q^{-d}. \end{align}

Now we return to the previous flow of the argument. Using the $q$-binomial theorem once more, exactly as in \eqref{refertome} with $\epsilon = 0$, and then adding signs using $q^{2d} = -1$, we get
\begin{equation} (-q)^{\beta} \magic(3d,4d-1,\beta,-1) = -q^{2m-1-d} (-1)^{d-1} q^{\binom{d}{2}} [d] \rho(d-1). \end{equation}
Noting that $q^{2m} = 1$ and $-q^{-d} = q^d$ and $[d] = \frac{q^d - q^{-d}}{q-q^{-1}}$, we get
\begin{equation} (-q)^{\beta} \magic(3d,4d-1,\beta,-1) = q^{-1} (-1)^{d-1} q^{\binom{d}{2} + d}  \frac{\rho(d)}{q-q^{-1}}. \end{equation}
We easily simplify to obtain \eqref{magic3}.
\end{proof}

Now we switch to the case when $m$ is odd. Now $k = 2m = 4d+2$. If $a$ and $b$ are odd then $\nu = 3d+2$, while if $a$ and $b$ are even then $\nu = 3d+1$, see \eqref{alphabetad}. We follow the same steps as in the previous lemma. The major differences between the arguments below and above are:
\begin{itemize}
	\item The mirror for symmetry between quantum numbers has shifted, so we use \eqref{mirror} instead of \eqref{mirror3} in step 3.
	\item Now $q^{2d+1} = -1$ instead of $q^{2d}$, which adjusts step 5 where we ``add signs and reindex.''
\end{itemize}

\begin{lem} Let $m = 2d+1$ be odd. For $\beta$ in the nonzero range\footnote{For \eqref{magicoddodd} this implies $\beta \ge \bottom = d$. For the other equations, $\beta \ge \bottom = d-1$. This matches the conditions needed for this instance of the function $\magic$ to relate to $\gamma_3$.}, we have 
\begin{equation} \label{magicoddodd} \magic(3d+2,4d+2,\beta,-1) = (-1)^{\beta+d}q^{\binom{d+1}{2}} \rho(d). \end{equation}
\begin{equation} \label{magiceveneven1} \magic(3d+1, 4d+2,\beta,+1) = (-1)^{\beta+d-1} q^{\binom{d}{2}} \rho(d-1). \end{equation}
\begin{equation} \label{magiceveneven2} \magic(3d+1, 4d+2, \beta,0) = (-1)^{\beta + d} q^{\binom{d+1}{2}} \frac{\rho(d)}{1-q^2}. \end{equation}
\begin{equation} \label{magiceveneven3} q^{\beta} \magic(3d+1, 4d+1, \beta,0) = (-1)^{\beta+d} q^{\binom{d+1}{2}} \frac{\rho(d)}{1-q^2}. \end{equation}	
\end{lem}

\begin{proof}
We begin with \eqref{magicoddodd}. We have
\begin{equation} \term(3d+2,4d+2,\beta,-1) = {2m-1 \brack \beta-j} {-d-1 \brack j} q^{-dj}. \end{equation}
In the third step we use \eqref{mirror} to identify ${d+j \brack j}$ with ${d \brack j}$. We get
\begin{equation} \term(3d+2,4d+2,\beta,-1) = (-1)^{\beta} {d \brack j} q^{-dj}. \end{equation}
This term vanishes for $j > d$. Note that this formula is being applied when $a$ and $b$ are odd, so that $\bottom = d$. The condition $j \le \beta$ is therefore redundant.

Using the $q$-binomial theorem with $n=d$ and $t = q^{-1}$, we have 
\begin{equation} \sum_j {d \brack j} q^{-dj} = \prod_{c=0}^{d-1} (1 + q^{-2c-1}). \end{equation}
Now $q^{2d+1} = -1$, so $q^{-2c-1} = -q^{2(d-c)}$. Reindexing, we deduce altogether that
\begin{equation} \magic(3d+2,4d+2,\beta,-1) =(-1)^{\beta} \prod_{c = 1}^d (1 - q^{2c})= (-1)^{\beta+d}q^{\binom{d+1}{2}} \rho(d). \end{equation}
	
Now consider \eqref{magiceveneven1} and \eqref{magiceveneven2}.
\begin{equation} \term(3d+1, 4d+2, \beta, \epsilon) = {2m-1 \brack \beta-j} {-d-2 \brack j} q^{-j(d-1)} q^{-2\epsilon j}. \end{equation}
Here are the variations. We use \eqref{mirror} we identify ${-d-2 \brack j}$ with $(-1)^j {d-1 \brack j}$. We use the $q$-binomial theorem with $n=d-1$ and $t = q^{-1-2\epsilon}$. The result is that
\begin{equation} \magic(3d+1,4d+2,\beta,\epsilon) = (-1)^{\beta} \prod_{c = 2 - \epsilon}^{d-\epsilon} (1 - q^{2c}). \end{equation}
The process to derive \eqref{magiceveneven1} and \eqref{magiceveneven2} from here is straightforward.

% DOING FOR MYSELF: \eqref{magiceveneven1} and \eqref{magiceveneven2}
% \begin{equation} \term(3d+1, 4d+2, \beta, \epsilon) = {2m-1 \brack \beta-j} {-d-2 \brack j} q^{-j(d-1)} q^{-2\epsilon j}. \end{equation}
% \begin{equation} = (-1)^{\beta} {d-1 \brack j} q^{-j(d-2)} (q^{-1-2 \epsilon})^j. \end{equation}
% \begin{equation} = (-1)^{\beta} \prod_{c = 0}^{d-2}	(1+q^{-2c-1-2\epsilon}) = (-1)^{\beta}  \prod_{c = 0}^{d-2}	(1 - q^{2(d-c - \epsilon)}) \end{equation}
% \begin{equation} = (-1)^{\beta} \prod_{c = 2 - \epsilon}^{d-\epsilon} (1 - q^{2c}). \end{equation}

Finally, let us analyze \eqref{magiceveneven3} following the same steps as for \eqref{magic3}. 	Our replacement for \eqref{alongtheway} is
\begin{equation} \label{alongtheway2} (-q)^{\beta} \magic(3d+1,4d+1,\beta,0) =  q^{\beta} \sum_{j=0}^{d} [\beta+1-j] {d \brack j} q^{-j(d+1)}. \end{equation}
% \begin{equation} (-q)^{\beta} \term(3d+1, 4d+1, \beta,0,j) = (-q)^{\beta} {2m-2 \brack \beta-j} {-d-1 \brack j} q^{-j(d+1)}. \end{equation}
% \begin{equation} = q^{\beta} [\beta-j+1] {d \brack j} q^{-j(d+1)}. \end{equation}
Again, the right side is independent of $\beta$, as the difference between two successive terms is
\begin{equation} \sum_{j=0}^d {d \brack j} q^{-j(d+1)} (q^{\beta} [\beta+1-j] - q^{\beta-1}[\beta-j]) = q^{2 \beta} \sum_{j=0}^d {d \brack j} q^{-j(d+2)}. \end{equation}
Using the $q$-binomial theorem \eqref{qbinomialthmredux} with $n=d$ and $t = q^{-3}$, this is 
\[ \prod_{c = 0}^{d-1} \left( 1 + q^{-2c-3} \right). \]
When $c = d-1$, the factor is $1 + q^{-2d-1} = 0$.

Plugging in $\beta = 2m-1$ and manipulating as before, we have
\begin{equation} (-q)^{\beta} \magic(3d+1,4d+1,\beta,0) = -q^{2m-1} [d] \sum_{j=0}^{d-1} {d-1 \brack j} q^{-j(d+1)} q^{-(d+1)}. \end{equation}
Using the $q$-binomial theorem with $n=d-1$ and $t = q^{-3}$, then adding signs using $q^{2d+1} = -1$, we get
\begin{align}\nonumber -q^{2m-d-2} [d] \prod_{c=0}^{d-2} (1 + q^{-2c-3}) & = -q^{2m-d-2} [d] \prod_{c=0}^{d-2} (1 - q^{2(d-c-1)})  \\ & = -q^{2m-d-2} [d] \prod_{c=1}^{d-1} (1 - q^{2c}). \end{align}
Noting $q^{2m} = 1$ and $q^{-d-1} = -q^{d}$ we get
\begin{equation} (-q)^{\beta} \magic(3d+1,4d+1,\beta,0) = q^{-1} (-1)^{d-1} q^{\binom{d+1}{2}}  \frac{\rho(d)}{q-q^{-1}}. \end{equation}
Again, it is easy to simplify to obtain \eqref{magiceveneven3}. \end{proof}

%========================================================
\subsection{Evaluation of all $\gamma$ factors} \label{ssec:gamma}
%========================================================

\begin{lem} \label{lem:gammafactorsrou} Fix $a$ and $i$, let $b = 3m-a-1$, and let $\gamma_j = \gamma_j(a,b,i,2m)$. Use the same notation for $\alpha, \beta,\nu, d, \bottom$ as above. For $\alpha$ in the nonzero range we have
\begin{equation} \label{gammafactorsrou} \thesign \gamma_1 \gamma_2 \gamma_3 = (-1)^{d+b+1} m^2 {m-1-\bottom \brack \alpha - \bottom}\cdot \begin{cases} q^{\binom{d}{2} - \binom{\alpha+1}{2} - \binom{\beta+1}{2}} & \text{ if $a$ and $b$ are even}, \\ q^{\binom{d+1}{2} - \binom{\alpha+1}{2} - \binom{\beta+1}{2}} & \text{ else}, \end{cases} \end{equation}
\end{lem}

\begin{proof} In all cases, $\thesign = (-1)^{\beta}$. Several times in this proof we tacitly use \eqref{defrhoprime} to compare $\rho$ and $\rho'$.
	
First suppose either $a$ is even and $b$ is odd, or $a$ is odd and $b$ is even and $i=1$. So $m = 2d$ and $\nu = 3d$ and $q^{2d} = -1$, which we use below without comment. Thus
\begin{equation} \gamma_2 = q^{-3d} (q-q^{-1})[4d-3d] (-1)^a = (q^{-2d} - q^{-4d}) \cdot (-1)^{a} = 2 \cdot (-1)^{a+1}. \end{equation}
Combining \eqref{gamma1ze} and \eqref{magic0} we get
\begin{equation} \thesign \gamma_1 \gamma_2 \gamma_3 = 2 \cdot (-1)^{a+d} q^{\binom{d}{2} - \binom{\alpha+1}{2} - \binom{\beta+1}{2}} \rho(d-1)^2 \rho(m-1) {m-1-\bottom \brack \alpha - \bottom}. \end{equation}
Using \eqref{rhom1even} and \eqref{rhom1evenagain} and then \eqref{rhom1trigeven} gives
\begin{equation} 2 \rho(d-1)^2 \rho(m-1) = q^{-d} \rho(m-1)^2 = q^{-d} m^2 q^{2d(m-1)} = q^{d} m^2 \end{equation}
%	\frac{\rho(m-1)^2}{(q^d - q^{-d})} = \frac{(-1)^{m-1} m^2}{-2 q^{-d}} = \frac{m^2}{2 q^{-d}}. \end{equation}
From here \eqref{gammafactorsrou} follows easily.

Now suppose $a$ is odd and $b$ is even and $i=2$. Then
\begin{equation} \gamma_2 = q^{-6d+1}(q-q^{-1})[2d-1] = (1-q^{-4d+2}) = 1-q^2. \end{equation}
This cancels the denominator in \eqref{magic1}. We have
\begin{equation} \thesign \gamma_1 \gamma_2 \gamma_3 = (-1)^{d} q^{\binom{d+1}{2} - \binom{\alpha+1}{2} - \binom{\beta+1}{2}} \rho(d)\rho(d-1) \rho(m-1) {m-1-\bottom \brack \alpha - \bottom}. \end{equation}
Using \eqref{rhom1even} and \eqref{rhom1trigeven} gives
\begin{equation} \label{barfoob} \rho(d) \rho(d-1) \rho(m-1) = \rho(m-1)^2 = (-1)^{m-1} m^2 = -m^2. \end{equation}
From here \eqref{gammafactorsrou} follows easily.

The case where $a$ is odd and $b$ is even and $i=3$ is very similar.

% Now suppose $a$ is odd and $b$ is even and $i=3$. Under these hypotheses $\gamma_2 = q^{-6d+1}(q-q^{-1})[1-4d] = -(q^{2d} - q^{2-2d}) = 1 - q^2$. Using \eqref{magic3} and \eqref{gamma1ze} and \eqref{barfoob} we get
% \begin{equation} \thesign \gamma_1 \gamma_2 \gamma_3 = (-1)^{d+1} m^2 q^{\binom{d+1}{2}- \binom{\alpha+1}{2} - \binom{\beta+1}{2}}  {m-1-\bottom \brack \alpha - \bottom}. \end{equation}
% This matches \eqref{gammafactorsrou}.

Now suppose $a$ and $b$ are both odd, so $m = 2d+1$ and $\nu = 3d+2$. Now $\gamma_2 = 1$. Using \eqref{magicoddodd} and noting that $\bottom = d$ we have
\begin{equation} \thesign \gamma_1 \gamma_2 \gamma_3 = (-1)^{d}q^{\binom{d+1}{2}- \binom{\alpha+1}{2} - \binom{\beta+1}{2}} \rho(d)^2 \rho(m-1) {m-1-\bottom \brack \alpha - \bottom}. \end{equation}
Now we use \eqref{rhom1odd} and \eqref{rhom1trigodd} to deduce that
\begin{equation} \label{barfooa} \rho(d)^2 \rho(m-1) = m^2,\end{equation}
from which we easily deduce \eqref{gammafactorsrou}

Suppose that $a$ and $b$ are both even, and $i=1$.  Then $m = 2d+1$ and $\nu = 3d+1$ and $q^m = -1$. So $q^d = - q^{-(d+1)}$. We have
\begin{equation} \gamma_2 = q^{-3m}(q^{d+1} - q^{-(d+1)})(q^{-d} - q^{d}) = (q^d - q^{-d})^2. \end{equation}
Using \eqref{magiceveneven1} and \eqref{gamma1ze} we have
\begin{equation} \thesign \gamma_1 \gamma_2 \gamma_3 =  (-1)^{d-1} q^{\binom{d}{2} - \binom{\alpha+1}{2} - \binom{\beta+1}{2}} (q^d - q^{-d})^2 \rho(d-1)^2 \rho(m-1) {m-1-\bottom \brack \alpha - \bottom}. \end{equation}
Note that $(q^d-q^{-d})^2 \rho(d-1)^2 = \rho(d)^2$. Using \eqref{barfooa}, we get the desired result.

The manipulations are very similar for the remaining two cases, where $a$ and $b$ are both even, and $i=2$ or $i=3$. The fact that $\binom{d+1}{2} - d = \binom{d}{2}$ is also used.

\end{proof}

%
% \BE{This uses old notation!!}
%
% Suppose that $a$ and $b$ are both even, and $i=2$. Now $\gamma_2 = -(1-q^{-2d})(1 - q^2) = -q^{-d}(q^d - q^{-d})(1-q^2)$, and the factor of $1-q^2$ cancels the denominator of \eqref{magiceveneven2}. So we have
% \begin{equation} \thesign \gamma_1 \gamma_2 \gamma_3 \lambda_6  = (-1)^{d+1} q^{\binom{d+1}{2} - d - \binom{\alpha+1}{2} - \binom{\beta+1}{2}} (q^d - q^{-d})\rho(d) \rho(d-1) \rho(m-1) {m-1-\bottom \brack \alpha - \bottom}. \end{equation}
% Again using \eqref{barfooa} we get the desired result.
%
% Suppose that $a$ and $b$ are both even, and $i=3$. Now $\gamma_2 = (1-q^{-2})(q^{2d}-1) = -q^{d-2}(1-q^2)(q^d - q^{-d})$, and $\lambda_6 = y^{\nu-1} q^{\beta} = q^{-2d+2} q^{\beta}$. Thus $\gamma_2 \lambda_6 = -q^{-d}(1-q^2)(q^d - q^{-d})$. Using \eqref{magiceveneven3} and \eqref{gamma1ze} we get
% \begin{equation} \thesign \gamma_1 \gamma_2 \gamma_3 \lambda_6  =  (-1)^{d+1} q^{\binom{d+1}{2} - d- \binom{\alpha+1}{2} - \binom{\beta+1}{2}} (q^d - q^{-d}) \rho(d) \rho(d-1) \rho(m-1) {m-1-\bottom \brack \alpha - \bottom}. \end{equation}
% Again, using \eqref{barfooa} we get the desired result.

%========================================================
\subsection{Evaluation of remaining factors at a root of unity} \label{ssec:therest}
%========================================================

We now finish the proof of Theorem \ref{thm:whatisw0} by examining the remaining factors in \eqref{eq:mainnasty}.

\begin{lem} \label{lem:easyfactorsrou} When $k = 2m$ and $z \mapsto \ze$, then
	\begin{equation} \kappa_2 = \lambda_5 = 1, \end{equation}
	\begin{equation} \kappa_1 = p^{4m}, \end{equation}
\end{lem}

\begin{proof} This is easy, since when $k=2m$ and $\ell = 3m$, then $q^k = z^{\ell} = 1$. \end{proof}
	
This factor of $p^{4m}$ explains the factor of $p^{4m}$ in Theorem \ref{thm:whatisw0}.

Let us examine the signs. From \eqref{gammafactorsrou} we have a factor of $(-1)^{d + b + 1}$. Within $\lambda_1$ we also have two terms which become signs when evaluated at a root of unity. Because $p^{\ell} = -1$, we deduce that $p^{-\ell(\beta+1)} = (-1)^{\beta+1}$. Also, we claim that
\begin{equation} z^{\binom{\ell+1}{2}} = p^{-3m(3m+1)} = (-1)^{m+1}. \end{equation} 
This is because $p^{3m} = -1$. Altogether, the sign is
\begin{equation} (-1)^{d + b + \beta + m + 3} = (-1)^{d + \beta + a}, \end{equation}
since $m = a+b+1$. Below, we ignore all terms which have been encapsulated into the sign above.

When $a$ is even and $b$ is odd, one can compute that
\begin{equation} \lambda_1 \lambda_2 \lambda_3 \lambda_4 = \pm p^{-2 \beta^2 - 6 \beta - 4}. \end{equation}
Altogether, we derive using \eqref{eq:mainnasty} that
\begin{equation} \xi_m = (-1)^{d + \beta + a} m^2 q^{\binom{d+1}{2} - \binom{\alpha+1}{2} - \binom{\beta+1}{2}} p^{4m} {m-1-\bottom \brack \beta - \bottom} p^{-2\beta^2 - 6 \beta - 4}, \end{equation}
as desired.
% Using that $\alpha = 3d - \beta - 2$, a tedious calculation yields
% \begin{equation} q^{\binom{d+1}{2}-\binom{\alpha+1}{2} -\binom{\beta+1}{2}} p^{4m -2 \beta^2 - 6 \beta - 4} = p^{12 d^2 - 9 d \beta + \beta^2 - 7d - 1}. \end{equation}
% Noting that $p^{12d} = 1$ and $p^{6d} = -1$, we have
% \begin{equation} \xi_m = (-1)^{d + \beta +1} m^2 p^{\beta^2 +  3 d \beta - d - 1} {m-1-\bottom \brack \alpha - \bottom}.
% \end{equation}

When $a$ is odd and $b$ is even, one can compute that
\begin{equation} \lambda_1 \lambda_2 \lambda_3 \lambda_4 = \pm p^{-2 \beta^2 - 2 \beta}, \end{equation}
which is the same as before times $p^{4 \beta + 4}$.

When $a$ and $b$ are odd, one can compute that
\begin{equation} \lambda_1 \lambda_2 \lambda_3 \lambda_4 = \pm p^{-2 \beta^2  - 3 \beta -1}.\end{equation}

% Using that $\alpha = 3d-\beta$,
% \begin{equation}
% p^{-2 \beta^2  - 3 \beta -1 + 4m}  q^{\binom{d+1}{2} - \binom{\alpha+1}{2} - \binom{\beta+1}{2}}  =
% %p^{-2 \beta^2  - 3 \beta -1 + 4m} q^{\frac{1}{2}(d^2 + d - (3d-\beta+1)(3d-\beta) - (\beta+1)(\beta))} =
% %p^{-2 \beta^2  - 3 \beta -1 + 8d+4} q^{\frac{1}{2}(-8 d^2 - 2d + 6 d \beta - 2 \beta^2)} =
% %p^{-2 \beta^2  - 3 \beta + 8d+3} q^{-4 d^2 - d + 3 d \beta - \beta^2)} =
% p^{\beta^2  - 3 \beta + 11d + 12 d^2 - 9 d \beta + 3}. \end{equation}
% Recalling that $p^{6m} = p^{12d+6} = 1$, we see that $p^{12 d^2 + 11d + 3} = p^{-7d-3} p^{(12d+6)(d + 1)} = p^{-7d-3}$. Altogther this gives
% \begin{equation} \xi_m(a,i) = (-1)^{\beta+1 + d} m^2 p^{\beta^2  - 3 \beta - 9 d \beta -7d-3}  {m-1-\bottom \brack \alpha - \bottom}\end{equation}
% as desired.

Finally, when $a$ and $b$ are even, one can compute that
\begin{equation} \lambda_1 \lambda_2 \lambda_3 \lambda_4 = \pm p^{-2 \beta^2 -5\beta -3}. \end{equation}
Altogether we derive that
\begin{equation} \xi_m = (-1)^{d+\beta} m^2 q^{\binom{d}{2} - \binom{\alpha+1}{2} - \binom{\beta+1}{2}} p^{4m} {m-1-\bottom \brack \beta - \bottom} p^{-2 \beta^2 -5\beta -3}, \end{equation}
The difference between $q^{\binom{d+1}{2}}$ and $q^{\binom{d}{2}}$ is an additional factor of $q^{-d} = p^{3d}$.